\documentclass[leqno, myheadings, twoside]{amsart}

\usepackage{amsmath, amsthm, amssymb, amscd, amsxtra,graphicx}
\usepackage{latexsym, amsfonts}
\usepackage{easybmat}
\usepackage{etex}
\usepackage{url}
\usepackage{texdraw}
\usepackage{epsfig}
\usepackage{tikz}
\allowdisplaybreaks[4]
\usepackage{xcolor, soul}

\usepackage{geometry}
\usepackage{geometry}
\makeatletter \@addtoreset{equation}{section}

\makeatletter \renewcommand{\@biblabel}[1]{#1.}
\theoremstyle{remark}

\begin{document}
\setcounter{page}{1}
\title[Determining the viscosity  for incompressible fluid]{Determining the viscosity from the boundary information for incompressible fluid}

\author{Genqian Liu}
\address{School of Mathematics and Statistics, Beijing Institute of Technology, Beijing 100081, China}
\email{liugqz@bit.edu.cn}
\subjclass[2010]{}
\keywords{}

\maketitle

\date{}
\protect\footnotetext{{MSC 2020:  35R30, 76D07, 76D05, 53C21, 86A22.}
\\
{ ~~Key Words: Stokes system; Viscosity coefficient; Cauchy data; Pseudodifferential operator  } }
\maketitle ~~~\\[-15mm]

\begin{center}
{\footnotesize   School of Mathematics and Statistics, Beijing Institute of Technology, Beijing 100081, China\\
 Emails:  liugqz@bit.edu.cn \\
 }
\end{center}


\vskip 0.59 true cm

\begin{abstract}
\  For the Stokes equations in a compact connected Riemannian $n$-manifold $(\Omega,g)$ with smooth boundary $\partial \Omega$, we give an equivalent new system of elliptic equations with $(n+1)$ independent unknown functions on $\Omega$. We show that the Dirichlet-to-Neumann map ${\tilde{\Lambda}}_{\tilde{\epsilon}, \mu,g}$ associated with this new system is also equivalent to the original Dirichlet-to-Neumann map $\Lambda_{\mu,g}$ associated with the Stokes equations. We explicitly give the full symbol expression for the ${\tilde{\Lambda}}_{\tilde{\epsilon},\mu,g}$ by a method of factorization, and prove that Dirichlet-to-Neumann map ${\tilde{\Lambda}}_{\tilde{\epsilon},\mu,g}$  (or equivalently, $\Lambda_{\mu, g}$) uniquely determines viscosity $\mu$ and all tangential and normal derivatives of $\mu$ on $\partial \Omega$. In particular, combining this result, Lai-Uhlmann-Wang's theorem and Heck-Li-Wang's theorem, we completely solve a long-standing open problem that asks whether one can determine the viscosity for the Stokes equations and for the Navier-Stokes equations by boundary measurements on an arbitrary bounded domain in ${\mathbb{R}}^n$, ($n=2,3$).
\end{abstract}

\maketitle 

\vskip 1.39 true cm

\section{Introduction}

\vskip 0.45 true cm

 Let $\Omega\subset {\mathbb{R}}^n$, ($n=2,3$), be a bounded domain with smooth boundary $\partial \Omega$. Assume that $\Omega$ is filled with an incompressible fluid. Let $u = (u^1,\cdots, u^n)^t$ be the velocity vector field satisfying the stationary Stokes equations
   \begin{eqnarray} \label{20200502-1}  \left\{\!\begin{array}{ll} \mbox{div}\; \sigma_\mu (u, p) =0 \;\; & \mbox{in} \;\; \Omega,\\
   \mbox{div}\; u=0 \;\;& \mbox{on}\;\; \Omega,\end{array} \right. \end{eqnarray}
    where $\sigma_{\mu} (u, p)=2\mu \,\mbox{Def}\, (u) - pI_n$ is the ``stress tensor'' and
    $\; \mbox{Def}\, (u) = ((\nabla u) +(\nabla u)^t )/2 $ is the ``deformation tensor'',
$\mu$ is the viscosity and $p$ is the pressure. Here $I_n$ is the $n\times n$ identity matrix, $A^t$ is the transpose of a matrix (or vector) $A$.
     Physically, most fluids have positive viscosities. (Zero viscosity is observed only in superfluids that have the ability to
self-propel and travel in a way that defies the forces of gravity and surface tension, see \cite{LUW}.) Thus, we can assume that $\mu > 0$ in
$\bar\Omega$. A fluid with nonconstant viscosity is called a non-Newtonian fluid which is relatively common, such as blood, shampoo,
 custard and salt water with varying salinity. The second equation of (\ref{20200502-1}) is the incompressibility condition.
Let $\phi\in H^{\frac{3}{2}} (\partial \Omega) $ satisfy  the standard flux compatibility condition
 \begin{eqnarray} \label{20200502-2} \int_{\partial \Omega} \phi\cdot \nu \,ds=0 \end{eqnarray}
 where $\nu$ is the unit outer normal field to $\partial \Omega$. This boundary condition leads to the uniqueness of (\ref{20200502-1}), that is, there exists a unique $(u,p)\in H^2(\Omega) \times H^1(\Omega)$ ($p$ is unique up to a constant) solving (\ref{20200502-1}) and $u\big|_{\partial \Omega} =\phi$ (see \cite{LUW}, \cite{HLW} or \cite{LiW}). Throughout this paper, we always take $\int_{\Omega} p\, dV=0$, where $dV$ denotes the volume element in $\Omega$. Thus we can define the Cauchy data of $(u,p)$ satisfying (\ref{20200502-1})
  \begin{eqnarray*} {\mathcal{C}}_\mu = \left\{ (u\big|_{\partial \Omega}, \sigma_\mu (u, p) \nu\big|_{\partial \Omega}) \right\} \subset H^{\frac{3}{2}}(\partial \Omega) \times H^{\frac{1}{2}} (\partial \Omega).\end{eqnarray*}
  In physical sense, $\sigma_\mu (u,p) \nu\big|_{\partial \Omega}$ represents the Cauchy force acting on $\partial \Omega$. We also call $\sigma_\mu(u,p)\nu$ the Neumann boundary condition for the Stokes equations; so we can define the Dirichlet-to-Neumann map $\Lambda_{\mu}: H^{\frac{3}{2}}(\partial \Omega) \to H^{\frac{1}{2}}(\partial \Omega)$, associated
  with the Stokes equations, given by:
  \begin{eqnarray} \label{20200502-3} \Lambda_{\mu} u= \sigma_\mu(u,p)\nu \;\;\; \mbox{for any}\;\; \int_{\partial \Omega} u\cdot \nu \, ds =0 \;\;\mbox{and}\;\; \int_{\Omega}p\, dV=0,\end{eqnarray}
   where $dV$ is the volume element in $\Omega$. It has been a very interesting and challenging open problem (see, for example, \cite{HLW}, \cite{LUW}, \cite{Kir} or \cite{IY}) that whether one can determine $\mu$ from the knowledge of ${\mathcal{C}}_\mu$ (or equivalently, from the Dirichlet-to-Neumann map $\Lambda_{\mu}$)?

Such problems have been studied for a long time since the publication of the paper
by Calder\'{o}n \cite{Cald} in 1980, in particular for the identification of the scalar parameter $a>0$ in
operators of the form $v\mapsto -\mbox{div}\,(a\nabla v)$. In Calder\'{o}n's problem, $v$ represents an electric potential
and one assumes that the Poincar\'{e}-Steklov operator (also called the Dirichlet-to-Neumann map)
$\Lambda_a : H^{\frac{1}{2}} (\partial \Omega) \to H^{-\frac{1}{2}} (\partial\Omega)$ is known ($\Lambda_a$ is defined by $\Lambda_a(\psi):= a \frac{\partial v}{\partial \nu} $ where $\mbox{div}\, (a\nabla v) = 0$ in $\Omega$ and $v =\psi$ on $\partial \Omega$). The interested reader is referred to the review by Uhlmann \cite{Uhl1} for key historical remarks on this matter and to the pioneering works by Kohn and Vogelius \cite{KV} and Sylvester and Uhlmann \cite{SU1} for early results on this theory. We also refer the reader to \cite{CZ} and \cite{Pich} for the isotropic electromagnetic parameter problem, and further to \cite{NU1}, \cite{NU2} or \cite{Isak} for the isotropic elastic parameter problem.　

 Since the Stokes equations and Navier-Stokes equations play a very important role in fluid mechanics and physics, the viscosity determination problem by boundary measurements has a high attention in the field of inverse problems (see \cite{Isak}, \cite{ACFKO}, \cite{IY}, \cite{{LiW}}, \cite{HLW} and \cite{LUW}). Some great breakthroughs have been made for the above open problem:　　

  \vskip 0.29 true cm

  \noindent{\bf Theorem 1.1 (Lai-Uhlmann-Wang \cite{LUW}).}  {\it $\,$Let $\Omega$ be a simply connected bounded domain in $R^2$ with smooth
boundary $\partial \Omega$. Suppose that $\mu_1$ and $\mu_2$ are two viscosity functions for the Stokes equations. Assume that $\mu_j \in  C^3(\bar \Omega)$  and $\mu_j> 0$  with \begin{eqnarray}\label{20200504-1}\frac{\partial^{|K|}\mu_1(x)}{\partial x_K} = \frac{\partial^{|K|}\mu_2(x)}{\partial x_K}, \;\,  \;\; \forall  x\in \partial \Omega, \;\,  |K|\le 1.\end{eqnarray}
Let ${\mathcal{C}}_{\mu_1} $ and ${\mathcal{C}}_{\mu_2} $
be the Cauchy data associated with $\mu_1$ and $\mu_2$, respectively. If
 ${\mathcal{C}}_{\mu_1}={\mathcal{C}}_{\mu_2} $, then $\mu_1=\mu_2$ in $\Omega$.}

   \vskip 0.29 true cm

\noindent{\bf Theorem 1.2 (Heck-Li-Wang \cite{HLW}).}  {\it $\,$Let $\Omega$ be a bounded domain in $R^3$ with smooth
boundary $\partial \Omega$. Assume that $\mu_1$ and $\mu_2$ are two viscosity functions satisfying $\mu_1, \mu_2\in C^{n_0} (\bar \Omega)$ for $n_0\ge 8$ and
\begin{eqnarray} \label{20200502-7} \frac{\partial^{|K|} \mu_1(x)}{\partial x_K} =\frac{\partial^{|K|}\mu_2(x)}{\partial x_K}, \quad \forall x\in \partial \Omega,\; \,|K|\le 1.\end{eqnarray}
 Let ${\mathcal{C}}_{\mu_1}$ and ${\mathcal{C}}_{\mu_2}$  be the Cauchy data associated with $\mu_1$ and $\mu_2$, respectively. If ${\mathcal{C}}_{\mu_1}={\mathcal{C}}_{\mu_2}$, then $\mu_1=\mu_2$.}

   \vskip 0.29 true cm

  Heck, Li and Wang \cite{HLW} further proved that if $\Omega\subset {\mathbb{R}}^3$ is convex with boundary $\partial \Omega$ having nonvanishing Gauss curvature, and if $\mu_1(x)$ and $\mu_2(x)$ belong to $C^8(\bar \Omega)$ and ${\mathcal{C}}_{\mu_1}={\mathcal{C}}_{\mu_2}$, then $\mu_1=\mu_2$ and $\nabla\mu_1(x)\cdot \nu =\nabla\mu_2(x)\cdot \nu$ for all $x\in \partial \Omega$. In other words, Heck, Li and Wang actually proved that for a bounded convex domain $\Omega\subset {\mathbb{R}}^3$ with $\partial \Omega$ having nonvanishing Gauss curvature, if $\mu(x)\in C^8(\bar \Omega)$, then
the Cauchy data ${\mathcal{C}}_\mu$ uniquely determines the viscosity function $\mu$ in $\Omega$.

It remains to ask whether one can remove the convex assumption in three-dimensional case or directly prove (\ref{20200504-1}) from ${\mathcal{C}}_{\mu_1}={\mathcal{C}}_{\mu_2}$ for two-dimensional case?

\vskip 0.01 true cm

In this paper, by establishing an equivalent new system of elliptic equations for the Stokes equations in a Riemannian manifold and by factoring this new system into a product of two operators of order $1$, we get a pseudodifferential operator ${\tilde{\Lambda}}_{\tilde{\epsilon}, \mu,g}$, which is equivalent to the Dirichlet-to-Neumann map $\Lambda_{\mu, g}$ associated with the Stokes equations. Furthermore, by calculating the full symbol of the operator ${\tilde{\Lambda}}_{\tilde{\epsilon},\mu, g}$ and by analysing its homogenous symbols of degree $1$ and $0$, we show that ${\tilde{\Lambda}}_{\tilde{\epsilon},\mu,g}$  uniquely determines $\mu$ and its all  derivatives of order $1$ on $\partial \Omega$.

\vskip 0.29 true cm

 \noindent{\bf Theorem 1.3.} \ {\it  Let $(\Omega,g)$ be a compact Riemannian $n$-manifold with smooth boundary $\partial \Omega$, $(n=2,3)$. Assume that $\mu_1, \mu_2\in C^8(\bar \Omega)$ for $n=3$ (respectively, $\mu_1,\mu_2\in C^3(\bar \Omega)$ for $n=2$). If ${\mathcal{C}}_{\mu_1}={\mathcal{C}}_{\mu_2}$, then $\mu_1=\mu_2$ and $\frac{\partial^{|K|} \mu_1}{\partial x_K} =
 \frac{\partial^{|K|} \mu_1}{\partial x_K}$  for all  $x\in \partial \Omega$ and all $|K|\le 1$.}

\vskip 0.19 true cm

In particular, when $n=2$ or $n=3$ and the metric $g$ of $\Omega$ is the standard Euclidean metric (i.e., $g_{jk}=\delta_{jk}$), our Theorem 1.3 implies that the Cauchy data ${\mathcal{C}}_\mu$ uniquely determines $\mu$ and its all derivatives of order $1$ on $\partial \Omega$. Combining this result,  Lai-Uhlmann-Wang theorem  and Heck-Li-Wang theorem,  we have the following global uniqueness result. Note that the following theorem also holds for the Navier-Stokes equations:

\vskip 0.29 true cm

 \noindent{\bf Theorem 1.4.} \ {\it Let $\Omega\subset {\mathbb{R}}^n$, $(n=2,3)$, be a bounded domain with smooth boundary $\partial \Omega\,$ $(\Omega$ is required to be simply connected when $n=2)$. Assume that $\mu_1(x)$ and $\mu_2(x)$ are two viscosity functions satisfying $\mu_1, \mu_2 \in C^8(\bar \Omega)$ for $n=3$ (respectively, $\mu_1, \mu_2 \in C^3(\bar \Omega)$ for $n=2$). If ${\mathcal{C}}_{\mu_1}={\mathcal{C}}_{\mu_2}$, then $\mu_1=\mu_2$ in $\Omega$.}

\vskip 0.19 true cm

\noindent Therefore, the global identifiability problem for the viscosity in a bounded three-dimensional (or two-dimensional) incompressible fluid by boundary measurement is completely answered.

\vskip 0.19 true cm

The main ideas of this paper is as follows. The Dirichlet-to-Neumann map associated with the Stokes equations is a pseudodifferential operator defined on the boundary. In order to study this kind of operator, an effective method is to explicitly calculate its full symbol (see, for example, \cite{Ho3} or \cite{Gr}). However, explicit symbol calculation must apply the knowledge of Riemannian manifold (by flatting the boundary and inducing a Riemannian metric in a neighborhood of the boundary). Thus, we first give the local expression of the stationary Stokes equation in a Riemannian manifold $(\Omega,g)$ in terms of vector field by applying a result of \cite{Liu1}, in which the author of this paper gave an exact expression for the elastic equations. We then propose a key transformation with a fixed constant $\rho$, from which the Stokes equations is transformed into an equivalent new system of elliptic partial differential equations (see second 2). Since this new system is a linear, second order elliptic matrix-valued equation with $n+1$ independent unknown functions, in local normal coordinates we can rewrite it as  \begin{eqnarray*}   \left\{ \bigg[\frac{\partial^2}{\partial x_n^2}I_{n+1}\bigg]+B\bigg[\frac{\partial}{\partial x_n}I_{n+1}\bigg] +{C} \right\}\begin{bmatrix}w^1\\
  \vdots\\
 w^n\\
 f\end{bmatrix}=0,\end{eqnarray*}
 where $B$ and $C$ are differential operators of order $1$ and order $2$, respectively.
    So we will look for the factorization $\frac{\partial^2}{\partial x_n^2}I_{n+1}+B\left[\frac{\partial}{\partial x_n}I_{n+1}\right] +{C}= \left(\left[\frac{\partial }{\partial x_n}I_{n+1}\right] +B+Q\right)\left(\left[\frac{\partial }{\partial x_n}I_{n+1}\right] +Q\right)$, where operator $Q$ will be determined late (in fact, $\big((\frac{\partial}{\partial x_n} I_{n+1})(w,f)^t\big)\big|_{\partial \Omega}= -(Q\,(w,f)^t)\big|_{\partial \Omega}$ modulo a soothing operator, and the symbol of $Q$ has the form $\sum_{j\le 1} q_j(x, \xi')$). Because the matrix $B$ is a differential operator of order one, we will encounter two major difficulties:

      i) \ \  How to solve the unknown $q_1$ from the following matrix equation?
       \begin{eqnarray} \label{19.4.30-1} q_1^2+b_1 q_1 -c_2=0,\end{eqnarray}
       where $q_1$, $b_1$ and $c_2$ are the principal symbols of the differential (or pseudodifferential) operators $Q$, $B$ and $C$, respectively.

        Generally, the quadratic matrix equation of the form (\ref{19.4.30-1}) can not be exactly solved (in other words, there is not a formula of the solution represented by the coefficients of matrix equation (\ref{19.4.30-1}). Fortunately, in our setting we get the exact solution by applying a method of algebra ring theory in this paper. More precisely, by observing the coefficients of matrix equation (\ref{19.4.30-1}) we define an invariant sub-ring $\mathfrak{F}$ which is generated by coefficients matrices of equation (\ref{19.4.30-1}). This implies that the $q_1$ has a special form (see section 3), and hence by solving a linear equations (a system of coefficient equations) for the unknown constants in $q_1$, we obtain an exact solution $q_1$ (a surprise result$\,$!). This method is inspired by Galois group theory to solve the polynomial equation (see, for example, \cite{Art} or \cite{HME}) and was recently established by the author of this paper in \cite{Liu1} for solving an elastic inverse problem.

       ii) \ \   How to solve Sylvester's equation: $(q_1-b_1) q_{j-1} +q_{j-1} q_1=E_j\,$? where $q_{j-1}$ ($j\le 1$) are the remain symbols of $q$ (here $q\sim \sum_{j\le 1} q_j$),  and $E_j$ can be seen in section 3.

In mathematics (more precisely, in the field of control theory), a Sylvester equation is a matrix equation of the form (see \cite{Syl} and \cite{BaS}):
\begin{eqnarray}\label{19.6.22-1}  LX+XM=V. \end{eqnarray}
Then given matrices $L$, $M$, and $V$, the problem is to find the possible matrices $X$ that obey this equation. A celebrated result (see \cite{BaS} or \cite{BhR}) states that Sylvester's equation (\ref{19.6.22-1}) has a unique solution $X$ for all $V$ if and only if $L$ and $-M$ have no common eigenvalues. Generally, it is a quite difficult or impossible task to obtain an explicit solution of the Sylvester's equation. However, by putting Sylvester's equation into an equivalent $(n+1)^2\times (n+1)^2$-matrix equation and by applying our invariants sub-ring method mentioned above, we can get the inverse $U^{-1}$ of
$(n+1)^2\times (n+1)^2$-matrix $U:= (I_{n+1}\otimes L) + (M^t \otimes I_{n+1})$ and further obtain the exact solution of $q_{j-1}$, $(j\le 1)$ (another surprise result!)  Therefore, the equivalent Dirichlet-to-Neumann map ${\tilde{\Lambda}}_{\tilde{\epsilon},\mu,g}$ on $\partial \Omega$ is obtained.
 By analysing the full symbol of ${\tilde{\Lambda}}_{\tilde{\epsilon},\mu,  g}$, we find that ${\tilde{\Lambda}}_{\tilde{\epsilon},\mu,g}$ (or equivalently, $\Lambda_{\mu, g}$) uniquely determines the viscosity function $\mu$ and all its tangent and normal derivatives of order $1$ at every point $x_0\in \partial \Omega$. By using some simple properties in Riemannian geometry, we get that $\Lambda_{\mu,g}$ uniquely determines $\mu$ and $\nabla_g \mu$ on $\partial \Omega$,  and Theorem 1.3 is proved.

The paper is organized as follows. We show the equivalence of the Stokes equations and a new introduced system of elliptic partial differential equations in Riemannian manifold in Section 2. In
Section 3, we derive a  pseudodifferential operator from the new system and show that it is equivalent to the original Dirichlet-to-Neumann map $\Lambda_{\mu,g}$ associated with the Stokes equations. We calculate the full symbol of ${\tilde{\Lambda}}_{\tilde{\epsilon},\mu,g}$ and then prove that the Cauchy data of the ${\tilde{\Lambda}}_{\tilde{\epsilon},\mu,g}$  uniquely determines the
viscosity $\mu$ and all its derivatives of order $1$ on the boundary. In Section 4, we study the same inverse problem for the Navier-Stokes
equations.

\vskip 1.49 true cm

\section{Stokes equations and its equivalent system on a Riemannian manifold}

\vskip 0.45 true cm

  Let $\Omega$ be an $n$-dimensional Riemannian  manifold with smooth  boundary $\partial \Omega$, and let $\Omega$ be equipped with a smooth metric tensor $g$ (also denoted by $\langle\;, \; \rangle$, which is a smoothly varying  inner product on the tangent space. We will often denote the metric $g$ (and tensors in general) by its components $g_{jk}$). Denote by $[g^{jk}]_{n\times n}$ the inverse of the matrix $[g_{jk}]_{n\times n}$ and set $|g|:= \mbox{det}\, [g_{jk}]_{n\times n}$. In
particular, $d\mbox{V}$, the volume element in $\Omega$ is locally given by $d\mbox{V} = \sqrt{|g|}\, dx_1\cdots dx_n$. By $T\Omega$ and $T^*\Omega$ we denote, respectively, the tangent and cotangent bundle on $\Omega$. We shall also denote by $T\Omega$ global ($C^\infty$) sections in $T\Omega$ (i.e., $T\Omega \equiv C^\infty (\Omega, T\Omega)$);
similarly, $T^*\Omega \equiv C^\infty(\Omega, T^*\Omega)$.
Throughout this paper, we will use the Einstein summation convention: if the same index name appears exactly twice in any monomial term, once as an upper index and once as a lower index, that term is understood to be summed over all possible values of that index, generally from $1$ to the dimension $n$ of the space in question unless otherwise indicated. Let $\{ x_j\}$ be local coordinates in a neighborhood $\mathcal{O}$ of some point of $\Omega$. In $\mathcal{O}$ the vector fields $\{\frac{\partial}{\partial x_j} \} $ form a local basis for $T\Omega$. A vector field $X$ in $T\Omega$  will be denoted as $X=X^j \frac{\partial }{\partial x_j}$, where $X^j$ is called the $j$th component of $X$ in given coordinates. Recall first that
\begin{eqnarray} \label{9-2.1} \mbox{div} \,{X} := \frac{1}{\sqrt{|g|}}\, \frac{\partial(\sqrt{|g|} \,{{X}}^j)}{\partial x_j}\quad \, \mbox{if}\;\; {{X}}= {{X}}^j \frac{\partial}{\partial x_j}\in T\Omega,\end{eqnarray}
and \begin{eqnarray} \label{9-2.2} \nabla_g v =  \bigg(g^{jk} \frac{\partial v}{\partial x_k}\bigg)\frac{\partial}{\partial x_j}\quad \, \mbox{if}\;\; v\in C^\infty(\Omega), \end{eqnarray}
are, respectively, the usual divergence and gradient operators. Accordingly, the Laplace-Beltrami
operator $\Delta_g$ is just given by
 \begin{eqnarray} \label{18/12/22-2} \Delta_g:= \mbox{div}\, \nabla_g  = \frac{1}{\sqrt{|g|}}  \frac{\partial}{\partial x_j} \bigg(\sqrt{g}\,g^{jk} \frac{\partial}{\partial x_k}\bigg). \end{eqnarray}
Next, let $\nabla$ be the associated Levi-Civita connection. For each ${X} \in T\Omega$, $\nabla X$ is the tensor of type $(0,2)$ defined by
\begin{eqnarray} \label{9-2.5} (\nabla {X})({Y},{Z}):= \langle \nabla_{{Z}} {X}, {Y}\rangle, \quad \; \forall\, {Y}, {Z}\in T\Omega. \end{eqnarray}
It is well-known that in a local coordinate system with the naturally associated frame field on the tangent bundle,
\begin{eqnarray*} \nabla_{\frac{\partial}{\partial x_k}} X =  \big(\frac{\partial X^j}{\partial x_k} + \Gamma_{lk}^j X^l  \big)\frac{\partial }{\partial x_j}\quad \; \mbox{for}\;\; X= X^j \frac{\partial}{\partial x_j}, \end{eqnarray*}
 where $\Gamma_{lk}^j= \frac{1}{2} g^{jm} \big( \frac{\partial g_{km}}{\partial x_l} +\frac{\partial g_{lm}}{\partial x_k} -\frac{\partial g_{lk}}{\partial x_m}\big)$ are the Christoffel symbols associated with the metric $g$ (see, for example, \cite{Ta2}). If we denote \begin{eqnarray*} {X^j}_{;k}= \frac{\partial X^j}{\partial x_k} + \Gamma_{lk}^j X^l,\end{eqnarray*}
 then     \begin{eqnarray*} \label{18/10/29}  \nabla_Y X = Y^k {X^j}_{;k} \,\frac{\partial}{\partial x_j} \;\; \mbox{for}\,\;X=X^j \frac{\partial }{\partial x_j}, \;\, Y= Y^{k} \frac{\partial}{\partial x_k}.\end{eqnarray*}
   The symmetric part of $\nabla {X}$ is $\mbox{Def}\, {X}$, the deformation of ${X}$, i.e.,
\begin{eqnarray} (\mbox{Def}\; {X})({Y},{Z}) =\frac{1}{2} \{ \langle \nabla_{{Y}} {X}, {Z}\rangle +\langle \nabla_{{Z}} {X}, {Y}\rangle \}, \quad \, \forall\, {Y}, {Z}\in T\Omega\end{eqnarray}
(whereas the antisymmetric part of $\nabla\, {X}$ is simply $d{X}$, i.e.,
\begin{eqnarray*} d{X}({Y}, {Z}) = \frac{1}{2} \{ \langle \nabla_{{Y}} {X}, {Z}\rangle -\langle \nabla_{{Z}} {X}, {Y}\rangle \}, \quad \, \forall \,{Y}, {Z}\in T\Omega.)\end{eqnarray*}
   The Riemann curvature tensor $\mathcal{R}$ of $\Omega$ is given by
\begin{eqnarray} \label{9-2.10} \mathcal{R}({X}, {Y}){Z} = [\nabla_{{X}}, \nabla_{{Y}}]{Z} -\nabla_{[{X},{Y}]} {Z}, \quad \, \forall\, {X}, {Y}, {Z} \in T\Omega, \end{eqnarray}
where $[{X}, {Y}] := {X}{Y} - {Y} {X}$ is the usual commutator bracket. It is convenient to change this
into a $(0, 4)$-tensor by setting
\begin{eqnarray*} \mathcal{R}({X}, {Y}, {Z}, {W}) := \langle \mathcal{R}({X}, {Y}){Z}, {W}\rangle, \quad \; \forall\, {X}, {Y}, {Z}, {W} \in T\Omega.
\end{eqnarray*}
In other words, in a local coordinate system such as that discussed above,
 \begin{eqnarray*} &&R_{jklm}= \bigg\langle \mathcal{R}\left(\frac{\partial}{\partial x_l}, \frac{\partial}{\partial x_m}\right) \frac{\partial}{\partial x_k}, \frac{\partial}{\partial x_j}\bigg\rangle.\end{eqnarray*}
The Ricci curvature $\mbox{Ric}$ on $\Omega$ is a $(0, 2)$-tensor defined as a contraction of $\mathcal{R}$:
\begin{eqnarray*} \mbox{Ric} ({X}, {Y} ):=  \bigg\langle \mathcal{R}\bigg(\frac{\partial }{\partial x_j}, {Y}\bigg) {X}, \frac{\partial }{\partial x_j}\bigg\rangle = \bigg\langle \mathcal{R}\bigg({Y},\frac{\partial }{\partial x_j}\bigg)\frac{\partial }{\partial x_j}, {X}\bigg\rangle, \quad \forall\, {X}, {Y} \in T\Omega.\end{eqnarray*}
  That is,  \begin{eqnarray}\label{19.10.3-2}   R_{jk} =  R^l_{jlk}= g^{lm}R_{ljmk}.\end{eqnarray}
 Note that \begin{eqnarray} \label{19.6.25-1} R^j_{klm}=\frac{\partial \Gamma^j_{km}}{\partial x_l}- \frac{\partial \Gamma^j_{kl}}{\partial x_m} +\Gamma^j_{sl} \Gamma^s_{km}- \Gamma^j_{sm}\Gamma^s_{kl}.\end{eqnarray}

\vskip 0.25 true cm

 Now, assume that the Riemannian manifold $\Omega$ is filled with an incompressible fluid. Let $u =  u^k \frac{\partial }{\partial x_k}\in T\Omega$ be the velocity vector field satisfying the stationary Stokes equations \begin{eqnarray} \label{2020.4.18-1}\left\{\begin{array}{ll}  \mbox{div}\; \sigma_{\mu} (u,p)=0 \quad &\mbox{in}\;\; \Omega,\\
  \mbox{div}\; u=0\quad &\mbox{in} \;\; \Omega,
  \end{array} \right.\end{eqnarray}
   where $\sigma_{\mu} (u, p)=2\mu \,\mbox{Def}\, (u) - pg$.
    Let us note that the deformation tensor is a symmetric tensor field of type $(0,2)$ defined by
   \begin{eqnarray*} (\mbox{Def}\; u) (Y,Z) = \frac{ 1}{2} \left( \langle \nabla_Y u,  Z\rangle + \langle \nabla_Z u, Y\rangle\right), \,\; \;\; \forall\, Y, Z\in T\Omega;\end{eqnarray*}
   in coordinate notation, $(\mbox{Def}\; u)_{jk} =\frac{1}{2} (u_{j;k}+u_{k;j})$, where $u_{j;k}= \frac{\partial u_j}{\partial x_k}-
   \Gamma_{jk}^l u_l$. We have $\mbox{Def}:  C^\infty (\bar \Omega, T) \to C^\infty (\bar \Omega, S^2 T^* )$ (see p.$\,$464 of \cite{Ta1}). This tensor was introduced in Chap.$\,$2, $\S 3$, cf (3.35) of \cite{Ta1}.
 The adjoint  ${\mbox{Def}}^*$ of $\mbox{Def}$ is defined in local coordinates by \begin{eqnarray} \label{20200507-2} ({\mbox{Def}}^* {w})^j=- {{w}}^{jk}_{\;\;\;\,;k}\end{eqnarray}  for each symmetric tensor field ${w}:=w_{jk}$ of type $(0,2)$.
 In particular, if $\nu\in T\Omega$ is the outward unit normal to $\partial \Omega\hookrightarrow \Omega$, then the integration by parts formula (see formula (2.16) of \cite{DMM})
\begin{eqnarray} \label{18/113-4} \int_{\Omega} \langle \mbox{Def}\; u, w\rangle dV= \int_\Omega \langle u, \mbox{Def}^* w\rangle dV+ \int_{\partial \Omega} w( \nu, u) \,ds\end{eqnarray}
holds for any $u\in T\Omega$ and any symmetric tensor field $w$ of type $(0,2)$.
   Setting $S:= 2\, \mbox{Def}\;u$, we have
   \begin{eqnarray*} \mbox{div}\; \sigma_\mu (u,p) = \mbox{div}\,( \mu S) -\mbox{div}\,  (pg).\end{eqnarray*}
   It follows from p.$\,$562 of \cite{Ta2} that   \begin{eqnarray} \label{20200507-1} \mbox{div}\, w= -\mbox{Def}^*w\end{eqnarray}
    for any (0,2) type tensor $w$. According to the definition of divergence, it is easy to verify that for any scalar function $\psi$ and any vector field $X=X^j\frac{\partial }{\partial x_j}$,
    \begin{eqnarray} \label{20200511-2}&&　\mbox{div}\,(\psi X) = \psi \,\mbox{div}\; X + \langle \nabla_g \psi, X\rangle\\
    && \qquad \quad \quad\, =\psi \,\mbox{div}\; X + \frac{\partial \psi}{\partial x_k} X^k.\nonumber\end{eqnarray}
   Thus, in index notation (see p.$\,$562 of \cite{Ta3})
   \begin{eqnarray*} S^{jk} = u^{j;k}+ u^{k;j}, \end{eqnarray*}
   where $u^{j;k}=  g^{jm}g^{kl}u_{m;l}$,  and the $j$th component $(\mbox{div}(\mu S))^j$
   of vector field $\mbox{div}(\mu S)$ is given by
\begin{eqnarray}\label{20200507-4}  (\mu S^{jk})_{; k}=\!\!\!\!\! \!\!&&\left(\mu ( u^{j;k} +u^{k;j})\right)_{;k}
\\
=\!\!\!\!\! \!\! &&\mu u^{j;k}_{\;\;\;\;\,;k} +\mu u^{k;j}_{\;\;\;\;\,;k}+ S^{jk} \frac{\partial \mu}{\partial x_k}.\nonumber\end{eqnarray}
The first term in the last equality of (\ref{20200507-4}) is $- \mu (\nabla^{*} \nabla u)^j$; and the second term can be written as (see, (3.16) on p.$\,$554 or p.$\,$562 in \cite{Ta3})
\begin{eqnarray} \mu \left( u_{\;\;;k}^{k\;\;\; ;j} + R_{\;lk}^{k\;\;j} u^l\right) =\mu\big( \nabla_g (\mbox{div}\; u) +\mbox{Ric}\, (u)\big)^j,\end{eqnarray}
where $\big(\mbox{Ric}\, (u)\big)^j:=  R^j_l u^l=R_{\;lk}^{k\;\;j} u^l$.
 Hence, the $j$th component $(\mbox{div}\,(\mu S))^j$ is just  \begin{eqnarray}
\mu \big(\!- \nabla^* \nabla u + \nabla_g( \mbox{div}\; u) +\mbox{Ric}\, (u) \big)^j+S^{jk}  \frac{\partial \mu}{\partial x_k},\end{eqnarray}
Thus, as long as $\mbox{div}\; u=0$, we have
\begin{eqnarray*}(\mbox{div}\, (\mu S))^j= \mu \big(\!- \nabla^* \nabla u  +\mbox{Ric}\, (u) \big)^j+   S^{jk}\frac{\partial \mu}{\partial x_k}.\end{eqnarray*}
Similarly, the $j$th component $(\mbox{div}\;\!(pg))^j$ of vector field $\mbox{div}\;\!(pg)$ is  \begin{eqnarray*} &&(pg^{jk})_{;k}=pg^{jk}_{\,\;\;\;;k}+g^{jk} \frac{\partial p}{\partial x_k}  =g^{jk} \frac{\partial p}{\partial x_k} \end{eqnarray*}  because $g$ is a tensor of type $(0,2)$ and $g^{jk}_{\;\;\,\;;k}=0$.
It follows that
\begin{eqnarray} \label{20200510-1} (\mbox{div}\; \sigma_\mu (u,p))^j = \mu \left( -\nabla^*\nabla u +\mbox{Ric}\, (u) \right)^j + S^{jk}\,\frac{\partial u}{\partial x_k} -g^{jk} \frac{\partial p}{\partial x_k}=0,\end{eqnarray}
or equivalently,
\begin{eqnarray} \label{2020.4.18-4} \mbox{div}\; \sigma_\mu (u,p) =\mu  \left( -\nabla^*\nabla u +\mbox{Ric}\, (u) \right) + S^{jk}\,\frac{\partial u}{\partial x_k}\,\frac{\partial }{\partial x_j} -\nabla_g p=0,\end{eqnarray}
provided that $\mbox{div}\; u=0$ in $\Omega$.
 For a Riemannian manifold $\Omega$, let $ \phi\in \big( H^{\frac{3}{2}} (\partial \Omega) \big)^n$ satisfy $\int_{\partial \Omega} \langle \phi, \nu\rangle ds =0$, then there exists a unique $(u,p) \in H^2(\Omega) \times H^1(\Omega)$ ($p$ is the unique up to a constant) solve (\ref{2020.4.18-1}) and $u\big|_{\partial \Omega} =\phi$ (see, for example, A of Chapter 17 in \cite{Ta3}). So we can naturally define the Cauchy data of $(u,p)$ satisfying (\ref{2020.4.18-1}) with $\int_\Omega p\, dV=0$:
 \begin{eqnarray} C_\mu =\left\{ (u\big|_{\partial \Omega}, \sigma_\mu (u, p) \nu \big|_{\partial\Omega} ) \right\} \subset
 H^{\frac{3}{2}} (\partial \Omega) \times H^{\frac{1}{2}}(\partial \Omega), \end{eqnarray}
 where  \begin{eqnarray*}  \sigma_\mu (u, p) \nu\big|_{\partial \Omega} = \left(2\mu \,\mbox{Def}\, u -p g\right) \nu\big|_{\partial \Omega}  \end{eqnarray*}
 is the Cauchy force acting on $\partial \Omega$ (or Neumann boundary condition for the Stokes equations (\ref{2020.4.18-4})). Here, we identity $(2\mu\,\mbox{Def}\; {u} -pg){\nu}$ with the vector field uniquely determined by the requirement that $\langle (2\mu\,\mbox{Def}\; {u}-pg){\nu}, {X}\rangle = (2\mu\,\mbox{Def}\; {u}-pg) ({\nu}, {X})$ for each ${X}\in T\Omega$.
 In (\ref{2020.4.18-4}), $-\nabla^*\nabla u$ can be written as (see \cite{Liu1} or \cite{Liu2})
 \begin{eqnarray}\label{2020.4.18-5} && -\nabla^*\nabla u = \left\{ \Delta_g u^j +2  g^{kl} \Gamma_{sk}^j \frac{\partial u^s}{\partial x_l}
 + \Big( g^{kl} \frac{\partial \Gamma_{sl}^j}{\partial x_k} + g^{kl} \Gamma_{hl}^j \Gamma_{sk}^h -
g^{kl}\Gamma_{sh}^j \Gamma_{kl}^h \Big) u^s\right\} \frac{\partial}{\partial x_j}.\nonumber\end{eqnarray}

\vskip 0.25 true cm

\noindent We need the following:

\vskip 0.25 true cm

\noindent{\bf Lemma 2.1.} \ {\it Let $(\Omega,g)$ be a smooth Riemannian manifold. Then, for any function $f\in C^3(\Omega)$, the following relation holds}:
\begin{eqnarray} \label{20200522-1} \Delta_g(f^{;j})= (\Delta_g f)^{;j} + R_l^{\;j} \, f^{;l},\end{eqnarray}
i.e.,
\begin{eqnarray*} \label{200419} \Delta_g( (\nabla_g f)^j)= (\nabla_g (\Delta_g f))^j + R_l^{\;j} \, (\nabla_g f)^l,\end{eqnarray*}
where  $(\nabla_g f)^j:=g^{jl} \frac{\partial f}{\partial x_l} =f^{;j}$
\vskip 0.25 true cm

 \noindent  {\it Proof.} \   Because we were not able to find an exact reference to this lemma, we provide a short proof here. For a function $f\in C^3(\Omega)$, we denote $f_{;j}= \frac{\partial f}{\partial x_j}$. Since $f_{;j;k}= \frac{\partial}{\partial x_k} \left(\frac{\partial f}{\partial x_j}\right) -\Gamma_{jk}^l \frac{\partial f}{\partial x_l}$ and $\Gamma_{jk}^l =\Gamma_{kj}^l$ we have (see also (3.27) on p.$\,$148 of \cite{Ta1}) $$f_{;j;k}= f_{;k;j},$$  so $$f^{;j;k} = f^{;k;j}$$ by raising indices twice.
 This leads to  \begin{eqnarray} \label{200419-7} f^{;j;k}_{\;\;\;\;\;\;;k} = f^{;k;j}_{\;\;\;\;\;\;; k}. \end{eqnarray}
 It is well-known that (see p.554, (3.16) of \cite{Ta2}) for a vector field $X=  X^j\frac{\partial }{\partial x_j}$, one has
$$X^k_{\;\;;j;k}= X^k_{\;\;;k;j} + R^k_{\;lkj}X^l,$$
so   \begin{eqnarray} \label{200419-5} X^{k;j}_{\;\;\;\;\;;k}= X^{k\;\;\;;j}_{\;\;;k} + R^{k\;\;\;j}_{\;\;lk}X^l=  X^{k\;\;\;;j}_{\;\;;k} + R^{\;j}_{l}X^l\end{eqnarray}
by raising an index.
Replacing $X^k$ by $f^{;k}$ in (\ref{200419-5}) we get
\begin{eqnarray} \label{200419-6} f^{;k;j}_{\;\;\;\;\;\, ; k}= f^{;k\;\;\;; j}_{\;\;\;; k} +   R^{\;j}_{l}f^{;l}.\end{eqnarray}
 Combining (\ref{200419-7})  and  (\ref{200419-6}) we obtain  \begin{eqnarray*}  &&f^{;j;k}_{\;\;\;\;\;\;;k}= f^{;k\;\;\;;j}_{\;\;\; ;k}+R^{\;j}_{l} f^{;l},\end{eqnarray*}
i.e.,
 \begin{eqnarray*} \Delta_g \big(f^{;j}\big) = (\Delta_g  f)^{;j} + R^{\;j}_{l} f^{;l}.\end{eqnarray*}  \qed

 \vskip 0.29 true cm

Next, we derive a new system of elliptic equations from the stationary Stokes equations in Riemannian manifold $\Omega$. Inspired
by \cite{HLW} for the Stokes equations in ${\mathbb{R}}^3$ (or earlier for
the isotropic elastic system \cite{AITY}, \cite{ER} and \cite{Uhl3}), we set \begin{eqnarray} \label{20200513-1}   u= (\mu+\rho)^{-\frac{1}{2}} w+\mu^{-1} \nabla_g f -f \nabla_g \mu^{-1},\end{eqnarray}
i.e., \begin{eqnarray} \label{20200512-3} u^j= (\mu+\rho)^{-\frac{1}{2}} w^j +\mu^{-1} f^{;j}-f (\mu^{-1})^{;j},\;\; \; \; j=1,\cdots, n,\end{eqnarray}
where $f^{:j} = g^{jl} \frac{\partial f}{\partial x_l}$ and $ (\mu^{-1})^{;j} =  g^{jl} \frac{\partial \mu^{-1}}{ \partial x_l}$ as before,
$\rho$ is a constant which will be determined late (Note that $\rho$ plays a key role in our discussion). We will look for equations for $(w, f)$ such that $u$ solves (\ref{2020.4.18-1}).
The last equation in (\ref{2020.4.18-1}) is divergence free condition. In other words, we find by (\ref{2020.4.18-1}), (\ref{20200513-1}) and (\ref{20200511-2}) that
\begin{eqnarray} \label{2020.4.18-6} 0= \mbox{div}\; u = \mu^{-1} \Delta_g f + (\mu+\rho)^{\!-\frac{1}{2}} \,\mbox{div}\; w -( \Delta_g\mu^{-1}) f + \frac{\partial ((\mu+\rho)^{-\frac{1}{2}})}{\partial x_k} w^k.\end{eqnarray}
Recall that  \begin{eqnarray}  \label{20200512-4} &&  (\mbox{div}\, ( \mu S))^j =
\mu \Big(\!- \nabla^*\nabla u +\mbox{ Ric}\, (u)\Big)^j + S^{jk} \frac{\partial \mu}{\partial x_k}  \\
 &&  \qquad \quad \quad\;\;\;\;=  \mu \bigg(\! \Delta_g u^j + 2  g^{kl} \Gamma_{sk}^j \frac{\partial u^s}{\partial x_l}  +  \big( g^{kl} \frac{\partial \Gamma_{sl}^j}{\partial x_k} + g^{kl} \Gamma_{hl}^j \Gamma_{sk}^h -
g^{kl}\Gamma_{sh}^j \Gamma_{kl}^h \big)u^s +  R_l^{\; j} u^l\bigg) \nonumber\\
&& \qquad \quad \quad\;\;\;\; \quad  + \big(u^{j;k} +u^{k;j}\big) \frac{\partial \mu}{\partial x_k}.\nonumber\end{eqnarray}
Inserting (\ref{20200512-3}) into (\ref{20200512-4}), we get
 \begin{eqnarray} \label{200419-1} &&\;\;\;\; (\mbox{div}\, ( \mu S))^j \! =\! \mu\bigg\{ \!\Delta_g\big( (\mu+\rho)^{-\frac{1}{2} } w^j \!+\!\mu^{\!-1} f^{;j} \!-\!f (\mu^{-1})^{;j}\big) \!+\!
2 g^{kl} \Gamma_{sk}^j \frac{\partial }{\partial x_l} \big( (\mu+\rho)^{-\frac{1}{2} } w^s \!+\!\mu^{-1} f^{;s} \!-\!f (\mu^{-1})^{;s}\big) \\
&& \quad   +  \Big( g^{kl} \frac{\partial \Gamma_{sl}^j}{\partial x_k} + g^{kl} \Gamma_{hl}^j \Gamma_{sk}^h -
g^{kl}\Gamma_{sh}^j \Gamma_{kl}^h \Big) \big( (\mu+\rho)^{-\frac{1}{2} } w^s +\mu^{-1} f^{;s} -f (\mu^{-1})^{;s}\big)\nonumber\\
 && \quad  +  R_l^{\; j}\big( (\mu+\rho)^{-\frac{1}{2} } w^l +\mu^{-1} f^{;l} -f (\mu^{-1})^{;l}\big)\bigg\}+
 \bigg\{\! \big((\mu+\rho)^{-\frac{1}{2} }  w^j \big)^{;k}  + \big( \mu^{-1} f^{;j}\big)^{;k}
 - \big(f (\mu^{-1})^{;j}\big)^{;k} \nonumber\\
 &&\quad + \big((\mu+\rho)^{-\frac{1}{2} }  w^k \big)^{;j}  + \big( \mu^{-1} f^{;k}\big)^{;j}
 - \big(f (\mu^{-1})^{;k}\big)^{;j} \!\bigg\}\frac{\partial \mu}{\partial x_k}
 .\nonumber
 \end{eqnarray}
Note  that, for any $\phi,\psi\in C^2(\Omega)$,
\begin{eqnarray} \label{200419-2}   \Delta_g ( \phi\psi ) =\!\!\!\!\!\!\!\! &&  \psi (\Delta_g \phi) + 2\langle \nabla_g \phi, \nabla_g  \psi\rangle
+\phi (\Delta_g \psi) \\
       =\!\!\!\!\!\!\!\! &&  \psi (\Delta_g \phi)+ \phi (\Delta_g \psi)  + 2   \frac{\partial \phi}{\partial x_m} g^{ml} \frac{\partial \psi}{\partial x_l}. \nonumber
        \end{eqnarray}
       We then have \begin{eqnarray*}   && \big( \mbox{div}\,(\mu S)\big)^j =  \mu (\Delta_g (\mu+\rho)^{-\frac{1}{2}}) w^j +
       \mu (\mu+\rho)^{-\frac{1}{2}} \Delta_g w^j + \mu
  (\Delta_g \mu^{-1} ) f^{;j}
        +  (\Delta_g f)^{;j} + R^j_l f^{;l}  \\
  && \quad -\mu (\Delta_g f) (\mu^{-1})^{;j}  -　 \mu (\Delta_g \mu^{-1})^{;j}  f - \mu R^{\,j}_{l} ( \mu^{-1})^{;l} f \\
  && \quad + 2 \mu  \frac{ \partial ((\mu+\rho)^{-\frac{1}{2}}\!)}{\partial x_m} g^{ml}\frac{\partial w^j}{\partial x_l} +
 2 \mu  \frac{ \partial \mu^{-1}}{\partial x_m} g^{ml}\frac{\partial (f^{;j}\!)}{\partial x_l}
 -  2 \mu  \frac{ \partial f}{\partial x_m} g^{ml}\frac{\partial ((\mu^{-1})^{;j})}{\partial x_l}\\
 && \quad + \mu \bigg\{ 2  g^{kl} \Gamma_{sk}^j \frac{\partial }{\partial x_l} \Big( (\mu+\rho)^{-\frac{1}{2} } w^s +\mu^{-1} f^{;s} -f (\mu^{-1})^{;s}\Big) \nonumber \\
  && \quad    +  \bigg( g^{kl} \frac{\partial \Gamma_{sl}^j}{\partial x_k} + g^{kl} \Gamma_{hl}^j \Gamma_{sk}^h -
g^{kl}\Gamma_{sh}^j \Gamma_{kl}^h \bigg) \Big( (\mu+\rho)^{-\frac{1}{2} } w^s +\mu^{-1} f^{;s} -f (\mu^{-1})^{;s}\Big)\nonumber\\
 && \quad  + R_l^{\; j}\big( (\mu+\rho)^{-\frac{1}{2} } w^l +\mu^{-1} f^{;l} -f (\mu^{-1})^{;l}\big)\bigg\} + \Big( ( (\mu+\rho)^{- \frac{1}{2}}w^j)^{;k} +( \mu^{-1}f^{;j})^{;k} -   ( f (\mu^{-1})^{;j} )^{;k} \Big. \nonumber\\
 && \Big. \quad +  ( (\mu+\rho)^{- \frac{1}{2}}w^k)^{;j} +( \mu^{-1}f^{;k})^{;j} -
  ( f (\mu^{-1})^{;k})^{;j} \Big)\frac{\partial \mu}{\partial x_k}\nonumber\\
  &&= \mu (\mu+\rho)^{-\frac{1}{2}} \Delta_g w^j +  \Big( \Delta_g f +\mu ( \Delta_g \mu^{-1}) f + (\mu+\rho)^{-\frac{1}{2}}  \frac{\partial \mu}{\partial x_k} w^k \Big)^{;j}  + \mu  \big(\Delta_g ((\mu+\rho)^{-\frac{1}{2}}) \big)w^j +  R^j_l  f^{;l} \nonumber\\
  && \quad -2\mu  (\Delta_g  \mu^{-1} )^{;j} f -  \mu^{;j} (\Delta_g \mu^{-1}) f - \mu (\Delta_g f) (\mu^{-1})^{;j} - \mu R^j_l ( \mu^{-1})^{;l} f \nonumber\\
  &&\quad  + 2\mu  \frac{\partial ((\mu+\rho)^{-\frac{1}{2}}\!) }{\partial x_m} g^{ml} \frac{\partial w^j}{\partial x_l} +
  2\mu  \frac{\partial \mu^{-1} }{\partial x_m} g^{ml} \frac{\partial (f^{;j}\!)}{\partial x_l}-
   2\mu  \frac{\partial f }{\partial x_m} g^{ml} \frac{\partial ((\mu^{-1})^{;j}\!)}{\partial x_l}\nonumber\\
    && \quad + \mu\bigg\{ 2  g^{kl} \Gamma_{sk}^j \frac{\partial}{\partial x_l} \Big( (\mu+\rho)^{-\frac{1}{2}} w^s +\mu^{-1} f^{;s} -f (\mu^{-1})^{;s}\Big) \nonumber\\
         && \quad    +  \bigg( g^{kl} \frac{\partial \Gamma_{sl}^j}{\partial x_k} + g^{kl} \Gamma_{hl}^j \Gamma_{sk}^h -
g^{kl}\Gamma_{sh}^j \Gamma_{kl}^h \bigg) \Big( (\mu+\rho)^{-\frac{1}{2} } w^s +\mu^{-1} f^{;s} -f (\mu^{-1})^{;s}\Big)\nonumber\\
&& \quad +  R_l^j \big((\mu+\rho)^{-\frac{1}{2}} w^l +\mu^{-1} f^{;l} -f(\mu^{-1} )^{;l}\big)\bigg\} + \Big(  ( (\mu+\rho)^{- \frac{1}{2}}w^j)^{;k} +( \mu^{-1}f^{;j})^{;k} -   ( f (\mu^{-1})^{;j} )^{;k}  \Big.\nonumber\\
 && \quad \Big. +  ( \mu^{-1} f^{;k} )^{;j} - (f(\mu^{-1})^{;k})^{;j} \Big)  \frac{\partial \mu}{\partial x_k}-  (\mu+\rho)^{-\frac{1}{2}} w^k  \big(\frac{\partial \mu}{\partial x_k} \big)^{;j} \nonumber
         \end{eqnarray*}
Also, from  (\ref{2020.4.18-6}) we have  $\Delta_g f= -\mu (\mu +{\rho})^{-\frac{1}{2}}\, \mbox{div}\; w +
   \mu (\Delta_g \mu^{-1} ) f - \mu  \frac{\partial ((\mu+{\rho})^{\!-\frac{1}{2}} )}{\partial x_k} w^k $. Substitute this into the seventh term in the last equality, we obtain
\begin{eqnarray} \label{200420-10}  \!\!\! &&\!\!\!\!\!\!\! \!\!\!\!\!\! \big( \mbox{div}\,(\mu S)\big)^j = \mu (\mu+\rho)^{-\frac{1}{2}} \Delta_g w^j\! + \! \Big( \Delta_g f +\mu ( \Delta_g \mu^{-1}) f \!+ \! (\mu+\rho)^{-\frac{1}{2}} \frac{\partial \mu}{\partial x_k} w^k \Big)^{;j} \! +\! \mu  (\Delta_g (\mu+\rho)^{-\frac{1}{2}}) w^j
\!+ \! R^j_l  f^{;l} \nonumber\\
  &&\!\! -2\mu  (\Delta_g  \mu^{-1} )^{;j} f -  \mu^{;j} (\Delta_g \mu^{-1}) f - \mu \Big(\! \!-\! \mu (\mu+\rho)^{-\frac{1}{2}}\, \mbox{div}\; w +
   \mu (\Delta_g \mu^{-1} ) f - \mu  \frac{\partial (\mu+\rho)^{-\frac{1}{2}} }{\partial x_k} w^k \Big) (\mu^{-1})^{;j} \nonumber \\
   &&\!\!　　 - \mu  R^j_l ( \mu^{-1})^{;l} f \!+\! 2\mu  \frac{\partial ((\mu\!+\!\rho)^{\!-\!\frac{1}{2}}) }{\partial x_m} g^{ml} \frac{\partial w^j}{\partial x_l}  \! - \!   2\mu  \frac{\partial f }{\partial x_m} g^{ml} \frac{\partial ((\mu^{\!-1})^{;j}\!)}{\partial x_l}\!+\! \mu\bigg\{\! 2  g^{kl} \Gamma_{\!sk}^j \frac{\partial}{\partial x_l} \!\Big( (\mu\!+\!\rho)^{\!-\frac{1}{2}} w^s \! -\!f (\mu^{\!-1})^{;s}\!\Big)\nonumber \\
         &&    +  \bigg( g^{kl} \frac{\partial \Gamma_{sl}^j}{\partial x_k} + g^{kl} \Gamma_{hl}^j \Gamma_{sk}^h -
g^{kl}\Gamma_{sh}^j \Gamma_{kl}^h \bigg) \Big( (\mu+\rho)^{-\frac{1}{2} } w^s +\mu^{-1} f^{;s} -f (\mu^{-1})^{;s}\Big)\nonumber\\
&&  +  R_l^j \big((\mu+\rho)^{-\frac{1}{2}} w^l +\mu^{-1} f^{;l} -f(\mu^{-1} )^{;l}\big)\bigg\} +  \Big( ( (\mu+\rho)^{- \frac{1}{2}}w^j)^{;k}  -   ( f (\mu^{-1})^{;j} )^{;k}    - (f(\mu^{-1})^{;k})^{;j} \Big)  \frac{\partial \mu}{\partial x_k}\nonumber\\
 && -  (\mu+\rho)^{-\frac{1}{2}} \big(\frac{\partial \mu}{\partial x_k} \big)^{;j} w^k
  +\bigg\{
  2\mu  \frac{\partial \mu^{-1} }{\partial x_m} g^{ml} \frac{\partial (f^{;j}\!)}{\partial x_l}
  +2 \mu   g^{kl} \Gamma_{sk}^j \frac{\partial (\mu^{-1}
f^{;s})}{\partial x_l}+ \Big( (\mu^{-1} f^{;j})^{;k} + (\mu^{-1} f^{;k})^{;j}\Big) \frac{\partial \mu}{\partial x_k}
   \bigg\}
  . \nonumber\end{eqnarray}
But \begin{eqnarray} && 2\mu
\frac{\partial \mu^{-1}}{\partial x_m} g^{ml} \frac{\partial (f^{;j}\!)}{\partial x_l}  +2 \mu   g^{kl} \Gamma_{sk}^j \frac{\partial (\mu^{-1}
f^{;s})}{\partial x_l}+ \Big( (\mu^{-1} f^{;j})^{;k} + (\mu^{-1} f^{;k})^{;j}\Big) \frac{\partial \mu}{\partial x_k}\nonumber\\
&& \qquad  = 2\mu  \frac{\partial \mu^{-1}}{\partial x_m} g^{ml} \bigg\{\Big( \frac{\partial f^{;j}}{\partial x_l}+  \Gamma_{ls}^j f^{;s} \Big)-  \Gamma_{ls}^j f^{;s} \bigg\}    +2 \mu   g^{kl} \Gamma_{sk}^j \frac{\partial (\mu^{-1}
f^{;s})}{\partial x_l}\nonumber \\
&& \qquad \quad + \Big( (\mu^{-1} )^{;k} f^{;j} + \mu^{-1} f^{;j;k}  +(\mu^{-1})^{;j} f^{;k} +\mu^{-1} f^{;k;j} \Big) \frac{\partial \mu}{\partial x_k}\nonumber\\
&& \qquad= 2 \mu  \frac{\partial \mu^{-1}}{\partial x_m} g^{ml} f^{;j}_{\;\;\;;l} - 2\mu   \frac{\partial \mu^{-1}}{\partial x_m} g^{ml}
 \Gamma_{ls}^j f^{;s} -    2\mu^{-1}   g^{kl} \Gamma_{sk}^j \frac{\partial \mu}{\partial x_l} f^{;s} \nonumber \\
  && \qquad \quad + 2    g^{kl} \Gamma_{sk}^j   \frac{\partial f^{;s}}{\partial x_l}
  +2 \mu^{-1}  f^{;j;k}  \frac{\partial \mu}{\partial x_k} +    \Big( (\mu^{-1} )^{;k} f^{;j}   +(\mu^{-1})^{;j} f^{;k}  \Big) \frac{\partial \mu}{\partial x_k}\nonumber\\
&&\qquad  =   2    g^{kl} \Gamma_{sk}^j \frac{\partial f^{;s}}{\partial x_l} +     \Big( (\mu^{-1} )^{;k} f^{;j}   +(\mu^{-1})^{;j} f^{;k}  \Big) \frac{\partial \mu}{\partial x_k}\nonumber  \\
&&\qquad  =   2     \Gamma_{sk}^j f^{;s;k}  - 2   g^{kl} \Gamma_{sk}^j \Gamma_{lr}^s f^{;r}
+     \Big( (\mu^{-1} )^{;k} f^{;j}   +(\mu^{-1})^{;j} f^{;k}  \Big) \frac{\partial \mu}{\partial x_k}
\nonumber  \\
&&\qquad  =   2 \Gamma_{sk}^j g^{sl} g^{km} \Big( \frac{\partial^2}{\partial x_l\partial x_m} -\Gamma_{lm}^r \frac{\partial f}{\partial x_r} \Big)
- 2 g^{kl} \Gamma_{sk}^j\Gamma_{lr}^s f^{;r}
+     \Big( (\mu^{-1} )^{;k} f^{;j}   +(\mu^{-1})^{;j} f^{;k}  \Big) \frac{\partial \mu}{\partial x_k},\nonumber
\end{eqnarray}
where the second equality used $f^{;k;j}=f^{;j;k}$.
It follows that      \begin{eqnarray}   &&    \big( \mbox{div}\,(\mu S)\big)^j  = \mu (\mu+\rho)^{-\frac{1}{2}} \Delta_g w^j +   \Big(\Delta_g f +\mu ( \Delta_g \mu^{-1}) f +(\mu+\rho)^{-\frac{1}{2}} \frac{\partial \mu}{\partial x_k} w^k \Big)^{;j}   + \mu  (\Delta_g (\mu+\rho)^{-\frac{1}{2}}) w^j+ 2  R^j_l  f^{;l}\nonumber\\
  && \quad -2\mu  (\Delta_g  \mu^{-1} )^{;j} f - \mu^{;j} (\Delta_g \mu^{-1}) f
  +\mu^2 (\mu+\rho)^{-\frac{1}{2}}\, (\mu^{-1} )^{;j}(\mbox{div}\; w )  - \mu^2 (\Delta_g \mu^{-1})(\mu^{-1})^{;j}  f  \nonumber \\
   &&\quad 　 + \mu^2  \frac{\partial ((\mu+\rho)^{-\frac{1}{2}})}{\partial x_k} (\mu^{-1})^{;j} w^k 　　 - 2\mu  R^j_l (\mu^{-1})^{;l} f+ 2\mu \frac{\partial ((\mu+\rho)^{\!-\frac{1}{2}}\!) }{\partial x_m} g^{ml} \frac{\partial w^j}{\partial x_l} \nonumber\\
  &&\quad   -    2\mu
   \frac{\partial ((\mu^{-1})^{;j}\!)}{\partial x_l} g^{ml} \frac{\partial f }{\partial x_m}
    +2\mu  g^{kl}\Gamma_{sk}^j\Big((\mu+\rho)^{-\frac{1}{2}}\frac{\partial w^s}{\partial x_l}
    +\frac{\partial ((\mu+\rho)^{-\frac{1}{2}})}{\partial x_l}\, w^s - (\mu^{-1})^{;s}\frac{\partial f}{\partial x_l} - \frac{\partial ((\mu^{-1})^{;s}\!)}{\partial x_l}\, f\Big) \nonumber\\
    && \quad  + \mu  \bigg( g^{kl} \frac{\partial \Gamma_{sl}^j}{\partial x_k} + g^{kl} \Gamma_{hl}^j \Gamma_{sk}^h -
g^{kl}\Gamma_{sh}^j \Gamma_{kl}^h \bigg) \Big( (\mu+\rho)^{-\frac{1}{2} } w^s +\mu^{-1} f^{;s} -f (\mu^{-1})^{;s}\Big)\nonumber\\
&& \quad + \mu  R_l^j \,(\mu+\rho)^{-\frac{1}{2}} w^l +  \Big( ( (\mu+\rho)^{- \frac{1}{2}}w^j)^{;k} -   ( f (\mu^{-1})^{;j} )^{;k}  -   ( f (\mu^{-1})^{;k} )^{;j} \Big)  \frac{\partial \mu}{\partial x_k}\nonumber\\
&& \quad -  (\mu+\rho)^{-\frac{1}{2}} \big(\frac{\partial \mu}{\partial x_k} \big)^{;j} w^k +2
\Gamma_{sk}^j g^{sl}g^{km} \Big( \frac{\partial^2 f}{\partial x_l\partial x_m} \Big)  -2
\Gamma_{sk}^j g^{sl}g^{km}\Gamma_{lm}^r \, \frac{\partial f}{\partial x_r}
\nonumber \\
&& \quad - 2
g^{kl}\Gamma_{sk}^j \Gamma_{lr}^s f^{;r} + \frac{\partial \mu}{\partial x_k} \Big( (\mu^{-1})^{;k} f^{;j} +(\mu^{-1})^{;j} f^{;k} \Big)  \nonumber\\
 && = \mu  (\mu+\rho)^{-\frac{1}{2}} \Delta_g w^j +  \Big(\Delta_g f +\mu ( \Delta_g \mu^{-1}) f + (\mu+\rho)^{-\frac{1}{2}}  \frac{\partial \mu}{\partial x_k} w^k \Big)^{;j}
  + \mu (\Delta_g (\mu+\rho)^{-\frac{1}{2}}) w^j  + 2 R^{\;j}_l  g^{lm} \frac{\partial f}{\partial x_m} \nonumber\\
  && \quad\,- 2 \mu (\Delta_g \mu^{\!-\!1})^{;j} f　\!-  \!\mu^{;j} (\Delta_g \mu^{\!-\!1}) f\! +\! \mu^2 (\mu\!+\!\rho)^{\!-\!\frac{1}{2}}  (\mu^{\!-\!1})^{;j}  \frac{\partial w^k}{\partial x_k}\! +\!
 \mu^2  (\mu\!+\!\rho)^{\!-\!\frac{1}{2}}  (\mu^{\!-\!1})^{;j}  \Gamma_{kl}^l w^k \! -\!
   \mu^{2}(\Delta_g \mu^{\!-\!1}) (\mu^{\!-\!1})^{;j} f\nonumber\\
    &&  \quad\,　 + \mu^{2} (\mu^{-1})^{;j} \frac{\partial ((\mu+\rho)^{-\frac{1}{2}})}{\partial x_k} w^k  - 2\mu R^j_l\,( \mu^{-1})^{;l}f +
  2\mu \frac{\partial ((\mu+\rho)^{-\frac{1}{2}})}{\partial x_m} g^{ml}\frac{\partial w^j}{\partial x_l}
    　　   -2\mu g^{ml}\frac{\partial ((\mu^{-1})^{;j})}{\partial x_l} \, \frac{\partial f}{\partial x_m}  \nonumber\\
    &&  \quad\,  +2 \mu (\mu+\rho)^{-\frac{1}{2}} g^{kl}\Gamma_{sk}^j \frac{\partial w^s}{\partial x_l} +2\mu  g^{kl}\Gamma_{sk}^j
    \frac{\partial ((\mu+\rho)^{-\frac{1}{2}})}{\partial x_l} w^s 　 -2\mu g^{kl}\Gamma_{sk}^j\,  (\mu^{-1})^{;s} \frac{\partial f}
   {\partial x_l}     -   2\mu g^{kl}\Gamma_{sk}^j\, \frac{\partial ((\mu^{-1})^{;s})}{\partial x_l}\, f \nonumber\\
    &&    \quad\,　+  \mu (\mu+{\rho})^{-\frac{1}{2}} \bigg( g^{kl} \frac{\partial \Gamma_{sl}^j}{\partial x_k} + g^{kl} \Gamma_{hl}^j \Gamma_{sk}^h -
g^{kl}\Gamma_{sh}^j \Gamma_{kl}^h \bigg) w^s 　 +
  \bigg( g^{kl} \frac{\partial \Gamma_{sl}^j}{\partial x_k} + g^{kl} \Gamma_{hl}^j \Gamma_{sk}^h -
g^{kl}\Gamma_{sh}^j \Gamma_{kl}^h \bigg) g^{sr} \frac{\partial f}{\partial x_r} \nonumber\\
&&  \quad\,　 - \mu\bigg( g^{kl} \frac{\partial \Gamma_{sl}^j}{\partial x_k} + g^{kl} \Gamma_{hl}^j \Gamma_{sk}^h -
g^{kl}\Gamma_{sh}^j \Gamma_{kl}^h \bigg)(\mu^{-1})^{;s} f   +  \mu (\mu+\rho)^{-\frac{1}{2}} R_l^j w^l
 +(\mu+\rho)^{-\frac{1}{2}}
 \frac{\partial \mu}{\partial x_k} g^{kl} \Big( \frac{\partial w^j}{\partial x_l} +\Gamma_{ls}^j w^s\Big) \nonumber\\
&&  \quad\,　 +   \frac{\partial \mu}{\partial x_k} \, ((\mu+\rho)^{-\frac{1}{2}})^{;k} w^j -
   \frac{\partial \mu}{\partial x_k} \, (\mu^{-1})^{;j;k} f    -        \frac{\partial \mu}{\partial x_k} \, (\mu^{-1})^{;k;j} f
 - (\mu+\rho)^{-\frac{1}{2}} \big(\frac{\partial \mu}{\partial x_k} \big)^{;j} w^k \nonumber\\
 && \quad\,　 +2
\Gamma_{sk}^j g^{sl}g^{km}  \frac{\partial^2 f}{\partial x_l\partial x_m}　 -2
\Gamma_{sk}^j g^{sl}g^{km}\Gamma_{lm}^r \, \frac{\partial f}{\partial x_r}
 - 2  g^{kl}\Gamma_{sk}^j \Gamma_{lr}^s  g^{rt} \frac{\partial f}{\partial x_t}. \nonumber\end{eqnarray}
 By virtue of $(\mu^{-1})^{;j}= -\mu^{-2} \mu^{;j}$ and $\frac{\partial((\mu+\rho)^{-\frac{1}{2}})}{\partial x_m} =-\frac{1}{2} (\mu+\rho)^{-\frac{3}{2}} \frac{\partial \mu}{\partial x_m}$, we finally obtain
     \begin{eqnarray}   &&    \big( \mbox{div}\,(\mu S)\big)^j  = \mu (\mu+\rho)^{-\frac{1}{2}} \Delta_g w^j +  \Big(\Delta_g f +\mu ( \Delta_g \mu^{-1}) f + (\mu+\rho)^{-\frac{1}{2}} \frac{\partial \mu}{\partial x_k} w^k \Big)^{;j}  \nonumber\\
       &&   \quad\, +\bigg[ \mu^2(\mu+\rho)^{-\frac{1}{2}}  (\mu^{-1})^{;j}  \frac{\partial w^k}{\partial x_k}   \bigg. +\Big((\mu+\rho)^{-\frac{1}{2}} -\mu(\mu+\rho)^{-\frac{3}{2}}\Big)  \frac{\partial \mu}{\partial x_m} g^{ml} \frac{\partial w^j}{\partial x_l} + 2\mu (\mu+\rho)^{-\frac{1}{2}}g^{ml} \Gamma_{km}^j   \frac{\partial w^k}{\partial x_l}\bigg] \nonumber \\
   &&   \quad\, + \bigg[ \!\Big( \mu \Delta_g( (\mu\!+\!\rho)^{-\frac{1}{2}} ) +\frac{\partial \mu}{\partial x_l} \big( (\mu\!+\!\rho)^{-\frac{1}{2}} \big)^{;l} \Big) w^j  + \mu^2(\mu\!+\!\rho)^{-\frac{1}{2}} (\mu^{\!-1})^{;j} \Gamma_{kl}^l w^k +
   \mu^2 (\mu^{\!-1})^{;j}  \frac{\partial ((\mu\!+\!\rho)^{-\frac{1}{2}})}{\partial x_k} w^k \nonumber\\
   &&   \quad\, +
   \frac{\rho}{(\mu+\rho)^{\frac{3}{2}}} g^{ml} \Gamma_{km}^j \frac{\partial \mu}{\partial x_l} w^k
    +
  \mu (\mu+\rho)^{-\frac{1}{2}}  \Big( g^{ml} \frac{\partial \Gamma_{kl}^j}{\partial x_m} + g^{ml} \Gamma_{hl}^j \Gamma_{km}^h -
g^{ml}\Gamma_{kh}^j \Gamma_{ml}^h \Big)  w^k  \nonumber \\
&&
  \quad\, +  \mu (\mu+\rho)^{-\frac{1}{2}}R^j_k \, w^k -   (\mu+\rho)^{-\frac{1}{2}} \big(\frac{\partial \mu}{\partial x_k}\big)^{;j} w^k \bigg]　 +2
\Gamma_{sr}^j g^{sl}g^{rm}  \frac{\partial^2 f}{\partial x_l\partial x_m}  \nonumber\\
    &&  \quad\, + \bigg[\! 2  R^j_mg^{lm}\frac{\partial f}{\partial x_l}  -2\mu g^{ml}  \frac{\partial ((\mu^{-1})^{;j}\!)}{\partial x_m} \, \frac{\partial f}{\partial x_l}  -2\mu g^{ml} \Gamma_{sm}^j (\mu^{-1})^{;s} \frac{\partial f}{\partial x_l}
 \nonumber\\
&&  \quad\, +    \Big( g^{mr} \frac{\partial \Gamma_{sr}^j}{\partial x_m}
-g^{mr}\Gamma_{hr}^j \Gamma_{sm}^h -g^{mr}\Gamma_{sh}^j \Gamma_{mr}^h \Big) g^{sl} \frac{\partial f}{\partial x_l} -2
\Gamma_{sh}^j g^{sr}g^{hm}\Gamma_{rm}^l \, \frac{\partial f}{\partial x_l}
 \bigg]  \nonumber
\\ &&
 \quad\, +\bigg[\!  - 2 \mu (\Delta_g \mu^{-1})^{;j} f　 - 2\mu R^j_l\,( \mu^{-1})^{;l}f
-   2\mu g^{ml}\Gamma_{sm}^j\, \frac{\partial ((\mu^{-1})^{;s})}{\partial x_l}\, f \nonumber\\
&&   \quad\, - \mu \Big( g^{ml} \frac{\partial \Gamma_{sl}^j}{\partial x_m} + g^{ml} \Gamma_{hl}^j \Gamma_{sm}^h -
g^{ml}\Gamma_{sh}^j \Gamma_{ml}^h \Big)(\mu^{-1})^{;s} f
  -   2
  \frac{\partial \mu}{\partial x_m} \, (\mu^{-1})^{;m;j} f \bigg],
   \nonumber
       \end{eqnarray}
      If we choose \begin{eqnarray} \label{200420-11} p= \Delta_g f +\mu ( \Delta_g \mu^{-1}) f + (\mu+\rho)^{-\frac{1}{2}} \frac{\partial \mu}{\partial x_k} w^k ,\end{eqnarray}
then, by the above calculation and (\ref{2020.4.18-6}), \begin{eqnarray} \label{200420-12}  \begin{pmatrix} u\\ p \end{pmatrix} =\begin{pmatrix} (\mu+\rho)^{-\frac{1}{2}} w +\mu^{-1} \nabla_g f -
 f\nabla_g \mu^{-1} \\ \Delta_g f +\mu ( \Delta_g \mu^{-1}) f + (\mu+\rho)^{-\frac{1}{2}} \frac{\partial \mu}{\partial x_k} w^k  \end{pmatrix}\end{eqnarray}
is a solution of the stationary Stokes equations (\ref{2020.4.18-1}) provided   $ \begin{pmatrix} w\\ f\end{pmatrix} $
satisfies
    \begin{eqnarray} \label{200421-13} \left\{ \begin{array}{ll} L_j(w, f)=0, \,\quad j=1,\cdots,n, \\
   \Delta_g f  -\mu (\Delta_g\mu^{-1}) f + \mu (\mu+\rho)^{-\frac{1}{2}} \frac{\partial w^k }{\partial x_k}  +
 \Big( \mu (\mu+\rho)^{-\frac{1}{2}} \Gamma_{kl}^l +\mu \frac{\partial ((\mu+\rho)^{-\frac{1}{2}})}{\partial x_k} \Big) w^k =0,\end{array} \right. \end{eqnarray}
where
   \begin{eqnarray}    &&   \label{200422-3}    L_j(w,f):    = \Delta_g w^j +\bigg[ \mu (\mu^{-1})^{;j} \, \frac{\partial w^k}{\partial x_k}  + \frac{\rho}{\mu(\mu+\rho)} \frac{\partial \mu}{\partial x_m} g^{ml} \frac{\partial w^j}{\partial x_l} +2 g^{ml} \Gamma_{km}^j \frac{\partial w^k}{\partial x_l}\bigg]  \\
   &&   + \bigg[ \Big(  (\mu+ \rho)^{\frac{1}{2}} \Delta_g ((\mu+\rho)^{-\frac{1}{2}} ) +\mu^{-1} (\mu+\rho)^{\frac{1}{2}}  \frac{\partial \mu}{\partial x_l} \big((\mu+\rho)^{-\frac{1}{2}}\big)^{;l} \Big) w^j +  \mu (\mu^{-1})^{;j}  \Gamma_{kl}^l w^k \nonumber\\
&& + \mu (\mu+\rho)^{\frac{1}{2}} \,  (\mu^{-1})^{;j} \,\frac{\partial ((\mu+\rho)^{-\frac{1}{2}})}{\partial x_k} w^k
+ \frac{\rho}{\mu(\mu+\rho)} g^{ml} \Gamma_{km}^j \frac{\partial \mu}{\partial x_l} w^k
 \nonumber
\\
&& +   \Big( g^{ml} \frac{\partial \Gamma_{kl}^j}{\partial x_m} + g^{ml} \Gamma_{hl}^j \Gamma_{km}^h -
g^{ml}\Gamma_{kh}^j \Gamma_{ml}^h \Big)  w^k  +R^j_k \, w^k  - \mu^{-1}  \big(\frac{\partial \mu}{\partial x_k}\big)^{;j} w^k \bigg]\nonumber \\
&&  　 +2\mu^{-1} (\mu+\rho)^{\frac{1}{2}} \Gamma_{sr}^j g^{sl}g^{rm}  \frac{\partial^2 f}{\partial x_l\partial x_m} + \bigg[\! 2\mu^{-1} (\mu+\rho)^{\frac{1}{2}}  R^j_m g^{lm}\frac{\partial f}{\partial x_l}   \nonumber\\
    &&  -2(\mu+\rho)^{\frac{1}{2}} g^{ml}  \frac{\partial ((\mu^{-1})^{;j}\!)}{\partial x_m} \, \frac{\partial f}{\partial x_l}  -2(\mu+\rho)^{\frac{1}{2}} g^{ml} \Gamma_{sm}^j (\mu^{-1})^{;s} \frac{\partial f}{\partial x_l}
 \nonumber\\
&& +   \mu^{-1} (\mu+\rho)^{\frac{1}{2} } \Big( g^{mr} \frac{\partial \Gamma_{sr}^j}{\partial x_m} - g^{mr} \Gamma_{hr}^j \Gamma_{sm}^h  -
g^{mr}\Gamma_{sh}^j \Gamma_{mr}^h \Big) g^{sl} \frac{\partial f}{\partial x_l}\nonumber
\\ &&
 -2\mu^{-1} (\mu+\rho)^{\frac{1}{2}}
\Gamma_{sh}^j g^{sr}g^{hm}\Gamma_{rm}^l \, \frac{\partial f}{\partial x_l}
 \bigg] +\bigg[\!  - 2  (\mu+\rho)^{\frac{1}{2}} ( \Delta_g \mu^{-1})^{;j} f　\nonumber\\
&&  - 2(\mu+\rho)^{\frac{1}{2}} R^j_l\, (\mu^{-1})^{;l} f
-   2(\mu+\rho)^{\frac{1}{2}} g^{ml}\Gamma_{sm}^j\, \frac{\partial ((\mu^{-1})^{;s}\!)}{\partial x_l}\, f  - (\mu+\rho)^{\frac{1}{2} } \Big( g^{ml} \frac{\partial \Gamma_{sl}^j}{\partial x_m} \nonumber\\
&& + g^{ml} \Gamma_{hl}^j \Gamma_{sm}^h -
g^{ml}\Gamma_{sh}^j \Gamma_{ml}^h \Big)(\mu^{-1})^{;s} f
  -   2 \mu^{-1} (\mu+\rho)^{\frac{1}{2} }
   \frac{\partial \mu}{\partial x_m} \, (\mu^{-1})^{;m;j} f \bigg].
   \nonumber
       \end{eqnarray}
Clearly, (\ref{200421-13}) is a system of second-order  linear elliptic equations in $\Omega$.
 We further consider the following two elliptic boundary value problems:
  \begin{eqnarray}\label{20200524-3}  \left\{\!\!\! \begin{array}{ll} L_j(w, f)=0, \quad j=1,\cdots,n \quad &\mbox{in}\;\; \Omega, \\
     \Delta_g f \! -\!\mu (\Delta_g\mu^{\!-1}) f \!+\! \mu (\mu+\rho)^{-\frac{1}{2}} \frac{\partial w^k }{\partial x_k} \! +\!
 \Big(\! \mu (\mu\!+\!\rho)^{\!-\frac{1}{2}} \Gamma_{kl}^l \!+\!\mu \frac{\partial ((\mu\!+\!\rho)^{\!-\frac{1}{2}})}{\partial x_k}\! \Big) w^k=0 \;\; &\mbox{in}\;\; \Omega,\\
 (\mu+\rho)^{-\frac{1}{2}} w +\mu^{-1} \nabla_g f -
 f\nabla_g \mu^{-1}=u_0 \;\; &\mbox{on}\;\, \partial \Omega \end{array} \right.
  \end{eqnarray}
 and
  \begin{eqnarray}\label{20200524-4}  \left\{\!\!\! \begin{array}{ll} L_j(w, f)=0, \quad j=1,\cdots,n \quad &\mbox{in}\;\; \Omega, \\
     \Delta_g f \! -\!\mu (\Delta_g\mu^{\!-1}) f \!+\! \mu (\mu+\rho)^{-\frac{1}{2}} \frac{\partial w^k }{\partial x_k} \! +\!
 \Big(\! \mu (\mu\!+\!\rho)^{\!-\frac{1}{2}} \Gamma_{kl}^l \!+\!\mu \frac{\partial ((\mu\!+\!\rho)^{\!-\frac{1}{2}})}{\partial x_k}\! \Big) w^k=0 \;\; &\mbox{in}\;\; \Omega,\\
 (w,f)= (w_0,f_0) \;\; &\mbox{on}\;\, \partial \Omega. \end{array} \right.
  \end{eqnarray}
  If we discuss the corresponding eigenvalue problems with vanishing boundary conditions for the above two systems, we see that all eigenvalues are discrete and any eigenvalue of each kind problem will continuously vary in $\rho$. Thus we can choose a suitable constant $\rho=\tilde{\rho}\ge 0$ such that $0$ is neither an eigenvalue of (\ref{20200524-3}) nor an eigenvalue of (\ref{20200524-4}) when $u_0$ and $(w_0,f_0)$  being replacing by vanishing boundary conditions  $u_0=0$ and $(w_0,f_0)=0$, respectively.
It follows that for any $(w_0,f_0)\in (H^{\frac{3}{2}} (\partial \Omega))^{n}\times H^{\frac{3}{2}}(\partial \Omega)$, there is a uniquely solution $(w,f)\in (H^{2}(\Omega))^n \times H^2(\Omega)$ of the system (\ref{20200524-4}) (when $\rho$ being replaced by $\tilde{\rho}$) satisfying $(w,f)\big|_{\partial \Omega}=(w_0,f_0)$. Thus, we can define the Dirichlet-to-Neumann map ${\tilde{\Lambda}}_{\tilde{\rho},\mu,g}: (H^{\frac{3}{2}} (\partial \Omega))^{n}\times H^{\frac{3}{2}}(\partial \Omega)\to (H^{\frac{1}{2}} (\partial \Omega))^{n}\times H^{\frac{1}{2}}(\partial \Omega)$  associated with
 new system (\ref{20200524-4} ) by
\begin{eqnarray} {\tilde{\Lambda}}_{\tilde{\rho}, \mu,g} (w_0,f_0) =\frac{\partial (w, f)}{\partial \nu}\big|_{\partial \Omega}  \quad \; \mbox{for any}\;\; (w_0,f_0)\in H^{\frac{3}{2}}(\partial \Omega)\times H^{\frac{3}{2}}(\partial \Omega),\end{eqnarray}
where $(w,f)$ satisfies (\ref{20200524-4}).
The corresponding Cauchy data is ${\tilde{\mathcal{C}}}_{\tilde{\rho},\mu}= \{ (w,f)\big|_{\partial \Omega},\frac{\partial (w, f)}{\partial \nu}\big|_{\partial \Omega}\}$.

 \vskip 0.32 true cm

\noindent{\bf Lemma 2.2.}  {\it  The Cauchy data corresponding to the $\Lambda_{\mu,g}$ is equivalent to the Cauchy data corresponding to ${\tilde{\Lambda}}_{\tilde{\rho},\mu, g}$}.

\vskip 0.25 true cm

 \noindent  {\it Proof.} \  Recall that for the chosen $\tilde{\rho}\ge 0$, the real number $0$ is neither an eigenvalue of (\ref{20200524-3}) nor an eigenvalue of (\ref{20200524-4}) with $u_0$ and $(w_0,f_0)$  being replacing by vanishing boundary conditions  $u_0=0$ and $(w_0,f_0)=0$, respectively. Suppose $(u,p)\in (H^{2}(\Omega))^n\times H^1 (\Omega)$ is a solution of the Stokes equations (\ref{2020.4.18-1}) satisfying $\int_\Omega p\,dV=0$ with the boundary condition
 $u\big|_{\partial \Omega}= u_0$. By the previous discussion, we see that $(w,f)$ must be a unique solution of (\ref{20200524-3}) with boundary condition
 $ \big((\mu+\tilde{\rho})^{-\frac{1}{2}} w +\mu^{-1} \nabla_g f -
 f\nabla_g \mu^{-1}\big)\big|_{\partial \Omega}=u_0$, when $\rho$ is replaced by $\tilde{\rho}$. For such a $(w,f)$, if we set $(w,f)\big|_{\partial \Omega} =(w_0,f_0)$, then $(w,f)$ is also a uniquely solution of (\ref{20200524-4}) for the same constant $\tilde{\rho}$. Therefore, by this way we get Cauchy datum $\big((w_0,f_0),  \frac{\partial (w,f)}{\partial \nu}\big|_{\partial \Omega}\big)$.

 Conversely, for any  $(w_0,f_0)\in (H^{\frac{3}{2}}(\partial \Omega))^n \times H^{\frac{3}{2}}(\partial \Omega)$, let $(w,f)\in (H^{2}(\Omega))^n\times H^2 (\Omega)$  be a unique solution  of (\ref{20200524-4}) (when $\rho$ is replaced by $\tilde{\rho}$). By (\ref{200420-12}) we immediately get $(u,p)$, which satisfies  the Stokes equations (\ref{2020.4.18-1}) with boundary value $u\big|_{\partial \Omega}$ because the last equation in (\ref{20200524-4}) is exactly $\mbox{div}\; u=0$ in $\Omega$
  and the first $n$ equations are just $\,\mbox{div}\,\sigma_\mu (u,p) =\mbox{div}\, (\mu S)- \nabla_g p-\mu(\mu+\tilde{\rho})^{\frac{1}{2}} L (w,f)=0$ in $\Omega$. We may add a suitable constant to the above $p$ such that $\int_\Omega p\, dV=0$. This gives a Cauchy datum $(u\big|_{\partial \Omega}, \sigma_\mu(u,p)\nu\big|_{\partial \Omega})$ associated with the Stokes equations. Hence, the Cauchy data ${\mathcal{C}}_\mu=\{(u\big|_{\partial \Omega},\sigma_\mu (u,p)\mu\big|_{\partial \Omega})\}$ associated with the Stokes equations is equivalent to the Cauchy data ${\tilde{\mathcal{C}}}_{\tilde{\rho}, \mu}=\{(w,f)\big|_{\partial \Omega}, \frac{\partial (w,f)}{\partial \nu}\big|_{\partial \Omega})\}$ associated with the new system, and the desired conclusion is proved.  \qed

\vskip 1.49 true cm

\section{Factorization of equivalent new system and new Dirichlet-to Neumann map}

\vskip 0.45 true cm

   From now on, we will denote by
\begin{equation*}
  \begin{bmatrix}
\begin{BMAT}(@, 15pt, 15pt){c.c}{c.c}
   [a_{jk}]_{n\times n}&[b_j]_{n\times 1}\\
   [c_k]_{1\times n}& d
  \end{BMAT}
\end{bmatrix}
\end{equation*} the block matrix
\[ \begin{bmatrix}
\begin{BMAT}(@, 15pt, 15pt){ccc.c}{ccc.c}
    a_{11} & \cdots & a_{1n} & b_{1}\\
 \cdots & \cdots &\cdots & \vdots \\
   a_{n1} & \cdots & a_{nn} & b_{n}\\
    c_{1} & \cdots & c_{n} & d
  \end{BMAT}
\end{bmatrix},  \]
 where $\big[a_{jk}\big]_{n\times n}$, $\big[ b_j\big]_{n\times 1}$ and $\big[c_k\big]_{1\times n}$ are the $n\times n$ matrix \begin{eqnarray*}\begin{bmatrix} a_{11}& a_{12} & \cdots &a_{1n}\\
a_{21}& a_{22} & \cdots &a_{2n}\\
\cdots & \cdots& \cdots & \cdots\\
a_{n1}& a_{n2} & \cdots &a_{nn}\end{bmatrix}, \end{eqnarray*}
the $n\times 1$ matrix
\begin{eqnarray*}\begin{bmatrix} b_1 \\
\vdots \\ b_n\end{bmatrix} \end{eqnarray*}
and the $1\times n$ matrix \begin{eqnarray*}\begin{bmatrix} c_1 & \cdots & c_{n} \end{bmatrix}, \end{eqnarray*} respectively.

In what follows, we will let Greek indices run from $1$ to $n-1$, Roman indices from $1$ to $n$.
Then, in the local coordinates,  we can rewrite (\ref{200421-13}) as
\begin{align*}
 & \left\{\begin{bmatrix}
\begin{BMAT}(@, 15pt, 15pt){c.c}{c.c}
   [\delta_{jk}\Delta_g ]_{n\times n}&\left[2\mu^{-1} (\mu+\tilde{\rho})^{\frac{1}{2}}\Gamma_{sr}^j g^{sl}g^{rm} \frac{\partial^2}{ \partial x_l \partial x_m}\right]_{n\times 1}\\
   [0]_{1\times n}& \Delta_g
  \end{BMAT}
\end{bmatrix}\right.
\\
&   \!\! +\!\begin{bmatrix}\!
\begin{BMAT}(@, 1pt, 5pt){c.c}{c.c}
   \!\left[\! \mu (\mu^{\!-\!1}\!)^{;j}\!\frac{\partial }{\partial x_k}\!+\! \frac{\tilde{\rho} \delta_{jk}}{\mu(\mu+\tilde{\rho})} \frac{\partial \mu}{\partial x_m} g^{ml} \frac{\partial}{\partial x_l}  \!+\!2 g^{ml}\Gamma_{\!km}^j\frac{\partial }{\partial x_l}\right]   &   \! \left[\!\frac{(\mu\!+\!\tilde{\rho})^{\!\frac{1}{2}}}{\mu}\!\! \left\{\!\!\! \begin{array} {ll} \! 2R^j_m g^{lm} \!-\!2\mu g^{ml}\frac{\partial ((\mu^{\!-\!1})^{;j}\!)}{\partial x_m} \!  -\!2 \mu g^{ml} \Gamma_{\!sm}^j (\mu^{\!-\!1})^{;s}\\
        \!+ ( g^{mr}\frac{\partial \Gamma_{sr}^j }{\partial x_m}  \!-\!  g^{mr} \Gamma_{\!hr}^j \Gamma_{\!sm}^h \!-\! g^{mr} \Gamma_{\!sh}^j \Gamma_{\!mr}^h) g^{sl}  \\
       \! -2\Gamma_{sh}^j g^{sr}g^{hm}\Gamma_{rm}^l\end{array}\!\!\!\! \!\right\}\!\!\frac{\partial}{\partial x_l}\!\right]_{\!n\!\times\! 1}\\
   \left[\mu(\mu+\tilde{\rho})^{\!-\frac{1}{2}} \frac{\partial}{\partial x_k} \right]_{1\times n}& 0
  \end{BMAT}\!
\end{bmatrix}
\\
&    +\begin{bmatrix}
\begin{BMAT}(@, 5pt, 5pt){c.c}{c.c}
   \!\left[\!\Big((\mu\!+\!\tilde{\rho})^{\!\frac{1}{2}} \Delta_g ((\mu+\tilde{\rho})^{\!-\frac{1}{2}}) +\mu^{\!-\!1} (\mu+\tilde{\rho})^{\frac{1}{2}}   \frac{\partial \mu}{\partial x_l}   ((\mu+\tilde{\rho})^{\!-\frac{1}{2}})^{;l}\Big) \delta_{jk}  \right]_{n\times n}   &   \! \left[\, 0\,\right]_{n\times 1}\\
   \left[0 \right]_{1\times n}&  -\mu (\Delta_g \mu^{-1})
  \end{BMAT}
\end{bmatrix}
\\
\!\!& \! \left. \! +\!\begin{bmatrix}
\begin{BMAT}(@, 1pt, 1pt){c.c}{c.c}
  \! \!\left[\!\!\!\left. \begin{array}{ll}\! \mu(\mu^{\!-\!1}\!)^{;j} \Gamma_{\!kl}^l \!+\!\mu(\mu\!+\!\tilde{\rho})^{\!\frac{1}{2}}(\mu^{\!-\!1})^{;j} \frac{\partial (\mu\!+\!\tilde{\rho})^{\!-\!\frac{1}{2}}}{\partial x_k}\\
 \! \!+\! \big( g^{ml} \frac{\partial \Gamma_{\!kl}^j}{\partial x_m}\! + \!g^{ml}\Gamma_{\!hl}^j \Gamma_{\!km}^h \!\!-\! g^{ml} \Gamma_{\!kh}^j \Gamma_{\!ml}^h \big)\\
   \frac{\tilde{\rho}}{\mu(\mu\!+\!\tilde{\rho})} g^{ml} \Gamma_{\!km}^j \frac{\partial \mu}{\partial x_l} \!+\!R^j_k \!-\! \mu^{\!
-\!1}\big( \!\frac{\partial \mu}{\partial x_k}\big)^{;j}\end{array} \!\!
\right.\!\! \! \right]_{\!n\!\times \!n}   &   \! \left[\!(\mu\!+\!\tilde{\rho})^{\!\frac{1}{2}}\!\!\left\{\!\!\! \begin{array}{ll} \!\!-2 ( \Delta_g \mu^{\!-\!1})^{;j} - 2  R_l^j \,(  \mu^{\!-1})^{;l}  \\ \!\!-2  g^{ml} \Gamma_{\!sm}^j \frac{\partial((\mu^{\!-\!1})^{;s})}{\partial x_l}
  -2 \mu^{\!-\!1}\frac{\partial \mu}{\partial x_m} (\mu^{\!-\!1}\!)^{;m;j}\\
    \!\!- \!\big(\! g^{ml} \frac{\partial \Gamma^j_{\!sl}}{\partial x_m}\! + \!g^{ml} \Gamma^j_{\!hl}\Gamma^h_{\!sm} \!\!-\! g^{ml} \Gamma_{\!sh}^j\Gamma^h_{\!ml}\! \big) \!(u^{\!-\!1}\!)^{;s}\end{array}\!\!\!\!\!\right\}\!\!\right]_{\!n\!\times \!1}\\
   \left[\mu(\mu+\tilde{\rho})^{-\frac{1}{2}} \Gamma_{kl}^l + \mu \frac{\partial ((\mu+\tilde{\rho})^{-\frac{1}{2}})}{\partial x_k} \! \right]_{\!1\times n}&  0
  \end{BMAT}
\!\end{bmatrix}\!\right\}\!\!\begin{bmatrix} w^1\\
\vdots\\
w^n\\
f\end{bmatrix}\!\!=\!0,
\end{align*}
 where  \[ \delta_{jk}= \left\{\begin{array}{ll} 1 \;\; &\mbox{for}\;\, j=k\\
 0 \;\; &\mbox{for} \,\; j\ne k \end{array} \right.\]
 is the standard Kronecker symbol.

In order to describe the Dirichlet-to-Neumann map associated with  the equivalent new system, we first recall the construction of usual geodesic coordinates with respect to the
boundary (see p.$\,$1101 of \cite{LU}). For each $x'\in \partial \Omega$, let $r_{x'}: [0, \tau)\to \bar \Omega$ denote the unit-speed geodesic
starting at $x'$ and normal to $\partial \Omega$. If $x':=\{x_1, \cdots, x_{n-1}\}$  are any local coordinates for
$\partial \Omega$ near $x_0\in \partial \Omega$, we can extend them smoothly to functions on a neighborhood
of $x_0$ in $\Omega$ by letting them be constant along each normal geodesic $r_{x'}$. If we then
define $x_n$ to be the parameter along each $r_{x'}$, it follows easily that $\{x_1, \cdots, x_{n}\}$
form coordinates for $\Omega$ in some neighborhood of $x_0$, which we call the boundary
normal coordinates determined by $\{x_1, \cdots, x_{n-1}\}$. In these coordinates $x_n>0$ in
 $\Omega$, and $\partial \Omega$ is locally characterized by $x_n= 0$. A standard computation shows
that the metric $g$ on $\bar \Omega$ then has the form
 (see p.$\,$1101 of \cite{LU} or p.$\,$532 of \cite{Ta2})
\begin{eqnarray} \label{18/a-1} \; \quad\;\big[g_{jk} (x',x_{n}) \big]_{n\times n} = \begin{bmatrix}
 g_{11} (x',x_{n}) & g_{12} (x',x_{n}) & \cdots & g_{1,n-1} (x',x_{n}) & 0\\
 \cdots\cdots& \cdots\cdots & \cdots &\cdots\cdots  & \cdots\\
 g_{n-1,1} (x',x_{n})  & g_{n-1,2} (x',x_{n}) & \cdots & g_{n-1,n-1} (x',x_{n}) & 0\\
 0& 0& 0& 0&1 \end{bmatrix}.  \end{eqnarray}
       Furthermore, we can take a geodesic normal coordinate system for $(\partial \Omega, g|_{\partial \Omega})$ centered at $x_0=0$, with respect to $e_1, \cdots, e_{n-1}$, where  $e_1, \cdots, e_{n-1}$ are the principal curvature vectors. As Riemann showed, one has (see p.$\,$555 of \cite{Ta2}, or \cite{Spi2})
            \begin{eqnarray} \label{18/7/14/1} \begin{split}& g_{jk}(x_0)= \delta_{jk}, \; \; \frac{\partial g_{jk}}{\partial x_l}(x_0)
 =0  \;\;  \mbox{for all} \;\; 1\le j,k,l \le n-1,\\
 & \frac{1}{2}\,\frac{\partial g_{jk}}{\partial x_n} (x_0) =\kappa_k\delta_{jk}  \;\;  \mbox{for all} \;\; 1\le j,k \le n-1,\end{split}\quad \qquad \qquad \quad
 \end{eqnarray}
where  $\kappa_1,\cdots, \kappa_{n-1}$ are the principal curvatures of $\partial \Omega$ at point $x_0=0$.
Under this normal coordinates, we take $-{\nu}(x)=[0,\cdots, 0,1]^t$.
By (\ref{18/a-1}) we immediately see that the inverse of metric tensor $g$ in the boundary normal coordinates has form:
    \begin{eqnarray*} g^{-1}(x',x_n) =\begin{bmatrix} g^{11}(x', x_n) & \cdots & g^{1,n-1} (x', x_n)& 0 \\
    \cdots & \cdots  & \cdots  & \cdots\\
    g^{n-1, 1}(x',x_n) & \cdots & g^{n-1,n-1}(x',x_n)& 0\\
    0&\cdots \cdots &0 &1\end{bmatrix}. \end{eqnarray*}
 Note that under  the boundary normal coordinates,  we have \begin{eqnarray} \label{19.10.5-1}\Gamma_{nk}^n &=&\frac{1}{2} \sum_{m=1}^n g^{nm}\bigg(\frac{\partial g_{nm}}{\partial x_k} +\frac{\partial g_{km}}{\partial x_n} -\frac{\partial g_{nk}}{\partial x_m}\bigg) \\ &=&\frac{1}{2} \bigg(\frac{\partial g_{nn}}{\partial x_k}+\frac{\partial g_{kn}}{\partial x_n} -\frac{\partial g_{nk}}{\partial x_n}\bigg)=0,\nonumber\end{eqnarray}
   \begin{eqnarray} \Gamma_{nn}^l =\frac{1}{2}\sum\limits_{m=1}^n g^{lm}\bigg(\frac{\partial g_{nm}}{\partial x_n} +\frac{\partial g_{nm}}{\partial x_n} -\frac{\partial g_{nn}}{\partial x_m}\bigg) =0.\end{eqnarray}

 Thus, in the boundary normal coordinates, the above system of equations can be written as
 \begin{align*} & \left\{\Big( \frac{\partial^2}{ \partial x_n^2}  + \Gamma_{n\beta}^\beta \, \frac{\partial }{\partial x_n} + g^{\alpha\beta} \frac{\partial^2}{\partial x_\alpha \partial x_\beta} + \big( g^{\alpha\beta} \Gamma_{\alpha\gamma}^\gamma +\frac{\partial g^{\alpha\beta}}{\partial x_\alpha} \big)\frac{\partial }{\partial x_\beta} \Big) I_{n+1} \right.\\
 &     + \begin{bmatrix}
\begin{BMAT}(@, 25pt, 25pt){c.c}{c.c}
   \!\left[0\right]_{n\times n}   &   \! \left[4 \mu^{-1} (\mu+\tilde{\rho})^{\frac{1}{2}} \Gamma_{\beta n}^j g^{\beta \alpha} \frac{\partial^2}{\partial x_\alpha\partial x_n}\right]_{n\times 1}\\
   \left[0 \right]_{1\times n}&  0  \end{BMAT}
\end{bmatrix}
+\begin{bmatrix}
\begin{BMAT}(@, 15pt, 15pt){c.c}{c.c}
   \!\left[0\right]_{n\times n}   &   \! \left[2 \mu^{-1} (\mu+\tilde{\rho})^{\frac{1}{2}} \Gamma_{\gamma \sigma}^j g^{\alpha \gamma} g^{\beta \sigma} \frac{\partial^2}{\partial x_\alpha\partial x_\beta} \right]_{n\times 1}\\
   \left[0 \right]_{1\times n}&  0  \end{BMAT}
\end{bmatrix}\\
 &      \!+\!\begin{bmatrix}
\begin{BMAT}(@, 15pt, 15pt){c.c}{c.c}
   \!\left[\Big(\!\delta_{kn\,} \mu(\!\,\mu^{\!-1}\!)^{;j}\!+\!\frac{\delta_{jk} \tilde{\rho}}{\mu(\mu\!+\!\tilde{\rho})} \frac{\partial \mu}{\partial x_n}
    \!+\! 2\Gamma_{\!kn}^j
   \!\Big)  \frac{\partial }{\partial x_n}\right]_{\!n\!\times \!n} &  \! \left[\!\frac{(\mu\!+\!\tilde{\rho})^{\frac{1}{2}}}{\mu}\!\!\left\{\!\!\!\begin{array} {ll} 2 R^j_n  - 2\mu  \frac{\partial((\mu^{-1})^{;j})}{\partial x_n} - 2\mu \Gamma^j_{\alpha n} (\mu^{-1})^{;\alpha} \\
\!+  ( g^{\alpha\beta} \frac{\partial \Gamma_{\!n\alpha}^j}{\partial x_\beta} \!-\! g^{\alpha\gamma} \Gamma_{\!\beta \gamma}^j \Gamma_{\!n\alpha}^\beta \! \!-\! g^{\alpha\beta} \Gamma_{\!n\gamma}^j \Gamma_{\!\alpha\beta}^\gamma ) \\
 - 2 g^{\alpha \gamma} g^{\beta \sigma} \Gamma_{\gamma \sigma}^j \Gamma_{\alpha\beta}^n
  \end{array}\!\!\!\! \right\}\!\!\frac{\partial }{\partial x_n} \! \right]_{\!n\!\times\! 1}\\
     \left[\delta_{nk}\, \mu(\mu+\tilde{\rho})^{-\frac{1}{2}} \frac{\partial }{\partial x_n}\right]_{1\times n} & 0 \end{BMAT}\!
\end{bmatrix}\\
&      \!+\!\begin{bmatrix}
\begin{BMAT}(@, 5pt, 5pt){c.c}{c.c}
  \! \!\left[\!(\!1\!-\!\delta_{nk}\!)\mu(\mu^{\!-\!1}\!)^{;j} \!\frac{\partial }{\partial x_k} \!+\!
  \big(\!\frac{\delta_{jk} \tilde{\rho}}{\mu(\mu\!+\!\tilde{\rho})}\frac{\partial \mu}{\partial x_\alpha}
  \!+\! 2  \Gamma^j_{\!k\alpha}\!\big)g^{\alpha\beta} \frac{\partial }{\partial x_\beta}  \! \right]_{\!n\!\times \!n}   \! &  \! \left[\!\frac{(\mu\!+\!\tilde{\rho})^{\frac{1}{2}}}{\mu}\!\!\left\{\!\!\!\begin{array} {ll} \!2 R^j_\alpha g^{\alpha\beta}  \!\!-\! 2\mu g^{\alpha\beta}  \frac{\partial((\mu^{\!-1}\!)^{;j}\!)}{\partial x_\alpha}\\
  \! -\! 2\mu g^{\alpha\beta} \Gamma^j_{s\alpha } \!(\mu^{\!-1})^{;s}  \!
 -\! 2\Gamma_{sh}^j   g^{sr}\! g^{hm} \Gamma_{rm}^\beta\\
  \!-\!( g^{mr}\! \frac{\partial \Gamma_{\!\alpha r}^j}{\partial x_m}
  \!-\!\! g^{mr} \Gamma_{\!hr}^j \Gamma_{\!\alpha m}^h \!
 \!-\!\! g^{mr} \!\Gamma_{\!\alpha h}^j\! \Gamma_{\!mr}^h )g^{\alpha\!\beta} \\
 \end{array}\!\!\!\!\! \right\}\!\!\frac{\partial }{\partial x_\beta} \!\! \right]_{\!n\!\times \!1}\!\\
     \left[(1-\delta_{nk})\mu(\mu+\tilde{\rho})^{-\frac{1}{2}}  \frac{\partial }{\partial x_k}\right]_{1\times n} & 0\!\end{BMAT}\!
\end{bmatrix}\\
&       +\begin{bmatrix}
\begin{BMAT}(@, 15pt, 15pt){c.c}{c.c}
   \!\left[\!\Big((\mu+\tilde{\rho})^{\!\frac{1}{2}}\Delta_g ((\mu+\tilde{\rho})^{\!-\frac{1}{2}}) +\mu^{\!-\!1} (\mu\!+\!\tilde{\rho})^{\frac{1}{2}}   \frac{\partial \mu}{\partial x_l}   ((\mu+\tilde{\rho})^{\!-\frac{1}{2}})^{;l}\Big) \delta_{jk}  \right]_{n\times n}   &   \! \left[\, 0\,\right]_{n\times 1}\\
   \left[0 \right]_{1\times n}&  -\mu (\Delta_g \mu^{-1})
  \end{BMAT}
\end{bmatrix}
\\
\!\!&          \! \left. \! +\!\begin{bmatrix}
\begin{BMAT}(@, 1pt, 1pt){c.c}{c.c}
  \! \!\left[\!\!\!\left. \begin{array}{ll}\! \mu(\mu^{\!-\!1}\!)^{;j} \Gamma_{\!kl}^l \!+\!\mu(\mu\!+\!\tilde{\rho})^{\!\frac{1}{2}}(\mu^{\!-\!1})^{;j} \frac{\partial ( (\mu\!+\!\tilde{\rho})^{\!-\!\frac{1}{2}})}{\partial x_k}\\
 \! \!+\! \big( g^{ml} \frac{\partial \Gamma_{\!kl}^j}{\partial x_m}\! + \!g^{ml}\Gamma_{\!hl}^j \Gamma_{\!km}^h \!\!-\! g^{ml} \Gamma_{\!kh}^j \Gamma_{\!ml}^h \big)\\
   \frac{\tilde{\rho}}{\mu(\mu\!+\!\tilde{\rho})} g^{ml} \Gamma_{\!km}^j \frac{\partial \mu}{\partial x_l} \!+\!R^j_k \!-\! \mu^{\!
-\!1}\big( \!\frac{\partial \mu}{\partial x_k}\big)^{;j}\end{array} \!\!
\right.\!\! \! \right]_{\!n\!\times \!n}   &   \! \left[\!(\mu\!+\!\tilde{\rho})^{\!\frac{1}{2}}\!\!\left\{\!\!\! \begin{array}{ll} \!\!-2 ( \Delta_g \mu^{\!-\!1})^{;j} - 2  R_l^j \,(  \mu^{\!-1})^{;l}  \\ \!\!-2  g^{ml} \Gamma_{\!sm}^j \frac{\partial((\mu^{\!-\!1})^{;s})}{\partial x_l}
  -2 \mu^{\!-\!1}\frac{\partial \mu}{\partial x_m} (\mu^{\!-\!1}\!)^{;m;j}\\
    \!\!- \!\big(\! g^{ml} \frac{\partial \Gamma^j_{\!sl}}{\partial x_m}\! + \!g^{ml} \Gamma^j_{\!hl}\Gamma^h_{\!sm} \!\!-\! g^{ml} \Gamma_{\!sh}^j\Gamma^h_{\!ml}\! \big) \!(u^{\!-\!1}\!)^{;s}\end{array}\!\!\!\!\!\right\}\!\!\right]_{\!n\!\times \!1}\\
   \left[\mu(\mu+\tilde{\rho})^{-\frac{1}{2}} \Gamma_{kl}^l + \mu \frac{\partial ((\mu+\tilde{\rho})^{-\frac{1}{2}})}{\partial x_k} \! \right]_{\!1\times n}&  0
  \end{BMAT}
\!\end{bmatrix}\!\right\}\!\!\begin{bmatrix} w^1\\
\vdots\\
w^n\\
f\end{bmatrix}\!\!=\!0.
\end{align*}
That is,
\begin{eqnarray} \label{200425-1}
 \frac{\partial^2}{\partial x_n^2} I_{n+1} +B \frac{\partial }{\partial x_n} +C=0, \end{eqnarray}
 where
  \begin{align} & B:=\begin{bmatrix}
\begin{BMAT}(@, 15pt, 15pt){c.c}{c.c}
  \! \left[\Gamma_{n\beta}^\beta \delta_{jk}    \right]_{n\times n}   &  \left[\!4 \mu^{-1} (\mu+\tilde{\rho})^{\frac{1}{2}} \Gamma_{\beta n}^j  g^{\alpha\beta}\frac{\partial }{\partial x_\alpha} \!\right]_{n\!\times \!1} \\
          \left[0 \right]_{1\!\times n}  & \Gamma_{n\beta}^\beta  \end{BMAT}
\end{bmatrix}\\
&       +\begin{bmatrix}
\begin{BMAT}(@, 15pt, 15pt){c.c}{c.c}
  \! \left[\delta_{kn} \mu (\mu^{-1})^{;j}\!+\! \frac{\delta_{jk} \tilde{\rho}}{\mu(\mu\!+\!\tilde{\rho})}\frac{\partial \mu}{\partial x_n}
  \!+\! 2\Gamma_{kn}^j
   \right]_{\!n\times \!n}   &\left[\!\frac{(\mu\!+\!\tilde{\rho})^{\frac{1}{2}}}{\mu}\!\!\left\{\!\!\!\begin{array} {ll} 2 R^j_n  - 2\mu  \frac{\partial((\mu^{-1})^{;j}\!)}{\partial x_n} - 2\mu \Gamma^j_{\alpha n} (\mu^{-1})^{;\alpha} \\
\!+  ( g^{\alpha\beta} \frac{\partial \Gamma_{\!n\alpha}^j}{\partial x_\beta} \!-\! g^{\alpha\gamma} \Gamma_{\!\beta\gamma}^j \Gamma_{\!n\alpha}^\beta \! \!-\! g^{\alpha\beta} \Gamma_{\!n\gamma}^j \Gamma_{\!\alpha\beta}^\gamma ) \\
 - 2 g^{\alpha \gamma} g^{\beta \sigma} \Gamma_{\gamma \sigma}^j \Gamma_{\alpha\beta}^n
  \end{array}\!\!\!\! \right\}\!\right]_{\!n\!\times \!1}
        \\
     \left[\delta_{nk} \,\mu(\mu+\tilde{\rho})^{-\frac{1}{2}}\right]_{1\!\times \!n}  &    0  \end{BMAT}
\end{bmatrix}\nonumber \end{align}
 and \begin{eqnarray}C=C_2+C_1+C_0,\end{eqnarray}
 where
 \begin{align} \label{ 200425-3}& C_2=g^{\alpha\beta} \frac{\partial^2}{\partial x_\alpha \partial x_\beta} I_{n+1} +\begin{bmatrix}
\begin{BMAT}(@, 15pt, 15pt){c.c}{c.c}
   \!\left[0\right]_{n\times n}   &   \! \left[2 \mu^{-1} (\mu+\tilde{\rho})^{\frac{1}{2}} \Gamma_{\gamma \sigma}^j g^{\alpha \gamma} g^{\beta \sigma} \frac{\partial^2}{\partial x_\alpha\partial x_\beta} \right]_{n\times 1}\\
   \left[0 \right]_{1\times n}&  0  \end{BMAT}
\end{bmatrix},  \end{align}
\begin{align}\label{200425-4}  & \,C_1:=
  \Big( \big( g^{\alpha\beta} \Gamma_{\alpha\gamma}^\gamma +\frac{\partial g^{\alpha\beta}}{\partial x_\alpha} \big)\frac{\partial }{\partial x_\beta} \Big) I_{n+1}\qquad \qquad \qquad \qquad\qquad \qquad \qquad \qquad\;\;\;\; \end{align}
\begin{align*} &       \!+\!\begin{bmatrix}
\begin{BMAT}(@, 5pt, 5pt){c.c}{c.c}
  \! \!\left[\!(\!1\!-\!\delta_{nk}\!)\mu(\mu^{\!-\!1}\!)^{;j} \!\frac{\partial }{\partial x_k} \!+\!
  \big(\!\frac{\delta_{jk} \tilde{\rho}}{\mu(\mu\!+\!\tilde{\rho})}\frac{\partial \mu}{\partial x_\alpha}
  \!+\! 2  \Gamma^j_{\!k\alpha}\!\big)g^{\alpha\beta} \frac{\partial }{\partial x_\beta}  \! \right]_{\!n\!\times \!n}   \! &  \! \left[\!\frac{(\mu\!+\!\tilde{\rho})^{\frac{1}{2}}}{\mu}\!\!\left\{\!\!\!\begin{array} {ll} \!2 R^j_\alpha g^{\alpha\beta}  \!\!-\! 2\mu g^{\alpha\beta}  \frac{\partial((\mu^{\!-1}\!)^{;j}\!)}{\partial x_\alpha}\\
  \! -\! 2\mu g^{\alpha\beta} \Gamma^j_{s\alpha } \!(\mu^{\!-1})^{;s}  \!
 -\! 2\Gamma_{sh}^j   g^{sr}\! g^{hm} \Gamma_{rm}^\beta\\
  \!-\!( g^{mr}\! \frac{\partial \Gamma_{\!\alpha r}^j}{\partial x_m}
  \!-\!\! g^{mr} \Gamma_{\!hr}^j \Gamma_{\!\alpha m}^h \!
 \!-\!\! g^{mr} \!\Gamma_{\!\alpha h}^j\! \Gamma_{\!mr}^h )g^{\alpha\!\beta} \\
 \end{array}\!\!\!\!\! \right\}\!\!\frac{\partial }{\partial x_\beta} \!\! \right]_{\!n\!\times \!1}\!\\
     \left[(1-\delta_{nk})\mu(\mu+\tilde{\rho})^{-\frac{1}{2}}  \frac{\partial }{\partial x_k}\right]_{1\times n} & 0\!\end{BMAT}\!
\end{bmatrix},\nonumber
  \end{align*}
 \begin{align} \label{200425-6} & C_0:=
 \begin{bmatrix}
\begin{BMAT}(@, 15pt, 15pt){c.c}{c.c}
   \!\left[\!\Big((\mu+\tilde{\rho})^{\!\frac{1}{2}}\Delta_g ((\mu+\tilde{\rho})^{\!-\frac{1}{2}}) +\mu^{\!-\!1} (\mu\!+\!\tilde{\rho})^{\frac{1}{2}}   \frac{\partial \mu}{\partial x_l}   ((\mu+\tilde{\rho})^{\!-\frac{1}{2}})^{;l}\Big) \delta_{jk}  \right]_{n\times n}   &   \! \left[\, 0\,\right]_{n\times 1}\\
   \left[0 \right]_{1\times n}&  -\mu (\Delta_g \mu^{-1})
  \end{BMAT}
\end{bmatrix} \qquad \end{align}
\begin{align*}
\!\!&          \!  +\!\begin{bmatrix}
\begin{BMAT}(@, 1pt, 1pt){c.c}{c.c}
  \! \!\left[\!\!\!\left. \begin{array}{ll}\! \mu(\mu^{\!-\!1}\!)^{;j} \Gamma_{\!kl}^l \!+\!\mu(\mu\!+\!\tilde{\rho})^{\!\frac{1}{2}}(\mu^{\!-\!1})^{;j} \frac{\partial ((\mu\!+\!\tilde{\rho})^{\!-\!\frac{1}{2}}\!)}{\partial x_k}\\
 \! \!+\! \big( g^{ml} \frac{\partial \Gamma_{\!kl}^j}{\partial x_m}\! + \!g^{ml}\Gamma_{\!hl}^j \Gamma_{\!km}^h \!\!-\! g^{ml} \Gamma_{\!kh}^j \Gamma_{\!ml}^h \big)\\
   \frac{\tilde{\rho}}{\mu(\mu\!+\!\tilde{\rho})} g^{ml} \Gamma_{\!km}^j \frac{\partial \mu}{\partial x_l} \!+\!R^j_k \!-\! \mu^{\!
-\!1}\big( \!\frac{\partial \mu}{\partial x_k}\big)^{;j}\end{array} \!\!
\right.\!\! \! \right]_{\!n\!\times \!n}   &   \! \left[\!(\mu\!+\!\tilde{\rho})^{\!\frac{1}{2}}\!\!\left\{\!\!\! \begin{array}{ll} \!\!-2 ( \Delta_g \mu^{\!-\!1})^{;j} - 2  R_l^j \,(  \mu^{\!-1})^{;l}  \\ \!\!-2  g^{ml} \Gamma_{\!sm}^j \frac{\partial((\mu^{\!-\!1})^{;s}\!)}{\partial x_l}
  -2 \mu^{\!-\!1}\frac{\partial \mu}{\partial x_m} (\mu^{\!-\!1}\!)^{;m;j}\\
    \!\!- \!\big(\! g^{ml} \frac{\partial \Gamma^j_{\!sl}}{\partial x_m}\! + \!g^{ml} \Gamma^j_{\!hl}\Gamma^h_{\!sm} \!\!-\! g^{ml} \Gamma_{\!sh}^j\Gamma^h_{\!ml}\! \big) \!(u^{\!-\!1}\!)^{;s}\end{array}\!\!\!\!\!\right\}\!\!\right]_{\!n\!\times \!1}\\
   \left[\mu(\mu+\tilde{\rho})^{-\frac{1}{2}} \Gamma_{kl}^l + \mu \frac{\partial ((\mu+\tilde{\rho})^{-\frac{1}{2}}\!)}{\partial x_k} \! \right]_{\!1\times n}&  0
  \end{BMAT}
\!\end{bmatrix}.\nonumber
 \end{align*}

\noindent Throughout this paper, we denote $\sqrt{-1}=\text{i}$.
\vskip 0.28 true cm

 \noindent{\bf Proposition 3.1.} \ {\it There exists a pseudodifferential operator $Q(x, D_{x'})$ of order one in $x'$ depending smoothly on $x_n$ such that \begin{eqnarray} \label{19.3.19-1:}  \frac{\partial^2}{\partial x_n^2} I_{n+1} +B\frac{\partial }{\partial x_n} +C = \left( \frac{\partial}{\partial x_n} \,I_n + B -Q \right) \left(\frac{\partial }{\partial x_n}\,I_n  +Q\right), \end{eqnarray}
 modulo a smoothing operator, where $D_{x'}=(D_{x_1}, \cdots, D_{x_{n-1}})$}, $\,D_{x_j}=\frac{1}{\text{i}}\,\frac{\partial }{\partial x_j}$.

 \vskip 0.25 true cm

 \noindent  {\it Proof.} \    We will divide this proof into several steps.

 Step 1.  \   Let us assume that we have a factorization \begin{eqnarray*}   \frac{\partial^2}{\partial x_n^2} I_{n+1} +B \frac{\partial }{\partial x_n} +C  =\left(\frac{\partial}{\partial x_n} \,I_{n+1}  +B - Q\right)\left( \frac{\partial}{\partial x_n}\,I_{n+1} +Q\right),\end{eqnarray*}
 i.e.,  \begin{eqnarray*} \frac{\partial^2}{\partial x_n^2} I_{n+1}  + B\frac{\partial}{\partial x_n} + C = \frac{\partial^2}{\partial x_n^2}
I_{n+1} +B\frac{\partial }{\partial x_n}  -Q\Big(\frac{\partial }{\partial x_n} \, I_{n+1}\Big) +\Big(\frac{\partial }{\partial x_{n} }I_{n+1}\Big) Q +B Q -Q^2.\end{eqnarray*}
This implies
\begin{eqnarray} \label{200425-8}  C-\left( \big( \frac{\partial}{\partial x_n} I_{n+1} \big) Q- Q \big(\frac{\partial }{\partial x_n}I_{n+1}\big) \right)
- BQ + Q^2 =0.\end{eqnarray}
    Let $q(x, \xi')$, $b(x,\xi')$ and $c(x,\xi')$ be the full symbols of $Q$ and $B$ and $C$, respectively. Clearly,  $q(x, \xi') \sim \sum_{j\ge 0} q_{1-j} (x, \xi')$,   $\;b(x,\xi')=b_1(x,\xi')+b_0(x,\xi')$ and $c(x, \xi')= c_2(x,\xi') +c_1(x, \xi') +c_0(x, \xi')$,
    where \begin{eqnarray}
  \label{19.3.19-4'} & b_1(x, \xi') = \begin{bmatrix}
\begin{BMAT}(@, 15pt, 15pt){c.c}{c.c}
  \! \left[0   \right]_{n\times n}   &  \left[4\,\text{i}\, \mu^{-1} (\mu+\tilde{\rho})^{\frac{1}{2}}\Gamma_{\beta n}^j g^{\alpha\beta}  \xi_\alpha\!\right]_{n\!\times \!1} \\
          \left[0 \right]_{1\!\times n}  & 0  \end{BMAT}
\end{bmatrix},\qquad \qquad \qquad \end{eqnarray}
 \begin{eqnarray}\label{19.3.19-5} \!\!\!\!\!\!&& \;\,\;\, b_0(x, \xi') \!=\! \Gamma_{n\beta}^\beta I_{n+1} \\
 &&\; +\!\begin{bmatrix}\!
\begin{BMAT}(@, 5pt, 5pt){c.c}{c.c}
  \! \left[\delta_{kn} \mu (\mu^{-1})^{;j}\!+\! \frac{\delta_{jk} \tilde{\rho}}{\mu(\mu\!+\!\tilde{\rho})} \frac{\partial \mu}{\partial x_n}\!
   +\! 2\Gamma_{\!kn}^j \!
  \right]_{\!n\!\times\! n}   &\left[\!\frac{(\mu\!+\!\tilde{\rho})^{\frac{1}{2}}}{\mu}\!\!\left\{\!\!\!\begin{array} {ll} 2 R^j_n  \!-\! 2\mu  \frac{\partial(\mu^{\!-\!1})^{\!;j}}{\partial x_n}\! -\! 2\mu \Gamma^j_{\alpha n} (\mu^{\!-\!1})^{\!;\alpha} \\
\!+ \! ( g^{\alpha\beta} \frac{\partial \Gamma_{\!n\alpha}^j}{\partial x_\beta} \!-\! g^{\alpha\gamma} \Gamma_{\!\beta\gamma}^j \Gamma_{\!n\alpha}^\beta \! \!-\! g^{\alpha\beta} \Gamma_{\!n\gamma}^j \Gamma_{\!\alpha\beta}^\gamma ) \\
 - 2 g^{\alpha \gamma} g^{\beta \sigma} \Gamma_{\gamma \sigma}^j \Gamma_{\!\alpha\beta}^n
  \end{array}\!\!\!\! \right\}\!\right]_{\!n\!\times \!1}
        \\
     \left[\delta_{nk} \,\mu(\mu+\tilde{\rho})^{-\frac{1}{2}}\right]_{1\!\times \!n}  &    0  \!\end{BMAT}\!
\end{bmatrix};\nonumber \end{eqnarray}
  \begin{eqnarray} \label{20200515-5} c_2(x, \!\xi')= -g^{\alpha\beta} \xi_\alpha \xi_\beta I_{n+1} +\begin{bmatrix}
\begin{BMAT}(@, 15pt, 15pt){c.c}{c.c}
   \left[0\right]_{n\times n}   &   \! \left[-2 \mu^{-1} (\mu+\tilde{\rho})^{\frac{1}{2}} \Gamma_{\gamma \sigma}^j g^{\alpha \gamma} g^{\beta \sigma} \xi_\alpha\xi_\beta \right]_{n\times 1}\\
   \left[0 \right]_{1\times n}&  0  \end{BMAT}
\end{bmatrix},\end{eqnarray}
\begin{eqnarray} & \label{20200515-6}  \!\!\!\! c_1(x, \!\xi')\! =
  \text{i} \, \big( g^{\alpha\beta} \Gamma_{\alpha\gamma}^\gamma +\frac{\partial g^{\alpha\beta}}{\partial x_\alpha} \big) \xi_\beta \,I_{n+1}\qquad \qquad \qquad \qquad \qquad \qquad \qquad \;\;\quad\quad \quad \;\, \;\; \end{eqnarray}
  \begin{align*}&  \!+\! \begin{bmatrix}
\begin{BMAT}(@, 5pt, 5pt){c.c}{c.c}
  \! \!\left[\!\text{i}(\!1\!-\!\delta_{nk}\!)\mu(\mu^{\!-\!1}\!)^{;j} \!\xi_k \!+\!
  \text{i}\big(\!\frac{\delta_{jk} \tilde{\rho}}{\mu(\mu\!+\!\tilde{\rho})}\frac{\partial \mu}{\partial x_\alpha}
  \!+\! 2  \Gamma^j_{\!k\alpha}\!\big)g^{\alpha\beta} \xi_\beta  \! \right]_{\!n\!\times \!n}   \! &  \! \left[\!\frac{(\mu\!+\!\tilde{\rho})^{\frac{1}{2}} \text{i}}{\mu}\!\left\{\!\!\!\begin{array} {ll} \!2 R^j_\alpha g^{\alpha\beta}  \!\!-\! 2\mu g^{\alpha\beta}  \frac{\partial(\mu^{\!-1}\!)^{;j}}{\partial x_\alpha}\\
  \! -\! 2\mu g^{\alpha\beta} \Gamma^j_{s\alpha } \!(\mu^{\!-1})^{;s}  \!
 -\! 2\Gamma_{sh}^j   g^{sr}\! g^{hm} \Gamma_{rm}^\beta\\
  \!-\!( g^{mr}\! \frac{\partial \Gamma_{\!\alpha r}^j}{\partial x_m}
  \!-\!\! g^{mr} \Gamma_{\!hr}^j \Gamma_{\!\alpha m}^h \!
 \!-\!\! g^{mr} \!\Gamma_{\!\alpha h}^j\! \Gamma_{\!mr}^h )g^{\alpha\!\beta} \\
 \end{array}\!\!\!\!\! \right\}\!\xi_\beta \!\! \right]_{\!n\!\times \!1}\!\\
     \left[\text{i} (1-\delta_{nk})\mu(\mu+\tilde{\rho})^{-\frac{1}{2}}  \xi_k\right]_{1\times n} & 0\!\end{BMAT}\!
\end{bmatrix},\nonumber\end{align*}
 \begin{eqnarray} && \quad\,c_0(x, \xi') \!=
 \begin{bmatrix}
\begin{BMAT}(@, 15pt, 15pt){c.c}{c.c}
   \!\left[\!\Big((\mu+\tilde{\rho})^{\!\frac{1}{2}}\Delta_g ((\mu+\tilde{\rho})^{\!-\frac{1}{2}}) +\mu^{\!-\!1} (\mu\!+\!\tilde{\rho})^{\frac{1}{2}}   \frac{\partial \mu}{\partial x_l}   ((\mu+\tilde{\rho})^{\!-\frac{1}{2}})^{;l}\Big) \delta_{jk}  \right]_{n\times n}   &   \! \left[\, 0\,\right]_{n\times 1}\\
   \left[0 \right]_{1\times n}&  -\mu (\Delta_g \mu^{-1})
  \end{BMAT}
\end{bmatrix} \qquad \end{eqnarray}
\begin{align*}
\!\!&          \!  +\!\begin{bmatrix}
\begin{BMAT}(@, 1pt, 1pt){c.c}{c.c}
  \! \!\left[\!\!\!\left. \begin{array}{ll}\! \mu(\mu^{\!-\!1}\!)^{;j} \Gamma_{\!kl}^l \!+\!\mu(\mu\!+\!\tilde{\rho})^{\!\frac{1}{2}}(\mu^{\!-\!1})^{;j} \frac{\partial (\mu\!+\!\tilde{\rho})^{\!-\!\frac{1}{2}}}{\partial x_k}\\
 \! \!+\! \big( g^{ml} \frac{\partial \Gamma_{\!kl}^j}{\partial x_m}\! + \!g^{ml}\Gamma_{\!hl}^j \Gamma_{\!km}^h \!\!-\! g^{ml} \Gamma_{\!kh}^j \Gamma_{\!ml}^h \big)\\
   \frac{\tilde{\rho}}{\mu(\mu\!+\!\tilde{\rho})} g^{ml} \Gamma_{\!km}^j \frac{\partial \mu}{\partial x_l} \!+\!R^j_k \!-\! \mu^{\!
-\!1}\big( \!\frac{\partial \mu}{\partial x_k}\big)^{;j}\end{array} \!\!
\right.\!\! \! \right]_{\!n\!\times \!n}   &   \! \left[\!(\mu\!+\!\tilde{\rho})^{\!\frac{1}{2}}\!\!\left\{\!\!\! \begin{array}{ll} \!\!-2 ( \Delta_g \mu^{\!-\!1})^{;j} - 2  R_l^j \,(  \mu^{\!-1})^{;l}  \\ \!\!-2  g^{ml} \Gamma_{\!sm}^j \frac{\partial((\mu^{\!-\!1})^{;s})}{\partial x_l}
  -2 \mu^{\!-\!1}\frac{\partial \mu}{\partial x_m} (\mu^{\!-\!1}\!)^{;m;j}\\
    \!\!- \!\big(\! g^{ml} \frac{\partial \Gamma^j_{\!sl}}{\partial x_m}\! + \!g^{ml} \Gamma^j_{\!hl}\Gamma^h_{\!sm} \!\!-\! g^{ml} \Gamma_{\!sh}^j\Gamma^h_{\!ml}\! \big) \!(u^{\!-\!1}\!)^{;s}\end{array}\!\!\!\!\!\right\}\!\!\right]_{\!n\!\times \!1}\\
   \left[\mu(\mu+\tilde{\rho})^{-\frac{1}{2}} \Gamma_{kl}^l + \mu \frac{\partial (\mu+\tilde{\rho})^{-\frac{1}{2}}}{\partial x_k} \! \right]_{\!1\times n}&  0
  \end{BMAT}
\!\end{bmatrix}.\nonumber.\end{align*}
 Note that for any smooth $(n+1)$-dimensional vector-valued function $v$,  \begin{eqnarray*} \left( \big(\frac{\partial}{\partial x_n} I_{n+1}\big)q - q\big(\frac{\partial}{\partial x_n}I_{n+1}\big) \right)v\!\!&\!\!=\!\!\!&
 \big(\frac{\partial}{\partial x_n} I_{n+1}\big)\big(qv) - q\big(\frac{\partial}{\partial x_n}I_{n+1}\big) v\\
 \!\!&=\!\!\!& \bigg(\frac{\partial q}{\partial x_n}\bigg) v + q \big(\frac{\partial }{\partial x_n} I_{n+1}\big) v - q \big(\frac{\partial }{\partial x_n} I_{n+1}\big)
 v= \bigg(\frac{\partial q}{\partial x_n}\bigg) v.\nonumber\end{eqnarray*}
 This implies that  \begin{eqnarray} \label{19.3.19-8} \big(\frac{\partial}{\partial x_n} I_{n+1}\big)q - q\big(\frac{\partial}{\partial x_n}I_{n+1}\big)
= \frac{\partial q}{\partial x_n},\end{eqnarray}
i.e., the symbol of $ (\frac{\partial }{\partial x_n} I_{n+1})Q - Q(\frac{\partial }{\partial x_n} I_{n+1})$ is $\frac{\partial q}{\partial x_n}$.
     Combining this, the left-hand side of (\ref{200425-8}) and symbol formula for product of two pseudodifferential operators (see p.$\,37$ of \cite{Tre}, p.$\,$13 of \cite{Ta2} or \cite{KN})\, we get the full symbol equation:
   \begin{eqnarray} \label{19.3.19-4} \sum_{\vartheta} \frac{(-i)^{|\vartheta|}}{\vartheta!} \big(\partial^{\vartheta}_{\xi'} q\big)\big(\partial^\vartheta_{x'}q\big)
 - \sum_{\vartheta} \frac{(-i)^{|\vartheta|}}{\vartheta!} \big( \partial^\vartheta_{\xi'}b\big) \big(\partial^\vartheta_{x'} q\big) -\frac{\partial q}{\partial x_n} +c=0,\end{eqnarray}
 where $\partial_{x'}^{\vartheta}= \frac{\partial^{|\vartheta|}}{\partial x_1^{\vartheta_1}\cdots \partial x_{n-1}^{\vartheta_{n-1}}}$, $\partial_{\xi'}^{\vartheta}= \frac{\partial^{|\vartheta|}}{\partial \xi_1^{\vartheta_1} \cdots \partial\xi_{n-1}^{\vartheta_{n-1}}}$, and $\vartheta=(\vartheta_1,\cdots, \vartheta_{n-1})$ is a $(n\!-\!1)$-tuple of nonnegative integers.

\vskip 0.2 true cm

Step 2. \   Group the homogeneous terms of degree two in (\ref{19.3.19-4}) we obtain the matrix equation \begin{eqnarray}\label{20200531-1} q_1^2 -b_1q_1 +c_2=0,\end{eqnarray}  i.e.,
\begin{align} \label{200425-14} \,\;  \; q^2\!- \!\!\begin{bmatrix}
\begin{BMAT}(@, 1pt, 1pt){c.c}{c.c}
   \!\left[0\right]_{n\times n}   &   \! \left[4\text{i} \mu^{\!-\!1} (\mu\!+\!\tilde{\rho})^{\!\frac{1}{2}} \Gamma_{\beta n}^j g^{\beta \alpha} \xi_\alpha\right]_{\!n\!\times \!1}\\
   \left[0 \right]_{\!1\!\times \!n}&  0  \end{BMAT}\!
\end{bmatrix}\!q_1 \!-\!g^{\alpha\beta} \xi_\alpha\xi_\beta  I_{\!n\!+\!1}\! -\!
 \begin{bmatrix}
\begin{BMAT}(@, 1pt, 1pt){c.c}{c.c}
   \!\left[0\right]_{\!n\!\times \!n}   &   \! \left[2 \mu^{\!-\!1} (\mu\!+\!\tilde{\rho})^{\!\frac{1}{2}} \Gamma_{\gamma \sigma}^j g^{\alpha \gamma} g^{\beta \sigma} \xi_\alpha\xi_\beta\right]_{\!n\!\times \!1}\\
   \left[0 \right]_{\!1\!\times \!n}&  0  \end{BMAT}\!
\end{bmatrix}\!\!=\!0.\end{align}
Our aim is to calculate the unknown $q_1$ by solving the matrix equation (\ref{20200531-1}) (i.e., (\ref{200425-14})). Generally, it is impossible to obtain an exact solution for a quadratic matrix equation. However, by observing the coefficient matrices of equation (\ref{200425-14}), we see that
 the following two matrices play a key role
\begin{eqnarray*} \sqrt{\sum_{\alpha, \beta} g^{\alpha\beta} \xi_\alpha\xi_\beta}\,I_{n+1}, \quad \,   \begin{small} \left[\begin{BMAT}(@, 0pt, 0pt){c.c}{c.c} [0]_{n\times n} & \begin{small}\left[ \Gamma_{\beta n}^j g^{\alpha\beta}  \xi_\alpha \right]_{n\times 1}\end{small} \\
  \left[0\right]_{1\times n} & 0  \end{BMAT}\right]\end{small}.\end{eqnarray*}
The set $F$ of above two matrices can generate a matrix ring $\mathfrak{F}$ according to the following two operations: we first define a multiplication operation between the ring $C^\infty(\Omega\times {\mathbb{R}}^{n-1})$ of all functions and $\mathfrak{F}$: for every $s\in C^\infty(\Omega \times {\mathbb{R}}^{n-1})$ and $A\in \mathfrak{F}$, we have $sA\in \mathfrak{F}$ by the usual multiplication (The element of $C^\infty(\Omega\times {\mathbb{R}}^{n-1})$ is said to be the ``coefficient'' of matrix ring $\mathfrak{F}$); we then define the addition and multiplication by using the usual matrix addition and multiplication of $\mathfrak{F}$. Clearly, $\mathfrak{F}$ is a two-dimensional matrix ring on ``coefficients'' $C^\infty(\Omega\times {\mathbb{R}}^{n-1})$, and  $F$ is a basis of $\mathfrak{F}$.
  This implies that the solution $q_1$ of  equation (\ref{200425-14}) must have the following form:
\begin{eqnarray}  q_1= \begin{bmatrix}
\begin{BMAT}(@, 15pt, 15pt){c.c}{c.c}
   \!\left[\delta_{jk} d(x,\xi') \right]_{n\times n}   &   \! \left[d_j(x,\xi') \right]_{n\times 1}\\
   \left[0 \right]_{1\times n}&  d(x,\xi')  \end{BMAT}
\end{bmatrix}, \end{eqnarray}
where $d(x,\xi')$ and $d_j(x,\xi')$ are the symbol of order $1$. The above idea is similar to Galois group theory for solving the polynomial equation (see \cite{Art} or \cite{HME}) and was recently established by the author of this paper in \cite{Liu1} for solving an elastic inverse problem.
It is clear that for such $q_1$ we have
\begin{eqnarray*} \begin{bmatrix}
\begin{BMAT}(@, 15pt, 15pt){c.c}{c.c}\!
   \!\left[0\right]_{n\times n}   &   \! \left[4\text{i}  \mu^{\!-\!1}\! (\mu\!+\!\tilde{\rho})^{\!\frac{1}{2}} \Gamma_{\beta n}^j g^{\beta \alpha} \xi_\alpha\right]_{n\times 1}\\
   \!\left[0 \right]_{1\times n}&  0 \! \end{BMAT}\!
\end{bmatrix} \! \begin{bmatrix}
\begin{BMAT}(@, 15pt, 15pt){c.c}{c.c}
   \!\left[\delta_{jk} d(x,\xi') \right]_{n\times n}   &   \! \left[d_j(x,\xi') \right]_{n\times 1}\\
   \left[0 \right]_{1\times n}&  d(x,\xi')  \!\end{BMAT}\!
\end{bmatrix}\!\!=\!\! \begin{bmatrix}
\begin{BMAT}(@, 15pt, 15pt){c.c}{c.c}
   \!\left[0\right]_{n\times n}   &   \! \left[ 4\text{i} \,d(x,\xi')  \mu^{\!-\!1} (\mu\!+\!\tilde{\rho})^{\!\frac{1}{2}} \Gamma_{\beta n}^j g^{\beta \alpha} \xi_\alpha\right]_{n\times 1}\\
   \left[0 \right]_{1\times n}&  0  \!\end{BMAT}\!
\end{bmatrix}
\end{eqnarray*}
and
\begin{eqnarray*} \begin{bmatrix}
\begin{BMAT}(@, 15pt, 15pt){c.c}{c.c}
   \!\left[\delta_{jk} d(x,\xi') \right]_{n\times n}   &   \! \left[d_j(x,\xi') \right]_{n\times 1}\\
   \left[0 \right]_{1\times n}&  d(x,\xi')  \!\!\end{BMAT}\!
\end{bmatrix}\!\begin{bmatrix}
\begin{BMAT}(@, 15pt, 15pt){c.c}{c.c}
   \!\left[0\right]_{n\times n}   &   \! \left[4\text{i}  \mu^{\!-\!1} (\mu\!+\!\tilde{\rho})^{\!\frac{1}{2}} \Gamma_{\beta n}^j g^{\beta \alpha} \xi_\alpha\right]_{n\times 1}\\
  \! \left[0 \right]_{1\times n}&  0 \! \end{BMAT}\!
\end{bmatrix} \!
\!\!=\!\! \begin{bmatrix}
\begin{BMAT}(@, 15pt, 15pt){c.c}{c.c}
   \!\left[0\right]_{n\times n}   &   \! \left[ 4\text{i} \,d(x,\xi')  \mu^{\!-\!1} (\mu\!+\!\tilde{\rho})^{\!\frac{1}{2}} \Gamma_{\beta n}^j g^{\beta \alpha} \xi_\alpha\right]_{n\times 1}\\
   \left[0 \right]_{1\times n}&  0  \!\end{BMAT}\!
\end{bmatrix},
\end{eqnarray*}
i.e.,
$b_1 q_1=q_1b_1$. Thus the matrix equation (\ref{20200531-1}) becomes the following equivalent matrix equation
\begin{eqnarray} \label{200426-16} q_1^2 -\frac{1}{2} b_1 q_1  -\frac{1}{2} q_1 b_1+c_2 =0.\end{eqnarray}
We can easily verify that (\ref{200426-16}) has the following two solutions
\begin{eqnarray} \label{200426-17} q_1= \frac{1}{2} \left( b_1\pm \sqrt{ b^2_1 -4c_2}\right) .\end{eqnarray}
In view of $b_1^2 =0$, we get
\begin{align} \label{200426-171} & q_1\!=\! \frac{1}{2}\!  \begin{bmatrix}\!
\begin{BMAT}(@, 15pt, 15pt){c.c}{c.c}
   \!\left[0\right]_{\!n\!\times\! n}   &   \! \left[4\text{i} \mu^{\!-\!1}\! (\mu\!+\!\tilde{\rho})^{\!\frac{1}{2}} \Gamma_{\beta n}^j g^{\beta \alpha} \xi_\alpha\right]_{\!n\!\times \!1}\\
   \left[0 \right]_{\!1\!\times \!n}&  0  \end{BMAT}
\end{bmatrix} \!\pm \frac{1}{2}\!
\sqrt{\!4 \!\left(\! g^{\alpha\beta} \xi_\alpha\xi_\beta \, I_{n+1}  \!+ \! \begin{bmatrix}
\begin{BMAT}(@, 15pt, 15pt){c.c}{c.c}
   \!\left[0\right]_{\!n\!\times \!n}   &   \! \left[2 \mu^{\!-\!1}\! (\mu\!+\!\tilde{\rho})^{\!\frac{1}{2}} \Gamma_{\gamma \sigma}^j g^{\alpha \gamma} g^{\beta \sigma} \xi_\alpha\xi_\beta\!\right]_{\!n\!\times \!1}\\
   \left[0 \right]_{\!1\!\times\! n}&  0  \end{BMAT}
\end{bmatrix}\!\right) } \nonumber\\
&  \quad \!= 2\text{i} \mu^{\!-\!1}\! (\mu\!+\!\tilde{\rho})^{\!\frac{1}{2}} \!\begin{bmatrix}
\begin{BMAT}(@, 1pt, 1pt){c.c}{c.c}
   \!\left[0\right]_{n\times n}   &   \! \left[ \Gamma_{\!\beta n}^j g^{\beta \alpha}\xi_\alpha\right]_{\!n\times 1}\\
   \left[0 \right]_{\!1\times n}&  0  \end{BMAT}\!
\end{bmatrix} \!\pm \!\sqrt{g^{\alpha\beta} \xi_\alpha \xi_\beta} \; \sqrt{ I_{n+1}\! +\! 2 \mu^{\!-\!1}\! (\mu\!+\!\tilde{\rho})^{\!\frac{1}{2}} \!\!
\begin{bmatrix}
\begin{BMAT}(@, 5pt, 5pt){c.c}{c.c}
   \!\left[0\right]_{n\times n}   &   \! \left[\frac{\Gamma_{\gamma \sigma}^j g^{\alpha \gamma} g^{\beta \sigma} \xi_\alpha\xi_\beta}{ g^{\alpha\beta}\xi_\alpha\xi_\beta }\right]_{n\times 1}\\
   \left[0 \right]_{1\times n}&  0  \end{BMAT}\!
\end{bmatrix}}.\nonumber  \end{align}
Since $$ I_{n+1} + 2 \mu^{\!-\!1}\! (\mu\!+\!\tilde{\rho})^{\!\frac{1}{2}} \!
\begin{bmatrix}
\begin{BMAT}(@, 15pt, 15pt){c.c}{c.c}
   \!\left[0\right]_{n\times n}   &   \! \left[\frac{\Gamma_{\gamma \sigma}^j g^{\alpha \gamma} g^{\beta \sigma} \xi_\alpha\xi_\beta}{ g^{\alpha\beta}\xi_\alpha\xi_\beta }\right]_{n\times 1}\\
   \left[0 \right]_{1\times n}&  0  \end{BMAT}
\end{bmatrix}$$ is  a positive-definite matrix and
\begin{align*} & \left(\!I_{n+1} \!+\! \mu^{\!-\!1}\! (\mu\!+\!\tilde{\rho})^{\!\frac{1}{2}} \!\! \begin{bmatrix}
\begin{BMAT}(@, 1pt, 1pt){c.c}{c.c}
   \!\left[0\right]_{n\times n}   &   \! \left[\!\frac{\Gamma_{\!\gamma \sigma}^j g^{\alpha \gamma} g^{\beta \sigma} \xi_\alpha\xi_\beta}{ g^{\alpha\beta}\xi_\alpha\xi_\beta }\!\right]_{n\times 1}\\
   \left[0 \right]_{1\times n}&  0  \end{BMAT}
\end{bmatrix}\!\right) \!\left(\! I_{n+1} \!+\! \mu^{\!-\!1}\! (\mu\!+\!\tilde{\rho})^{\!\frac{1}{2}} \!\! \begin{bmatrix}
\begin{BMAT}(@, 1pt, 1pt){c.c}{c.c}
   \!\left[0\right]_{n\times n}   &   \! \left[\!\frac{\Gamma_{\!\gamma \sigma}^j g^{\alpha \gamma} g^{\beta \sigma} \xi_\alpha\xi_\beta}{ g^{\alpha\beta}\xi_\alpha\xi_\beta }\right]_{\!n\!\times \!1}\\
   \left[0 \right]_{1\times n}&  0  \end{BMAT}
\end{bmatrix}\!\right)\\
&=   I_{n+1} +2 \mu^{\!-\!1} (\mu\!+\!\tilde{\rho})^{\!\frac{1}{2}} \!  \begin{bmatrix}
\begin{BMAT}(@, 15pt, 15pt){c.c}{c.c}
   \!\left[0\right]_{n\times n}   &   \! \left[\frac{\Gamma_{\gamma \sigma}^j g^{\alpha \gamma} g^{\beta \sigma} \xi_\alpha\xi_\beta}{ g^{\alpha\beta}\xi_\alpha\xi_\beta }\right]_{n\times 1}\\
   \left[0 \right]_{1\times n}&  0  \end{BMAT}
\end{bmatrix},\nonumber\end{align*}
we have \begin{eqnarray} & q_1 = 2 \text{i}\,  \mu^{\!-\!1} (\mu\!+\!\tilde{\rho})^{\!\frac{1}{2}} \!
 \begin{bmatrix}
\begin{BMAT}(@, 15pt, 15pt){c.c}{c.c}
   \!\left[0\right]_{n\times n}   &   \! \left[\!  \Gamma_{\beta n}^j g^{\beta \alpha}\xi_\alpha\right]_{\!n\times 1}\\
   \left[0 \right]_{1\times n}&  0  \end{BMAT}
\end{bmatrix} \!\pm\! \sqrt{ g^{\alpha\beta} \xi_\alpha\xi_\beta} \!\left(\! \!I_{n+1} \!+\! \mu^{\!-\!1} (\mu\!+\!\tilde{\rho})^{\!\frac{1}{2}} \! \begin{bmatrix}
\begin{BMAT}(@, 15pt, 15pt){c.c}{c.c}
   \!\left[0\right]_{n\times n}   &   \! \left[\frac{\Gamma_{\gamma \sigma}^j g^{\alpha \gamma} g^{\beta \sigma} \xi_\alpha\xi_\beta}{ g^{\alpha\beta}\xi_\alpha\xi_\beta }\right]_{n\times 1}\\
   \left[0 \right]_{1\times n}&  0  \end{BMAT}
\end{bmatrix}\!\right)\nonumber \\
 \!\!\!\!\!\! & \!\!\! =\pm\sqrt{ g^{\alpha\beta} \xi_\alpha\xi_\beta} \, I_{n+1}\!+\!  \mu^{\!-\!1} (\mu\!+\!\tilde{\rho})^{\!\frac{1}{2}} \! \begin{bmatrix}
\begin{BMAT}(@, 15pt, 15pt){c.c}{c.c}
   \!\left[0\right]_{n\times n}   &   \! \left[ 2\text{i}\, \Gamma_{\beta n}^j  g^{\beta \alpha}\xi_\alpha\!\pm\!
   \frac{\Gamma_{\gamma \sigma}^j g^{\alpha \gamma} g^{\beta \sigma} \xi_\alpha\xi_\beta}{\sqrt{ g^{\alpha\beta}\xi_\alpha\xi_\beta }}\right]_{\!n\times 1}\\
   \left[0 \right]_{1\times n}&  0  \end{BMAT}
\end{bmatrix}.\qquad \qquad \qquad \quad \quad \;\;\,\quad \nonumber\end{eqnarray}
Because we have chosen the outward normal $\nu$ of $\partial \Omega$, we should take
\begin{eqnarray}\label{20200516-1}  q_1\! =\sqrt{ g^{\alpha\beta} \xi_\alpha\xi_\beta}\,  I_{n+1}\!+\!  \mu^{\!-\!1}\! (\mu\!+\!\tilde{\rho})^{\!\frac{1}{2}} \! \begin{bmatrix}
\begin{BMAT}(@, 15pt, 15pt){c.c}{c.c}
   \!\left[0\right]_{n\times n}   &   \! \left[ 2\text{i}\, \Gamma_{\beta n}^j  g^{\beta \alpha}\xi_\alpha+
   \frac{\Gamma_{\gamma \sigma}^j g^{\alpha \gamma} g^{\beta \sigma} \xi_\alpha\xi_\beta}{\sqrt{ g^{\alpha\beta}\xi_\alpha\xi_\beta }}\right]_{\!n\times 1}\\
   \left[0 \right]_{1\times n}&  0  \end{BMAT}
\end{bmatrix},\end{eqnarray}
which is a positive-definite matrix.
\vskip 0.2 true cm

Step 3. \  The terms of degree one in (\ref{19.3.19-4}) are
\begin{eqnarray}  \label{19.3.22-10} \quad  q_1 q_0 +q_0q_1  -i\sum_{l=1}^{n-1}  \frac{\partial q_1}{\partial \xi_l} \, \frac{\partial q_1}{\partial x_l}  - b_1 q_0 - b_0 q_1 -\frac{1}{i} \sum_{l=1}^{n-1} \frac{\partial b_1}{\partial \xi_l} \, \frac{\partial q_1}{\partial x_l}   -\frac{\partial q_1}{\partial x_n} +c_1=0, \end{eqnarray}
i.e., \begin{eqnarray} \label{19.4.17-1} (q_1- b_1) q_0 +q_0 q_1= E_1,\end{eqnarray}
where \begin{eqnarray} \label{19.3.23-1} E_1:= i\sum_{l=1}^{n-1} \frac{\partial q_1}{\partial \xi_l} \, \frac{\partial q_1}{\partial x_l} + b_0q_1 -i \sum_{l=1}^{n-1} \frac{\partial b_1}{\partial \xi_l}\, \frac{\partial q_1}{\partial x_l} +\frac{\partial q_1}{\partial x_n}  -c_1,\end{eqnarray}
 and $b_0$ and $c_1$ are given in (\ref{19.3.19-5}) and (\ref{20200515-6}).
More precisely, \begin{eqnarray}  \label{19.3.22-11} &&  \begin{small} \left(\sqrt{ g^{\alpha\beta} \xi_\alpha\xi_\beta}  I_{n+1}\!+\! \mu^{\!-\!1}\! (\mu\!+\!\tilde{\rho})^{\!\frac{1}{2}} \! \begin{bmatrix}
\begin{BMAT}(@, 5pt, 5pt){c.c}{c.c}
   \!\left[0\right]_{n\times n}   &   \! \left[ -2\text{i}\, \Gamma_{\beta n}^j  g^{\beta \alpha}\xi_\alpha\!+\!
   \frac{1}{\sqrt{ g^{\alpha\beta}\xi_\alpha\xi_\beta }}\Gamma_{\gamma \sigma}^j g^{\alpha \gamma} g^{\beta \sigma} \xi_\alpha\xi_\beta\right]_{\!n\times 1}\\
   \left[0 \right]_{1\times n}&  0  \end{BMAT}
\end{bmatrix} \right)\end{small}q_0 \\
 && +
 q_0 \begin{small} \left(\sqrt{ g^{\alpha\beta} \xi_\alpha\xi_\beta}  I_{n+1}\!+\!  \mu^{\!-\!1}\! (\mu\!+\!\tilde{\rho})^{\!\frac{1}{2}} \! \begin{bmatrix}
\begin{BMAT}(@, 15pt, 15pt){c.c}{c.c}
   \!\left[0\right]_{n\times n}   &   \! \left[ 2\text{i}\,  \Gamma_{\beta n}^jg^{\beta \alpha} \xi_\alpha\!+\!
   \frac{1}{\sqrt{ g^{\alpha\beta}\xi_\alpha\xi_\beta }}\Gamma_{\gamma \sigma}^j g^{\alpha \gamma} g^{\beta \sigma} \xi_\alpha\xi_\beta\right]_{\!n\times 1}\\
   \left[0 \right]_{1\times n}&  0  \end{BMAT}
\end{bmatrix} \right)\end{small}=E_1.\nonumber \end{eqnarray}
   Now, we calculate $q_0$ by solving Sylvester's matrix equation (\ref{19.3.22-11}).
It is well-known that Sylvester's equation of the form
$LX+XM=E$ can be put into the form (see \cite{BaS} or \cite{BhR})
\begin{eqnarray} \label{19.10.6-1} U(\mbox{vec}\,X)= \mbox{V}\end{eqnarray}  for larger matrices $U$ and $V$.
 Here  $\mbox{vec}\, X$ is a stack of all columns of matrix $X$ (see, for example,  Chapter 4 of \cite{HoJ}).
Indeed, $U=(I_{n+1}\otimes L)+(M^t \otimes I_{n+1})$, and $V=\mbox{vec}\,E$, where $\otimes$ denotes the Kronecker product.
Thus, if we can obtain the inverse $U^{-1}$ of the matrix $U$, then we have $\mbox{vec}\,X=U^{-1} (\mbox{vec}\,V)$, and the corresponding solution $X$ will immediately be obtained.
From (\ref{19.3.22-11}) we see that
 \begin{eqnarray*}  \label{19.3.22-12} L:=& \begin{small} \!\!\left(\!\sqrt{ g^{\alpha\beta} \xi_\alpha\xi_\beta} \, I_{n+1}\!+\!  \mu^{\!-\!1} (\mu\!+\!\tilde{\rho})^{\!\frac{1}{2}} \begin{bmatrix}
\begin{BMAT}(@, 15pt, 15pt){c.c}{c.c}
   \!\left[0\right]_{n\times n}   &   \! \left[ -2\text{i}\,\Gamma_{\beta n}^j g^{\beta \alpha}  \xi_\alpha\!+\!
   \frac{1}{\sqrt{ g^{\alpha\beta}\xi_\alpha\xi_\beta }}\Gamma_{\gamma \sigma}^j g^{\alpha \gamma} g^{\beta \sigma} \xi_\alpha\xi_\beta\right]_{\!n\times 1}\\
   \left[0 \right]_{1\times n}&  0  \end{BMAT}
\end{bmatrix} \right)\end{small}\end{eqnarray*} and
 \begin{eqnarray*}  M^t :=\begin{small} \left(\sqrt{ g^{\alpha\beta} \xi_\alpha\xi_\beta} \, I_{n+1}\!+\!  \mu^{\!-\!1} (\mu\!+\!\tilde{\rho})^{\!\frac{1}{2}} \begin{bmatrix}
\begin{BMAT}(@, 15pt, 15pt){c.c}{c.c}
   \!\left[0\right]_{n\times n}   &   \! \left[0 \right]_{n\times 1}\\
   \left[2\text{i}\, \Gamma_{\beta n}^k  g^{\beta \alpha}\xi_\alpha\!+\!
   \frac{1}{\sqrt{ g^{\alpha\beta}\xi_\alpha\xi_\beta }}\Gamma_{\gamma \sigma}^k g^{\alpha \gamma} g^{\beta \sigma} \xi_\alpha\xi_\beta \right]_{1\times n}&  0  \end{BMAT}
\end{bmatrix} \right)\end{small}\nonumber \end{eqnarray*}
Thus the  $U$ has the form:
\begin{eqnarray} \label{19.3.23-15} &&\,\; U = (I_n\otimes L) + (M^t \otimes I_n) \\
&& \;\quad \;= 2 \sqrt{g^{\alpha\beta} \xi_\alpha\xi_\beta} \, I_{(n+1)^2} +  \mu^{\!-\!1} (\mu\!+\!\tilde{\rho})^{\frac{1}{2}}\,  (I_{n+1}\otimes A_1) + \mu^{\!-\!1} (\mu\!+\!\tilde{\rho})^{\frac{1}{2}}\,( A_2^t \otimes I_{n+1}),\nonumber\end{eqnarray}
where \begin{eqnarray} \label{19.3.24-1} A_1 &=&  \begin{small} \begin{bmatrix}
\begin{BMAT}(@, 15pt, 15pt){c.c}{c.c}
   \!\left[0\right]_{n\times n}   &   \! \left[ -2\text{i}\,  \Gamma_{\beta n}^jg^{\beta \alpha} \xi_\alpha\!+\!
   \frac{1}{\sqrt{ g^{\alpha\beta}\xi_\alpha\xi_\beta }}\Gamma_{\gamma \sigma}^j g^{\alpha \gamma} g^{\beta \sigma} \xi_\alpha\xi_\beta\right]_{\!n\times 1}\\
   \left[0 \right]_{1\times n}&  0  \end{BMAT}
\end{bmatrix} \end{small} \\
    A_2^t &=&\begin{small}\begin{bmatrix}
\begin{BMAT}(@, 15pt, 15pt){c.c}{c.c}
   \!\left[0\right]_{n\times n}   &   \! \left[0 \right]_{n\times 1}\\
   \left[2\text{i}\,  \Gamma_{\beta n}^k g^{\beta \alpha}\xi_\alpha\!+\!
   \frac{1}{\sqrt{ g^{\alpha\beta}\xi_\alpha\xi_\beta }}\Gamma_{\gamma \sigma}^k g^{\alpha \gamma} g^{\beta \sigma} \xi_\alpha\xi_\beta \right]_{1\times n}&  0  \end{BMAT}
\end{bmatrix}.\end{small}  \end{eqnarray}
Setting $\Upsilon_j=2\;\!\text{i}\;\! \Gamma_{\beta n}^j  g^{\alpha\beta}\xi_\alpha$ and $\Theta_j=  \Gamma_{\gamma\sigma}^j g^{\alpha\gamma} g^{\beta\sigma} \frac{\xi_\alpha\xi_\beta}{\sqrt{g^{\alpha\beta} \xi_\alpha\xi_\beta}}$, we see that
\begin{align*} \!\!\!& \!\!
         \begin{small} \begin{bmatrix} \begin{small} \left[\!\begin{BMAT}(@, 10pt, 16pt){c.c}{c.c} \left[0\right]_{n\times n}& \left[-\Upsilon_j
+ \Theta_j \right]_{n\times 1}\\
\left[0\right]_{1\times n} & 0\end{BMAT} \right] \end{small}
 \\
 \!\!&\! \!\ddots\! \\
\! \!&\!&\!  \! \begin{small} \left[\!\begin{BMAT}(@, 10pt, 1pt){c.c}{c.c} \left[0\right]_{n\times n}& \left[-\Upsilon_j
+ \Theta_j \right]_{n\times 1}\\
\left[0\right]_{1\times n} & 0\end{BMAT} \right] \end{small}\!\\
  \!&\!&\!& \!
   \begin{small} \left[\!\begin{BMAT}(@, 10pt, 16pt){c.c}{c.c} \left[0\right]_{n\times n}& \left[-\Upsilon_j
+ \Theta_j \right]_{n\times 1}\\
\left[0\right]_{1\times n} & 0\end{BMAT} \right] \end{small}
\end{bmatrix}\end{small}\!
 \\
  & \quad\; \times
  \begin{small} \begin{bmatrix}\!\!\begin{small} \begin{bmatrix}0\!\!\\
\! \!&\! \!\ddots \!\!&\!\!&\! \\
 \!\! &\!\!& \!\!0\!\!\!\!\!\end{bmatrix}_{\!(\!n\!+\!1)\!\times \!(\!n\!+\!1\!)}\end{small}\!
  \!\! & \!\!\cdots \!\!& \!\!\!\begin{small} \begin{bmatrix}0\!\!\\
\! \!&\! \!\ddots \!\!&\!\!&\! \\
 \!\! &\!\!& \!\!0\!\!\!\!\!\end{bmatrix}_{\!(\!n\!+\!1)\!\times \!(\!n\!+\!1\!)}\end{small} \!\! &\! \!\begin{small} \begin{bmatrix}0\!\!\\
\! \!&\! \!\ddots \!\!&\!\!&\! \\
 \!\! &\!\!& \!\!\!0\!\!\!\!\!\!\end{bmatrix}_{\!(\!n\!+\!1)\!\times \!(\!n\!+\!1\!)}\end{small}\! \\
\!\!\!\!\!\!\!\!\!\!\!\!\!\! \vdots \!\!\!\!\!\!\!\!\!& \!\!\!\!\!\!\!\! \!\!\!& \!\!\!\!\!\!\!\!\!\!\!\vdots\!\!\!\!\!\!\! \!\!\!\!&\!\!\!\! \!\!\!\!\!\!\!\vdots \!\\
\!\begin{small} \begin{bmatrix}0\!\!\\
\! \!&\! \!\ddots \!\!&\!\!&\! \\
 \!\! &\!\!& \!\!0\!\!\!\!\end{bmatrix}_{\!(\!n\!+\!1)\!\times \!(\!n\!+\!1\!)}\end{small} \!\!&\cdots \!\!& \!\!\begin{small} \begin{bmatrix}\!0\!\!\\
\! \!&\! \!\ddots \!\!&\!\!&\! \\
 \!\! &\!\!& \!\!0\!\!\!\!\end{bmatrix}_{\!(\!n\!+\!1)\!\times \!(\!n\!+\!1\!)}\end{small}\! \!&\!\! \begin{small} \begin{bmatrix}0\!\!\\
\! \!&\! \!\ddots \!\!&\!\!&\! \\
 \!\! &\!\!& \!\!\!0\!\!\!\!\!\end{bmatrix}_{\!(\!n\!+\!1)\!\times \!(\!n\!+\!1\!)}\end{small}\\
\!\begin{small} \begin{bmatrix}\Upsilon_1
+ \Theta_1\!\\
\! \!\!\!\!\!\!&\!\!\! \!\!\!\!\ddots\!\!\!\! \!\!\!&\!\!\!\!\!\!&\!\!\! \\
 \! \!\!&\!\!\!&\!\! \!\!\Upsilon_1
+ \Theta_1\!\!\end{bmatrix}_{\!(\!n\!+\!1)\!\times \!(\!n\!+\!1\!)}\end{small}\!\!\!\!\!&\!\!\!\!\!\cdots\!\!\!\! \!\!&\!\! \!\!\!\begin{small} \begin{bmatrix}\Upsilon_n
+ \Theta_n\!\\
\!\! \!\!\!\!\!&\!\!\!\! \!\!\ddots\! \!\!\!\!\!\!\!&\!\!\!\!\!\!\!&\!\! \\
 \!\!\!\! \!\!&\!\!\!\!\!\!&\!\!\!\! \!\!\Upsilon_n
+ \Theta_n\!\end{bmatrix}_{\!(\!n\!+\!1)\!\times \!(\!n\!+\!1\!)}\end{small} \!\!& \!\!\begin{small} \begin{bmatrix}0\!\!\\
\! \!&\! \!\ddots \!\!&\!\!&\! \\
 \!\! &\!\!& \!\!0\!\!\!\!\end{bmatrix}_{\!(\!n\!+\!1)\!\times \!(\!n\!+\!1\!)}\end{small}\end{bmatrix} \end{small}\! \\
 \!\!\!&         = \begin{small} \begin{bmatrix}\!\!\begin{small} \begin{bmatrix}0\!\!\\
\! \!&\! \!\ddots \!\!&\!\!&\! \\
 \!\! &\!\!& \!\!0\!\!\!\!\!\end{bmatrix}_{\!(\!n\!+\!1)\!\times \!(\!n\!+\!1\!)}\end{small}\!
  \!\! & \!\!\cdots \!\!& \!\!\!\begin{small} \begin{bmatrix}0\!\!\\
\! \!&\! \!\ddots \!\!&\!\!&\! \\
 \!\! &\!\!& \!\!0\!\!\!\!\!\end{bmatrix}_{\!(\!n\!+\!1)\!\times \!(\!n\!+\!1\!)}\end{small} \!\! &\! \!\begin{small} \begin{bmatrix}0\!\!\\
\! \!&\! \!\ddots \!\!&\!\!&\! \\
 \!\! &\!\!& \!\!\!0\!\!\!\!\!\!\end{bmatrix}_{\!(\!n\!+\!1)\!\times \!(\!n\!+\!1\!)}\end{small}\! \\
\!\vdots \!\!& \!\! & \!\!\vdots \!\!& \!\!\vdots \!\\
\!\begin{small} \begin{bmatrix}0\!\!\\
\! \!&\! \!\ddots \!\!&\!\!&\! \\
 \!\! &\!\!& \!\!0\!\!\!\!\end{bmatrix}_{\!(\!n\!+\!1)\!\times \!(\!n\!+\!1\!)}\end{small} \!\!&\cdots \!\!& \!\!\begin{small} \begin{bmatrix}\!0\!\!\\
\! \!&\! \!\ddots \!\!&\!\!&\! \\
 \!\! &\!\!& \!\!0\!\!\!\!\end{bmatrix}_{\!(\!n\!+\!1)\!\times \!(\!n\!+\!1\!)}\end{small}\! \!&\!\! \begin{small} \begin{bmatrix}0\!\!\\
\! \!&\! \!\ddots \!\!&\!\!&\! \\
 \!\! &\!\!& \!\!\!0\!\!\!\!\!\end{bmatrix}_{\!(\!n\!+\!1)\!\times \!(\!n\!+\!1\!)}\end{small}\\
\begin{small} \left[\!\begin{BMAT}(@, 10pt, 1pt){c.c}{c.c} \left[0\right]_{n\times n}& \left[(-\Upsilon_j\!
+\! \Theta_j)(\Upsilon_1\!
+\! \Theta_1) \right]_{n\times 1}\\
\left[0\right]_{1\times n} & 0\end{BMAT} \right] \end{small}\!\!\!\!\!&\!\!\!\!\!\cdots\!\!\!\! \!\!&\!\! \!\!\!\begin{small} \left[\!\begin{BMAT}(@, 0pt, 0pt){c.c}{c.c} \left[0\right]_{n\times n}& \left[(-\Upsilon_j
+ \Theta_j)(\Upsilon_n\!
+\! \Theta_n) \right]_{n\times 1}\\
\left[0\right]_{1\times n} & 0\end{BMAT} \right] \end{small} \!\!& \!\!\begin{small} \begin{bmatrix}0\!\!\\
\! \!&\! \!\ddots \!\!&\!\!&\! \\
 \!\! &\!\!& \!\!0\!\!\!\!\end{bmatrix}_{\!(\!n\!+\!1)\!\times \!(\!n\!+\!1\!)}\end{small}\end{bmatrix} \end{small}\! \\
\!\!\! & =
  \begin{small} \begin{bmatrix}\!\!\begin{small} \begin{bmatrix}0\!\!\\
\! \!&\! \!\ddots \!\!&\!\!&\! \\
 \!\! &\!\!& \!\!0\!\!\!\!\!\end{bmatrix}_{\!(\!n\!+\!1)\!\times \!(\!n\!+\!1\!)}\end{small}\!
  \!\! & \!\!\cdots \!\!& \!\!\!\begin{small} \begin{bmatrix}0\!\!\\
\! \!&\! \!\ddots \!\!&\!\!&\! \\
 \!\! &\!\!& \!\!0\!\!\!\!\!\end{bmatrix}_{\!(\!n\!+\!1)\!\times \!(\!n\!+\!1\!)}\end{small} \!\! &\! \!\begin{small} \begin{bmatrix}0\!\!\\
\! \!&\! \!\ddots \!\!&\!\!&\! \\
 \!\! &\!\!& \!\!\!0\!\!\!\!\!\!\end{bmatrix}_{\!(\!n\!+\!1)\!\times \!(\!n\!+\!1\!)}\end{small}\! \\
\!\vdots \!\!& \!\! & \!\!\vdots \!\!& \!\!\vdots \!\\
\!\begin{small} \begin{bmatrix}0\!\!\\
\! \!&\! \!\ddots \!\!&\!\!&\! \\
 \!\! &\!\!& \!\!0\!\!\!\!\end{bmatrix}_{\!(\!n\!+\!1)\!\times \!(\!n\!+\!1\!)}\end{small} \!\!&\cdots \!\!& \!\!\begin{small} \begin{bmatrix}\!0\!\!\\
\! \!&\! \!\ddots \!\!&\!\!&\! \\
 \!\! &\!\!& \!\!0\!\!\!\!\end{bmatrix}_{\!(\!n\!+\!1)\!\times \!(\!n\!+\!1\!)}\end{small}\! \!&\!\! \begin{small} \begin{bmatrix}0\!\!\\
\! \!&\! \!\ddots \!\!&\!\!&\! \\
 \!\! &\!\!& \!\!\!0\!\!\!\!\!\end{bmatrix}_{\!(\!n\!+\!1)\!\times \!(\!n\!+\!1\!)}\end{small}\\
\!\begin{small} \begin{bmatrix}\Upsilon_1
+ \Theta_1\!\\
\! \!\!\!\!\!\!&\!\!\! \!\!\!\!\ddots\!\!\!\! \!\!\!&\!\!\!\!\!\!&\!\!\! \\
 \! \!\!&\!\!\!&\!\! \!\!\Upsilon_1
+ \Theta_1\!\!\end{bmatrix}_{\!(\!n\!+\!1)\!\times \!(\!n\!+\!1\!)}\end{small}\!\!\!\!\!&\!\!\!\!\!\cdots\!\!\!\! \!\!&\!\! \!\!\!\begin{small} \begin{bmatrix}\Upsilon_n
+ \Theta_n\!\\
\!\! \!\!\!\!\!&\!\!\!\! \!\!\ddots\! \!\!\!\!\!\!\!&\!\!\!\!\!\!\!&\!\! \\
 \!\!\!\! \!\!&\!\!\!\!\!\!&\!\!\!\! \!\!\Upsilon_n
+ \Theta_n\!\end{bmatrix}_{\!(\!n\!+\!1)\!\times \!(\!n\!+\!1\!)}\end{small} \!\!& \!\!\begin{small} \begin{bmatrix}0\!\!\\
\! \!&\! \!\ddots \!\!&\!\!&\! \\
 \!\! &\!\!& \!\!0\!\!\!\!\end{bmatrix}_{\!(\!n\!+\!1)\!\times \!(\!n\!+\!1\!)}\end{small}\end{bmatrix} \end{small}\! \\
   \!\!\!&        \quad \;\times
   \begin{small} \begin{bmatrix} \begin{small} \left[\!\begin{BMAT}(@, 0pt, 0pt){c.c}{c.c} \left[0\right]_{n\times n}& \left[-\Upsilon_j
+ \Theta_j \right]_{n\times 1}\\
\left[0\right]_{1\times n} & 0\end{BMAT} \right] \end{small}
 \\
 \!\!&\! \!\ddots\! \\
\! \!&\!&\!  \! \begin{small} \left[\!\begin{BMAT}(@, 10pt, 10pt){c.c}{c.c} \left[0\right]_{n\times n}& \left[-\Upsilon_j
+ \Theta_j \right]_{n\times 1}\\
\left[0\right]_{1\times n} & 0\end{BMAT} \right] \end{small}\!\\
  \!&\!&\!& \!
   \begin{small} \left[\!\begin{BMAT}(@, 0pt, 0pt){c.c}{c.c} \left[0\right]_{n\times n}& \left[-\Upsilon_j
+ \Theta_j \right]_{n\times 1}\\
\left[0\right]_{1\times n} & 0\end{BMAT} \right] \end{small}
\end{bmatrix}\end{small},  \end{align*}
i.e., $$(I_{n+1} \otimes A_1)(A_2^t \otimes I_{n+1}) = (A_2^t \otimes I_{n+1}) (I_{n+1} \otimes A_1).$$
 Also $$(I_{n+1} \otimes A_1)(I_{n+1} \otimes A_1)=0$$ because of $A_1^2=0$, and $(A_2^t \otimes I_{n+1}) (A_2^t \otimes I_{n+1})=0$.
 We find that the following four matrices are linearly independent and generate a matrix ring $\mathfrak{X}$ on the ring $C^\infty (\Omega\times {\mathbb{R}}^{n-1})$ of all functions about the addition and multiply of matrices:
\begin{eqnarray*} \label{19.4.6-p2} &&H=\bigg\{ \sqrt{ g^{\alpha \beta}\; \xi_\alpha \xi_\beta}\; I_{(n+1)^2},\; \quad I_{n+1}\otimes A_1, \quad   A_2^t \otimes I_{n+1}, \quad (I_{n+1} \otimes A_1) (A_2^t \otimes I_{n+1})\bigg\}, \end{eqnarray*}
This implies that we should look for the inverse $U^{-1}$ of the form:
\begin{eqnarray} \label{19.3.26-6} && \;\;\quad U^{-1} \!=\! {\tilde{s}}_1 \sqrt{ g^{\alpha\beta} \xi_\alpha \xi_\beta}\, I_{(n+1)^2}
 +{\tilde{s}}_2(I_{n+1}\otimes A_1) +{\tilde{s}}_3(A_2^t \otimes I_{n+1} )\! +\! {\tilde{s}}_4  \frac{1}{\sqrt{ g^{\alpha\beta} \xi_\alpha \xi_\beta}} \!(I_{n+1}\otimes A_1)(A_2^t \otimes I_{n+1} ), \end{eqnarray}
 where ${\tilde{s}}_1, {\tilde{s}}_2, {\tilde{s}}_3, {\tilde{s}}_4$ are the undetermined functions.
By inserting (\ref{19.3.26-6}) into $UU^{-1}=I_{(n+1)^2}$, we have
\begin{eqnarray*} &\left(2 \sqrt{g^{\alpha\beta} \xi_\alpha\xi_\beta} \, I_{(n+1)^2} +\mu^{\!-\!1} (\mu\!+\!\tilde{\rho})^{\frac{1}{2}} (I_{n+1}\otimes A_1) + \mu^{\!-\!1} (\mu\!+\!\tilde{\rho})^{\frac{1}{2}}( A_2^t \otimes I_{n+1})\right)\left(
 {\tilde{s}}_1 \sqrt{ g^{\alpha\beta} \xi_\alpha \xi_\beta}\, I_{(n+1)^2}
  \right.\\
  & \left.\;\;\;+{\tilde{s}}_2(I_{n+1}\otimes A_1)+{\tilde{s}}_3(A_2^t \otimes I_{n+1} )\! +\! {\tilde{s}}_4  \frac{1}{\sqrt{ g^{\alpha\beta} \xi_\alpha \xi_\beta}} \!(I_{n+1}\otimes A_1)(A_2^t \otimes I_{n+1} )
\right)=I_{(n+1)^2}, \end{eqnarray*}
i.e.,
\begin{eqnarray} \label{200427-1,} &  2{\tilde{s}}_1 g^{\alpha\beta} \xi_\alpha\xi_\beta \, I_{(n+1)^2}  \!+\! (2 {\tilde{s}}_2 \!+\!{\tilde{s}}_1 \mu^{\!-\!1} (\mu\!+\!\tilde{\rho})^{\frac{1}{2}} ) \sqrt{ g^{\alpha\beta} \xi_\alpha\xi_\beta} (I_{n+1} \otimes A_1) \! +\! (2{\tilde{s}}_3 \!+\! {\tilde{s}}_1 \mu^{\!-\!1} (\mu\!+\!\tilde{\rho})^{\frac{1}{2}}) \sqrt{ g^{\alpha\beta} \xi_\alpha \xi_\beta} (A_2^t \otimes  I_{n+1})\nonumber\\
 & + (2{\tilde{s}}_4 + {\tilde{s}}_2 \mu^{\!-\!1} (\mu\!+\!\tilde{\rho})^{\frac{1}{2}} + {\tilde{s}}_3 \mu^{\!-\!1} (\mu\!+\!\tilde{\rho})^{\frac{1}{2}}) (I_{n+1}\otimes A_1)(A_2^t \otimes I_{n+1} ).\nonumber \end{eqnarray}
This implies
\begin{eqnarray*} \label{200427-2,}  \left\{ \begin{array}{ll} 　2  {\tilde{s}}_1 g^{\alpha\beta} \xi_\alpha\xi_\beta =1,\\
  2 {\tilde{s}}_2 + {\tilde{s}}_1 \mu^{\!-\!1} (\mu\!+\!\tilde{\rho})^{\frac{1}{2}}=0,\\
  2{\tilde{s}}_3 + {\tilde{s}}_1  \mu^{\!-\!1} (\mu\!+\!\tilde{\rho})^{\frac{1}{2}}=0,\\
2{\tilde{s}}_4 + {\tilde{s}}_2  \mu^{\!-\!1} (\mu\!+\!\tilde{\rho})^{\frac{1}{2}} + {\tilde{s}}_3  \mu^{\!-\!1} (\mu\!+\!\tilde{\rho})^{\frac{1}{2}}=0,\end{array} \right.\end{eqnarray*}
i.e.,
\begin{eqnarray} \label{200427-3,}  \left\{ \begin{array}{ll} 　  {\tilde{s}}_1=\frac{1}{2 g^{\alpha\beta} \xi_\alpha\xi_\beta},\\
   {\tilde{s}}_2= -\frac{(\mu+\tilde{\rho})^{\!\frac{1}{2}}}{ 4\mu g^{\alpha\beta} \xi_\alpha\xi_\beta},\\
  {\tilde{s}}_3 = -\frac{(\mu+\tilde{\rho})^{\!\frac{1}{2}}}{ 4\mu g^{\alpha\beta} \xi_\alpha\xi_\beta},\\
{\tilde{s}}_4=\frac{(\mu+\tilde{\rho})}{4\mu^2 g^{\alpha\beta} \xi_\alpha\xi_\beta }.\end{array} \right.\end{eqnarray}
It follows that
\begin{eqnarray} \label{19.3.27-4,} && \;\;\quad U^{-1} \!=\! \frac{1}{2  \sqrt{ g^{\alpha\beta} \xi_\alpha \xi_\beta}}\, I_{(n+1)^2}
  -\frac{(\mu+\tilde{\rho})^{\!\frac{1}{2}}}{ 4\mu g^{\alpha\beta} \xi_\alpha\xi_\beta}(I_{n+1}\otimes A_1) -\frac{(\mu+\tilde{\rho})^{\!\frac{1}{2}}}{ 4\mu g^{\alpha\beta} \xi_\alpha\xi_\beta}(A_2^t \otimes I_{n+1} )
  \\
 && \qquad \quad\quad\;\;  + \frac{(\mu+\tilde{\rho})}{4\mu^2 \big(g^{\alpha\beta} \xi_\alpha\xi_\beta\big)^{\frac{3}{2}} }   \!(I_{n+1}\otimes A_1)(A_2^t \otimes I_{n+1} ). \nonumber \end{eqnarray}
 We find from (\ref{19.10.6-1}) that
 \begin{eqnarray*} \label{200427-5,} &&\mbox{vec}\; X= \frac{1}{2  \sqrt{ g^{\alpha\beta} \xi_\alpha \xi_\beta}}(\mbox{vec}\; E)   -\frac{(\mu+\tilde{\rho})^{\!\frac{1}{2}}}{ 4\mu g^{\alpha\beta} \xi_\alpha\xi_\beta}(I_{n+1}\otimes A_1) (\mbox{vec}\; E) \\
&& \qquad \quad \quad   -\frac{(\mu+\tilde{\rho})^{\!\frac{1}{2}}}{ 4\mu g^{\alpha\beta} \xi_\alpha\xi_\beta}(A_2^t \otimes I_{n+1} )(\mbox{vec}\; E) +  \frac{(\mu+\tilde{\rho})}{4\mu^2 \big(g^{\alpha\beta} \xi_\alpha\xi_\beta\big)^{\frac{3}{2}} }   \!(I_{n+1}\otimes A_1)(A_2^t \otimes I_{n+1} )(\mbox{vec}\; E),\nonumber\end{eqnarray*}
so that
\begin{eqnarray} \label{200427-6,} && X= \frac{1}{2  \sqrt{ g^{\alpha\beta} \xi_\alpha \xi_\beta}} E  -\frac{(\mu+\tilde{\rho})^{\!\frac{1}{2}}}{ 4\mu g^{\alpha\beta} \xi_\alpha\xi_\beta} A_1 E \\
&& \qquad \;\, -\frac{(\mu+\tilde{\rho})^{\!\frac{1}{2}}}{ 4\mu g^{\alpha\beta} \xi_\alpha\xi_\beta} E A_2 +  \frac{(\mu+\tilde{\rho})}{4\mu^2 \big(g^{\alpha\beta} \xi_\alpha\xi_\beta\big)^{\frac{3}{2}} }  \! A_1 E A_2.\nonumber\end{eqnarray}
Therefore, when replacing the matrix $E$ in (\ref{200427-6,}) by $E_1$, we immediately get $q_0$, i.e.,
\begin{eqnarray}\label{200427-7,}   q_0=  \frac{1}{2  \sqrt{ g^{\alpha\beta} \xi_\alpha \xi_\beta}} E_1 -\frac{(\mu+\tilde{\rho})^{\!\frac{1}{2}}}{ 4\mu g^{\alpha\beta} \xi_\alpha\xi_\beta}  A_1 E_1 -\frac{(\mu+\tilde{\rho})^{\!\frac{1}{2}}}{ 4\mu g^{\alpha\beta} \xi_\alpha\xi_\beta} E_1 A_2 +    \frac{(\mu+\tilde{\rho})}{4\mu^2 \big(g^{\alpha\beta} \xi_\alpha\xi_\beta\big)^{\frac{3}{2}} }  \! A_1 E_1 A_2.\nonumber\end{eqnarray}

Step 4. \  Furthermore, by considering the terms of degree zero in (\ref{19.3.19-4}), we have
\begin{eqnarray} \label{19.3.28-6}& \big(q_1-b_1\big)  q_{-1} +q_{-1}q_1 =E_0,\end{eqnarray}
where \begin{eqnarray} \label{19.6.1-1}&& E_0:=-q_0^2 +i\sum_{l=1}^{n-1} \big(\frac{\partial q_1}{\partial \xi_l} \frac{\partial q_0}{\partial x_l} +\frac{\partial q_0}{\partial \xi_l}\, \frac{\partial q_1}{\partial x_l}\big) \\
 && \qquad \quad + \frac{1}{2} \sum_{l,\gamma=1}^{n-1} \frac{\partial^2 q_1}{\partial \xi_l\partial \xi_\gamma}   \, \frac{\partial^2 q_1}{\partial x_l \partial x_\gamma}   +b_0 q_0-i\sum_{l=1}^{n-1} \frac{\partial b_1}{\partial \xi_l} \, \frac{\partial q_0}{\partial x_l} +\frac{\partial q_0}{\partial x_n} -c_0.\nonumber\end{eqnarray}
Generally, for $m\ge 1$ we get
\begin{eqnarray} \label{19.3.28-7} \qquad \quad (q_1-b_1)  q_{-m-1} +q_{-m-1}q_1 =E_{-m},\end{eqnarray}
 where \begin{eqnarray} \label{19.6.1-2} E_{-m}:= \begin{small}\sum\limits_{\underset{ |\vartheta| = j+k+m}{-m\le j,k\le 1}} \frac{(-i)^{|\vartheta|}}{\vartheta!}\end{small} ( \partial_{\xi'}^{\vartheta} q_j ) ( \partial_{x'}^{\vartheta} q_k) +b_0 q_{-m}  -i\sum_{l=1}^{n-1} \frac{\partial b_1}{\partial \xi_l}  \frac{\partial q_{-m}}{\partial x_l} +\frac{\partial q_{-m}}{\partial x_n}. \end{eqnarray}
Replacing the matrices $E$ and $X$ by the above $E_{-m}$ and $q_{-m-1}$ in (\ref{200427-6,}), respectively, we explicitly get all $q_{-m-1}$, $m\ge 0$. \qed

\vskip 0.18 true cm

We have obtained the full symbol $q(x, \xi')\sim \sum_{l\le 1} q_l(x,\xi')$ of $Q$ from above proposition 3.1. This implies that modulo a smoothing operator, the pseudodifferential operator $Q$ have been obtained on $\partial \Omega$. Thus we have the following:

\vskip 0.25 true cm

\noindent{\bf Proposition 3.2.} \ {\it In the local boundary normal coordinates, the Dirichlet-to-Neumann map ${\tilde{\Lambda}}_{\tilde{\rho},\mu,g}$ associated with (\ref {20200524-4})   can be represented as:
 \begin{eqnarray}\label{19.3.28-10}&&
 {\tilde{\Lambda}}_{\tilde{\rho},\mu,g} \begin{bmatrix}w_1\\
  \vdots\\
 w_n\\
 f\end{bmatrix}=   Q \begin{bmatrix}w_1\\
  \vdots\\
 w_n\\
 f\end{bmatrix}  \quad  \mbox{on}\;\, \partial \Omega\end{eqnarray}
modulo a  smoothing operator.}

\vskip 0.26 true cm

 \noindent  {\it Proof.} \  Let $(x', x_n)$ be local boundary normal coordinates, for $x_n\in [0,T]$.
 Set $\tilde{\mathcal{S}}:= \frac{\partial^2}{\partial x_n^2} I_{n+1} + B \frac{\partial }{\partial x_n} +C$. Since the principal symbol of $\tilde{{\mathcal{S}}}$  is a negative-definite matrix, the hyperplane $\{x_n=0\}$ is non-characteristic, and hence $\tilde{\mathcal{S}}$ is partially hypoelliptic with respect to this boundary (see p.$\,$107 of \cite{HormL}). Therefore, the solution $(w,f)$ of the equivalent new system of equations $\tilde{\mathcal{S}}(w,f)=0$ is smooth in normal variable, i.e., in boundary normal coordinates $(x',x_n)$ with $x_n\in [0,T]$, $(w,f)\in (C^\infty([0,T]; {\mathfrak{D}}' ({\Bbb R}^{n-1})))^{n+1}$ locally. From Proposition 3.1, we see that the equivalent new system of equations $\tilde{\mathcal{S}}(w,f)=0$  is also locally equivalent to the following system of equations for $(v,h)\in (C^\infty([0,T]; {\mathfrak{D}}' ({\Bbb R}^{n-1})))^{n+1}$: \begin{eqnarray*} &&\bigg(\frac{\partial }{\partial x_n}\, I_{n+1} + Q\bigg) (w,f)=(v,h) \quad \; (w,f)\big|_{x_n=0}=(w^0,f^0),\\
 &&\bigg(\frac{\partial }{\partial x_n} I_{n+1} +B-Q\bigg)(v,h) =(\psi, r)\in (C^\infty ([0,T]\times R^{n-1}))^{n+1}. \end{eqnarray*}
 Making the substitution $t=T-x_n$ for the second equation mentioned above (as done in \cite{LU}), we get a backwards generalized heat equation system:
 \begin{eqnarray*} \bigg(\frac{\partial}{\partial t} I_{n+1}\bigg) (v,h) -(-Q+B)(v,h) =-(\psi,r).\end{eqnarray*}
    Since $(w,f)$ is smooth in the interior of $\Omega$ by interior regularity for the system of elliptic equations $\tilde{\mathcal{S}}(w,f)=0$, it follows that
    $(v,h)$ is smooth in the interior of $\Omega$, and hence
    $(v,h)\big|_{x_n=T}$ is smooth. In view of the principal symbol of $Q$ is strictly positive for any $\xi'\ne 0$, we get that the solution operator for this heat equation system is smooth for $t>0$ (see p.$\,$134 of \cite{Tre}).
    Therefore,
    \begin{eqnarray*} \left(\frac{\partial }{\partial x_n}I_{n+1}\right) (w,f) +Q(w,f)= (v,h)\in (C^\infty ([0,T]\times R^{n-1}))^{n+1} \end{eqnarray*}
    locally.  Setting $J(w^0,f^0)= (v,h)\big|_{\partial \Omega} $, we immediately see that $J$ is a smoothing operator and
    \begin{eqnarray*} \bigg(\big(\frac{\partial}{\partial x_n} I_{n+1} \big)(w,f)\bigg)\bigg|_{\partial \Omega} = -Q(w,f)\big|_{\partial \Omega} +J(w^0,f^0). \end{eqnarray*}
From this,  we obtain (\ref{19.3.28-10}). \qed

 \vskip 0.54 true cm

 \noindent  {\it Proof of Theorem 1.3.} \   Since the Cauchy data ${\mathcal{C}}_{\mu}=\{ (u\big|_{\partial \Omega}, \sigma_\mu (u,p)\nu\big|_{\partial \Omega})\}$  for the Stokes equations is equivalent to the new Cauchy data ${\tilde{\mathcal{C}}}_{\tilde{\rho},\mu}=\{\big( (w,f)\big|_{\partial \Omega}, \frac{\partial (w,f)}{\partial \nu} \big|_{\partial \Omega}\big)\}$ for the new system (\ref{20200524-4}), it suffice to show that the new Dirichlet-to-Neumann map ${\tilde{\Lambda}}_{\tilde{\rho},\mu,g}$ uniquely determines $\mu$ and its derivatives up to order $1$ on $\partial \Omega$.
 Clearly, the Dirichlet-to-Neumann map ${\tilde{\Lambda}}_{\tilde{\rho},\mu,g}$ uniquely determines the principal symbol $q_1$ on $\partial \Omega$.  Now, for given manifold $(\Omega,g)$, we immediately see from (\ref{20200516-1}) in section 3 that $q_1$ uniquely determines $\mu^{-1} (\mu+\tilde{\rho})^{\frac{1}{2}}$ on $\partial \Omega$, and hence $q_1$ uniquely determines $\mu$ on $\partial \Omega$.  Therefore, ${\tilde{\Lambda}}_{\tilde{\rho},\mu,g}$ (only by $q_1$) uniquely determines $\mu$ on $\partial \Omega$.

   Obviously, $q_1$ also uniquely determines  $U^{-1}$ on $\partial \Omega$, where $U^{-1}$ is given in (\ref{19.3.27-4,}) of  section 3. Since ${\tilde{\Lambda}}_{\tilde{\rho},\mu,g}$ uniquely determines $q_0$, it follows that ${\tilde{\Lambda}}_{\tilde{\rho},\mu,g}$ uniquely determines $\mbox{vec}\, (q_0)$, so that ${\tilde{\Lambda}}_{\tilde{\rho},\mu,g}$ uniquely determines $\mbox{vec}\,( E_1)$  by $q_0= U^{-1} E_1$. In other words, $q_1$ and $q_0$ uniquely determine $E_1$.  Note that
$E_1=b_0 q_1 + \frac{\partial q_1}{\partial x_n} -c_0 +M_0$, where $M_0$ is a matrix expression involving only $\mu$ and its all tangential derivatives of order $1$ along $\partial \Omega$. It follows that
\begin{align*}&  b_0q_1 + \frac{\partial q_1}{\partial x_n} -c_1 = \sqrt{g^{\alpha\beta} \xi_\alpha\xi_\beta} \,\, \Gamma_{n\gamma}^\gamma \,I_{n+1} \\
&    +\!
 \sqrt{g^{\alpha\beta} \xi_\alpha\xi_\beta}  \begin{bmatrix}
\begin{BMAT}(@, 15pt, 15pt){c.c}{c.c}
  \! \left[\delta_{kn} \mu (\mu^{\!-\!1})^{;j}\!+\! \frac{\delta_{jk} \tilde{\rho}}{\mu(\mu\!+\!\tilde{\rho})}\frac{\partial \mu}{\partial x_n}
  \!+\! 2\Gamma_{kn}^j  \right]_{\!n\times \!n}   &\left[\!\frac{(\mu\!+\!\tilde{\rho})^{\frac{1}{2}}}{\mu}\!\!\left\{\!\!\!\begin{array} {ll} 2 R^j_n  - 2\mu  \frac{\partial(\mu^{-1})^{;j}}{\partial x_n} - 2\mu \Gamma^j_{\alpha n} (\mu^{-1})^{;\alpha} \\
\!+  ( g^{\alpha\beta} \frac{\partial \Gamma_{\!n\alpha}^j}{\partial x_\beta} \!-\! g^{\alpha\gamma} \Gamma_{\!\beta\gamma}^j \Gamma_{\!n\alpha}^\beta \! \!-\! g^{\alpha \beta} \Gamma_{\!n\gamma}^j \Gamma_{\!\alpha\beta}^\gamma ) \\
 - 2 g^{\alpha \gamma} g^{\beta \sigma} \Gamma_{\gamma \sigma}^j \Gamma_{\alpha\beta}^n
  \end{array}\!\!\!\! \right\}\!\right]_{\!n\!\times \!1}
        \\
     \left[\delta_{nk} \,\mu(\mu+\tilde{\rho})^{-\frac{1}{2}}\right]_{1\!\times \!n}  &    0 \! \end{BMAT}\!
\end{bmatrix}\\
&     +\Big(\mu^{\!-1}(\mu+\tilde{\rho})^{\frac{1}{2}}\Gamma_{n\gamma}^\gamma
+\frac{\partial (\mu^{-1} (\mu+\tilde{\rho})^{\frac{1}{2}})}{\partial x_n} \Big)    \begin{small} \begin{bmatrix}
\begin{BMAT}(@, 15pt, 15pt){c.c}{c.c}
   \!\left[0\right]_{n\times n}   &   \! \left[ 2\text{i}\,  \Gamma_{\beta n}^j g^{\beta \alpha} \xi_\alpha\!+\!
   \frac{1}{\sqrt{ g^{\alpha\beta}\xi_\alpha\xi_\beta }}\Gamma_{\gamma \sigma}^j g^{\alpha \gamma} g^{\beta \sigma} \xi_\alpha\xi_\beta\right]_{\!n\times 1}\\
   \left[0 \right]_{1\times n}&  0  \end{BMAT}
\end{bmatrix} \end{small}\\
&       +  \begin{bmatrix}
\begin{BMAT}(@, 15pt, 15pt){c.c}{c.c}
  \! \left[0  \right]_{n\times n}   \!& \left[ \left(\delta_{ln} \mu(\mu^{-1})^{;j}+\frac{\delta_{jl}\tilde{\rho}}{\mu(\mu+\tilde{\rho})}\frac{\partial \mu}{\partial x_n}
 \!+\! 2\Gamma^j_{ln}
  \right)  \Big(2\text{i} g^{\alpha\beta} \Gamma_{\beta n}^l \xi_\alpha +
   \frac{1}{ \sqrt{ g^{\alpha \beta} \xi_\alpha\xi_\beta} } \Gamma_{\gamma \sigma}^l g^{\alpha \gamma}g^{ \beta\sigma} \xi_\alpha\xi_\beta \Big)\right]_{n\times 1} \\
\left[0\right]_{1\times n} & \delta_{nl} \,\mu(\mu+\tilde{\rho})^{-\frac{1}{2}} \Big(2\text{i} \, g^{\alpha\beta} \Gamma_{\beta n}^l \xi_\alpha +
   \frac{1}{ \sqrt{ g^{\alpha \beta} \xi_\alpha\xi_\beta} } \Gamma_{\gamma \sigma}^l g^{\alpha \gamma}g^{ \beta\sigma} \xi_\alpha\xi_\beta \Big)
  \end{BMAT}
\end{bmatrix}\\
&    + \frac{\frac{\partial g^{\alpha\beta}}{\partial x_n} \xi_\alpha\xi_\beta} {2\sqrt{g^{\alpha\beta} \xi_\alpha\xi_\beta}} I_{n+1} + \mu^{-1}(\mu+\tilde{\rho})^{\!\frac{1}{2}}  \frac{\partial }{\partial x_n} \! \begin{small} \begin{bmatrix}
\begin{BMAT}(@, 15pt, 15pt){c.c}{c.c}
   \!\left[0\right]_{n\times n}   &   \! \left[ 2\text{i}\, g^{\beta \alpha} \Gamma_{\beta n}^j \xi_\alpha\!+\!
   \frac{1}{\sqrt{ g^{\alpha\beta}\xi_\alpha\xi_\beta }}\Gamma_{\gamma \sigma}^j g^{\alpha \gamma} g^{\beta \sigma} \xi_\alpha\xi_\beta\right]_{\!n\times 1}\\
   \left[0 \right]_{1\times n}&  0  \end{BMAT}
\end{bmatrix} \end{small}\\
& -  \text{i} \, \big( g^{\alpha\beta} \Gamma_{\alpha\gamma}^\gamma +\frac{\partial g^{\alpha\beta}}{\partial x_\alpha} \big) \xi_\beta \,I_{n+1}\qquad \qquad \qquad \qquad \qquad \qquad \qquad \end{align*}
  \begin{align*}\!&  \!-\! \begin{bmatrix}
\begin{BMAT}(@, 5pt, 5pt){c.c}{c.c}
  \! \!\left[\!\text{i}(\!1\!-\!\delta_{nk}\!)\mu(\mu^{\!-\!1}\!)^{;j} \!\xi_k \!+\!
  \text{i}\big(\!\frac{\delta_{jk} \tilde{\rho}}{\mu(\mu\!+\!\tilde{\rho})}\frac{\partial \mu}{\partial x_\alpha}
  \!+\! 2  \Gamma^j_{\!k\alpha}\!\big)g^{\alpha\beta} \xi_\beta  \! \right]_{\!n\!\times \!n}   \! &  \! \left[\!\frac{(\mu\!+\!\tilde{\rho})^{\frac{1}{2}} \text{i}}{\mu}\!\left\{\!\!\!\begin{array} {ll} \!2 R^j_\alpha g^{\alpha\beta}  \!\!-\! 2\mu g^{\alpha\beta}  \frac{\partial((\mu^{\!-1}\!)^{;j})}{\partial x_\alpha}\\
  \! -\! 2\mu g^{\alpha\beta} \Gamma^j_{s\alpha } \!(\mu^{\!-1})^{;s}  \!
 -\! 2\Gamma_{sh}^j   g^{sr}\! g^{hm} \Gamma_{rm}^\beta\\
  \!-\!( g^{mr}\! \frac{\partial \Gamma_{\!\alpha r}^j}{\partial x_m}
  \!-\!\! g^{mr} \Gamma_{\!hr}^j \Gamma_{\!\alpha m}^h \!
 \!-\!\! g^{mr} \!\Gamma_{\!\alpha h}^j\! \Gamma_{\!mr}^h )g^{\alpha\!\beta} \\
 \end{array}\!\!\!\!\! \right\}\!\xi_\beta \! \right]_{\!n\!\times \!1}\!\\
     \left[\text{i} (1-\delta_{nk})\mu(\mu+\tilde{\rho})^{-\frac{1}{2}}  \xi_k\right]_{1\times n} & 0\!\end{BMAT}\!
\end{bmatrix}.
\end{align*}
Obviously,  the $(n,n)$ entry of the matrix on the right-hand side of the above equality is
 \begin{eqnarray}   \label{200428-1}
&& \sqrt{g^{\alpha\beta} \xi_\alpha\xi_\beta} \; \Gamma_{n\beta}^\beta +  \sqrt{g^{\alpha\beta} \xi_\alpha\xi_\beta}\,\Big( \mu (\mu^{\!-1})^{;n} +\frac{\tilde{\rho}}{\mu(\mu+\tilde{\rho})}\frac{\partial \mu}{\partial x_n}\Big) +\frac{\frac{\partial g^{\alpha\beta}}{\partial x_n} \xi_\alpha\xi_\beta}{2\sqrt{g^{\alpha\beta} \xi_\alpha\xi_\beta}}  \\
&&\qquad\quad \qquad \quad \; \;- \text{i} \,( g^{\alpha\beta} \Gamma_{\alpha\gamma}^\gamma +\frac{\partial g_{\alpha\beta}}{\partial x_\alpha}) \xi_\beta - \frac{\text{i}  \,\tilde{\rho}}{\mu(\mu+\tilde{\rho})}\frac{\partial \mu}{\partial x_\alpha} g^{\alpha\beta} \xi_\beta \nonumber \\
&&\qquad\quad \qquad \; \,= \sqrt{g^{\alpha\beta} \xi_\alpha\xi_\beta}\,\Big( \mu (\mu^{\!-1})^{;n} +\frac{\tilde{\rho}}{\mu(\mu+\tilde{\rho})}\frac{\partial \mu}{\partial x_n}\Big) + M'_0 \nonumber\\
&&\qquad\quad  \qquad\; \,=  \sqrt{g^{\alpha\beta} \xi_\alpha\xi_\beta}\,\Big( - \mu^{-1} \frac{\partial \mu} {\partial x_n}  +\frac{\tilde{\rho}}{\mu(\mu+\tilde{\rho})}\frac{\partial \mu}{\partial x_n}\Big) + M'_0\nonumber \\
&&\qquad\quad  \qquad\; \,= -\frac{1}{\mu+\tilde{\rho}}\sqrt{g^{\alpha\beta} \xi_\alpha\xi_\beta}\,  \frac{\partial \mu} {\partial x_n}  + M'_0, \nonumber\end{eqnarray}
where $M'_0$ is a matrix expression involving only $\mu$ and its all tangential derivatives of order $1$ along $\partial \Omega$.
Since $-\frac{1}{\mu+\tilde{\rho}} \sqrt{g^{\alpha\beta} \xi_\alpha\xi_\beta}$ has been determined by $q_1$ on $\partial \Omega$,
  we immediately find by (\ref{200428-1}) that $q_0$ determines the $\frac{\partial \mu}{\partial x_n}$ along on $\partial \Omega$. Note that $\frac{\partial \mu}{\partial x_\alpha}$, ($1\le  \alpha \le n-1$), can be determined by $\mu$ along on $\partial \Omega$. Thus $\frac{\partial^{|K|} \mu}{\partial x_K}$, $\,\forall x\in \partial \Omega$, $|K|\le 1$ can uniquely be  determined  by  $q_1$ and $q_0$. Of course, $\mu$ and $\nabla_g \mu $ can be uniquely determined by $\Lambda_{\mu,g}$ on $\partial \Omega$. \qed

 \vskip 0.29 true cm

\noindent{\bf Remark 3.3.} \  {\it Similar to \cite{Liu1} for the elastic Dirichlet-to-Neumann map  (or \cite{Liu2} for the electromagnetic Dirichlet-to-Neumann map), we can further show that ${\tilde{\Lambda}}_{\tilde{\rho}, \mu, g}$  uniquely determine $\mu$ and $\frac{\partial^{|K|} \mu}{\partial x_K}$  for all multi-index $K$ with
$|K|\ge 0$ on $\partial \Omega$. This kinds of results are called  Kohn-Vogelius type theorem because R.  Kohn  and M. Vogelius first showed that Dirichlet-to-Neumann map associated with the equation $-\mbox{div}\; (a\nabla v)=0$ uniquely determines $a$ and its all order derivatives on $\partial \Omega$ in \cite{KV}. Our method is based on theory of pseudodifferential operators and some exact calculations (see, \cite{Liu1} and \cite{Liu2} for more details).}

\vskip 1.49 true cm

\section{ Determining the viscosity  for the stationary Navier-Stokes equations}

\vskip 0.45 true cm

In this section we consider the unique determination of the viscosity in an incompressible fluid described by the stationary Navier-Stokes equations. Under some additional assumptions, this problem has been solved by Li, Uhlmann and Wang \cite{LUW} in two dimensions  and  by Li and Wang \cite{LiW} in three dimensions using the linearization technique. (The linearization technique was first introduced by Isakov in \cite{Isak1}, which allows for the reduction of the semilinear inverse boundary problem to the corresponding linear one, see also \cite{HeS}, \cite{IsS}, \cite{Sun} or \cite{SunU}.)  We will apply their techniques in \cite{LiW} and \cite{LUW} and our method to show the uniqueness result of the viscosity $\mu$ for the Navier-Stokes equations in any bounded domain $\Omega\subset {\mathbb{R}}^n$ for any dimensional case because our method and their technique are both independent of spatial dimensions. Let $u = (u_1,\cdots, u_n)^t$ be the velocity vector field satisfying the stationary Navier-Stokes equations
\begin{eqnarray} \label{20200526-1} \left\{ \begin{array}{ll} \mbox{div}\; \sigma_\mu(u,p)-(u\cdot \nabla)u=0 \;\;\, &\mbox{in}\;\; \Omega,\\
\mbox{div}\; u=0 \;\;\, &\mbox{in}\;\; \Omega,\end{array}\right. \end{eqnarray}
and the corresponding Cauchy data is denoted by  \begin{eqnarray*} \label{20200526-2}  {\mathcal{CN}}_{\mu} = \Big\{(u,\sigma_\mu (u,p)\nu)\big|_{\partial \Omega} \big| (u,p) \;\, \mbox{satisfies}\;\,  (\ref{20200526-1})\;\, \mbox{with}\;\, \int_{\Omega} p\, dV=0\Big\}.\end{eqnarray*}
Let $u\big|_{\partial \Omega} = \phi \in H^{\frac{3}{2}} (\partial \Omega)$ satisfy (\ref{20200502-2}). In order to study inverse problem for Navier-Stokes equations,  in \cite{LiW} Li and Wang took $\phi= \epsilon \varphi \in H^{\frac{3}{2}} (\partial \Omega)$ with $|\epsilon|$ sufficiently small and let $(u_\epsilon, p_\epsilon) = (\epsilon v_\epsilon, \epsilon r_\epsilon)$ satisfy (\ref{20200526-1}). The problem (\ref{20200526-1}) is reduced to
\begin{eqnarray}\label{20200526-3} \left\{ \begin{array}{ll} \mbox{div}\;\sigma_\mu (v_\epsilon,r_\epsilon) -\epsilon (v_\epsilon \cdot \nabla)v_\epsilon =0\;\; &\mbox{in}\;\; \Omega,\\
   \mbox{div}\; v_\epsilon  = 0\;\; &\mbox{in}\;\; \Omega,\\
    v_\epsilon = \varphi \;\; &\mbox{on}\;\; \partial \Omega.\end{array} \right. \end{eqnarray}
    Look for a solution of (\ref{20200526-3}) with the form $v_\epsilon = v_0 + \epsilon v$ and $r_\epsilon = r_0 + \epsilon r$, where $(v_0,r_0)$ satisfies the Stokes equations
    \begin{eqnarray} \label{20200526} \left\{ \begin{array} {ll} \mbox{div}\;\sigma_{\mu}(v_0,r_0)=0\;\;&\mbox{in}\;\;\Omega, \\
    \mbox{div}\; v_0 =0 \;\;&\mbox{in}\,\; \Omega, \\
    v_0 = \varphi \;\; &\mbox{on}\;\; \partial \Omega,\end{array} \right. \end{eqnarray}
     and $(v,r)$ satisfies
     \begin{eqnarray} \label{20200526-4} \left\{ \begin{array}{ll}    -\mbox{div}\; \sigma_{\mu} (v,r) + \epsilon (v_0 \cdot \nabla)v + \epsilon (v \cdot \nabla )v_0 + \epsilon^2 (v\cdot \nabla)v = h \;\;&\mbox{in}\;\;  \Omega,\\
      \mbox{div}\; v = 0 \;\;&\mbox{in}\;\; \Omega,\\
       v=0  \;\;&\mbox{on}\;\; \partial\Omega\end{array} \right. \end{eqnarray}  with $h=-(v_0 \cdot \nabla)v_0$. In \cite{LiW}, it is shown that for any $\varphi \in H^{\frac{3}{2}} (\partial \Omega)$,  let $(v_0,r_0)\in H^2 (\Omega) \times H^1 (\Omega)$ be the unique solution ($r_0$ is unique up to a constant) of the Stokes equations (\ref{20200526}). There exists a solution $(u_\epsilon,p_\epsilon)$ of (\ref{20200526-1}) of the form $u_\epsilon = \epsilon v_0 + \epsilon^2 v$, $p_\epsilon = \epsilon r_0 + \epsilon^2 r$
with the boundary data $u_\epsilon\big|_{\partial \Omega} = \epsilon \varphi$ for all $|\epsilon|\le \epsilon_0$, where $\epsilon_0$ depends on $\|\varphi\|_{H^{\frac{3}{2}} (\partial \Omega)}$. Moreover, it is proved  in  \cite{LiW} that as $\epsilon\to 0$,
 \begin{eqnarray} \label{20200526-6} \left.\begin{array}{ll}
 \big\| \epsilon^{-1} u_\epsilon\big|_{\partial \Omega} - v_0\big|_{\partial \Omega} \big\|_{H^{\frac{3}{2}} (\partial \Omega)} \to 0,\\
     \big\| \epsilon^{-1} \sigma_\mu(u_\epsilon,p_\epsilon)\nu\big|_{\partial \Omega} -\sigma_\mu (v_0,r_0)\nu\big|_{\partial \Omega}\big\|_{H^{\frac{1}{2}}(\partial \Omega)}\to 0,\end{array} \right.\end{eqnarray}
          provided $$\int_{\Omega}  p_\epsilon \,dx  =\int_\Omega r_0 \,dx = 0.$$  This implies  that the Cauchy data ${\mathcal{NC}}_{\mu}$  of the Navier-Stokes equations uniquely determines the Cauchy data ${\mathcal{C}}_\mu$ of the Stokes equations. In other words, ${\mathcal{NC}}_{\mu_1}  ={\mathcal{NC}}_{\mu_1}$ implies ${\mathcal{C}}_{\mu_1} = {\mathcal{C}}_{\mu_1}$. Therefore, the uniqueness of the viscosity for the Navier-Stokes equations follows from our Theorem 1.4. We have the following theorem:

        \vskip 0.29 true cm

\noindent{\bf Theorem 4.1.} \  {\it   Let $\Omega$ be a simply connected bounded domain in ${\mathbb{R}}^n$, ($n=2,3$), with smooth boundary. Suppose that $\mu_1$  and $\mu_2$ are two viscosity functions for the Navier-Stokes equations. Assume that $\mu_j>$ in $\bar \Omega$ and $\mu_j\in C^3 (\bar \Omega)$ for $n=2$ and $\mu_j\in C^8(\bar \Omega)$ for $n=3$. Let ${\mathcal{NC}}_{\mu_1}$ and ${\mathcal{NC}}_{\mu_2}$ be the Cauchy data associated with $\mu_1$  and $\mu_2$, respectively. If ${\mathcal{NC}}_{\mu_1}  ={\mathcal{NC}}_{\mu_1}$ , then $\mu_1 = \mu_2$  in $\Omega$.}

\vskip 1.68 true cm

\centerline {\bf  Acknowledgments}

\vskip 0.59 true cm
   This research was supported by NNSF of China (11671033/A010802).

  \vskip 1.68 true cm

\begin{center}

\end{center}

\end{document}